\documentclass[preprint,sort&compress,12pt]{elsarticle}
\usepackage{bbm}
\usepackage{amsmath}
\usepackage{mathrsfs}
\usepackage{stmaryrd}
\usepackage{mathrsfs}
\usepackage{bm}
\usepackage{amsmath}
\usepackage{amsfonts}
\usepackage{color}
\usepackage{amssymb}
\usepackage{amsthm}
\usepackage{lineno}
 \usepackage{leftidx}

\allowdisplaybreaks
\providecommand{\abs}[1]{\left\vert#1\right\vert}

\newcommand{\mf}{\mathbf}
\newcommand{\mm}{\mathrm}
\newcommand{\ml}{\mathcal}

\newcommand{\be}{\begin{equation}}
\newcommand{\bea}{\begin{equation}\begin{aligned}}
\newcommand{\beas}{\begin{equation*}\begin{aligned}}
\newcommand{\eeas}{\end{aligned}\end{equation*}}
\newcommand{\eea}{\end{aligned}\end{equation}}
\newcommand{\ee}{\end{equation}}

\begin{document}
\begin{frontmatter}
\title{Rayleigh--Taylor Instability  in  Stratified Compressible Fluids with/without the Interfacial Surface Tension}
\author[FJ]{Fei Jiang
}
\ead{jiangfei0591@163.com}
\author[FJ]{Han Jiang}
\ead{225410005@fzu.edu.cn}
\author[sJ]{Song Jiang}
\ead{jiang@iapcm.ac.cn}
\address[FJ]{School of Mathematics and Statistics, Fuzhou University, Fuzhou, 350108, China.}
\address[sJ]{Institute of Applied Physics and Computational Mathematics, 
 Beijing, 100088, China.}
\begin{abstract}
  Guo--Tice formally  established in 2011 that the Rayleigh--Taylor instability inevitably occurs
  within stratified compressible viscous fluids in a slab domain $\mathbb{R}^2\times (h_-,h_+)$,
  irrespecive of the presence of interfacial surface tension,
   where the instability solutions are non-periodic with respect to both horizontal spacial variables $x_1$ and $x_2$,
   by applying a so-called ``normal mode'' method and a modified variational method to the linearized (motion) equations \cite{GYTI2}.
  It is a long-standing open problem, however, whether Guo--Tice's conclusion can be rigorously verified by the (original)
  nonlinear equations.  This challenge arises due to the failure of constructing a growing mode solution,
  which is non-periodic with respect to both horizontal spacial variables, to the linearized equations defined on a slab domain.
  In the present work, we circumvent the difficulty related to growing mode solutions by developing an alternative approximate scheme.
  In essence, our approach hinges on constructing the horizontally periodic growing mode solution of the linearized equations
  to approximate the {\it nonlinear} Rayleigh--Taylor instability solutions, which do not exhibit horizontal periodicity.
     Thanks to this new approximate scheme, we can apply Guo--Hallstrom--Spirn's bootstrap instability method in \cite{GYHCSDDC} to the nonlinear equations in Lagrangian
      coordinates, and thus prove  Guo--Tice's conclusion.  In particular,
     our approximate method could also be applied to other instability solutions characterized by non-periodic motion in a slab domain,
     such as the Parker instability and thermal instability.
 \end{abstract}
\end{frontmatter}

\newtheorem{thm}{Theorem}[section]
\newtheorem{lem}{Lemma}[section]
\newtheorem{pro}{Proposition}[section]
\newtheorem{cor}{Corollary}[section]
\newproof{pf}{Proof}
\newdefinition{rem}{Remark}[section]
\newtheorem{definition}{Definition}[section]
\section{Introduction}
\label{Intro} \numberwithin{equation}{section}

Consider two  perfectly  plane-parallel layers of immiscible fluids, the heavier on top of the lighter one,
and both subject to the earth's gravity. In this case, the equilibrium state is unstable to sustain a small disturbance,
and this unstable disturbance will grow and lead to a release of potential energy, as the heavier fluid moves down under the gravitational
force, and the lighter one is displaced upwards. This phenomenon was first studied by Rayleigh \cite{RLIS}
and  subsequently by Taylor \cite{TGTP}, and is called therefore the Rayleigh--Taylor (RT) instability.
In the last decades, this phenomenon has been extensively investigated from both physical and numerical  perspectives, see \cite{CSHHSCPO,WJH,GYTI1,desjardins2006nonlinear,GBKJSAN} for examples.
 It has been also widely investigated how the RT instability evolves under  the influence of various physical factors,
 such as elasticity \cite{JFJWGCOSdd,FJWGCZXOE}, rotation \cite{CSHHSCPO,BKASMMHRJA}, surface tension \cite{GYTI2,WYJTIKCT,JJTIWYJ,wilke2022rayleigh}, and magnetic fields \cite{JFJSWWWOA,JFJSJMFM,WYJ2019ARMA,JFJSARMA2019,JFJSNS}.

  In this article, we are interested in the existence of RT instability solutions within  stratified compressible viscous
  fluids with/without interfacial surface tension. Before introducing our main result and the relevant progress in Section \ref{202306071229},
  first we shall formulate the RT instability problem.

\subsection{ Model  for stratified compressible viscous fluids}
  The three-dimensional motion equations of compressible fluids
 under a uniform gravitational field (directed along the negative $x_3$-axis)  can be described as follows:
\begin{equation}\label{0101fMnew}
  \begin{cases}
  \rho_t+\mm{div}(\rho{  v})=0 ,\\[1mm]\rho v_t+\rho v\cdot\nabla v+\mm{div}\mathcal{S}=-\rho g \mathbf{e}^3.
  \end{cases}\end{equation}
Here $\rho:=\rho(x,t)$, $v:=v(x,t) $,  $\mathcal{S} $  and  $g$ denote the density, velocity,
  stress tensor and gravitational constant, respectively. The vector $\mathbf{e}^3:=(0,0,1)^{\top }$, where the superscript $\top$ denotes
  the transposition.
  In the above system we consider that the stress tensor $\mathcal{S}$  is characterized by the expression:
\begin{equation}
\label{Cauchy2219}
\mathcal{S}:= P\mathbb{I}- \mathbb{S}(v) ,
 \end{equation}
where  $\mathbb{I}$  denotes  a $3\times 3$ identity matrix. $P:= P(\tau)|_{\tau=\rho}$ and  $\mathbb{S}(v)$ denote
the hydrodynamic pressure and viscosity tensor, respectively. In this article, \emph{the pressure function $P(\tau)$ is always assumed to be
smooth, positive, and strictly increasing with respect to $\tau$}, and the viscosity tensor is given by
\begin{equation}
\mathbb{S}(v) := \mu  \mathbb{D}v+\left(\varsigma-{2\mu}/{3}\right)\mm{div} v\, \mathbb{I},\end{equation}
 where $\mathbb{D}v:= \nabla v+ \nabla v^\top  $, and $\mu>0$ resp. $\varsigma\geqslant  0$ denote the shear resp. bulk viscosity coefficients.
  For the sake of simplicity, in this paper we restrict our consideration to the case where $\mu$ and $\varsigma$ are constants.

To investigate the  RT instability, we shall further consider two distinct, immiscible, compressible viscous fluids
evolving in a moving domain $\Omega(t)$, where $\Omega(t):=\Omega_+(t)\cup\Omega_-(t)$ for $t\geqslant 0$. The upper fluid
 occupies  the upper domain:
$$  \Omega_+(t):=\{(x_{\mm{h}},x_3)^{\top }~|~x_{\mm{h}}:=(x_1,x_2)^{\top }\in \mathbb{R}^2 ,\ d(x_{\mm{h}},t)< x_3<h_+\},
$$
while the lower fluid occupies the lower domain:
$$
\Omega_-(t):=\{(x_{\mm{h}},x_3)^{\top }~|~x_{\mm{h}}\in \mathbb{R}^2 ,\ h_-< x_3<d(x_{\mm{h}},t)\}.
$$
Here  $h_-$ and $h_+$ are fixed constants that satisfy $h_-<h_+$, whereas the internal surface function
$d:=d(x_{\mm{h}},t)$ is unspecified and unknown.  The internal surface $\Sigma(t):=\{x_3=d\}$ moves between the two fluids,
 and $\Sigma_\pm:=\{x_3 =h_\pm\}$ represent the fixed upper and lower boundaries of $\Omega(t)$, respectively.

 In the following, we employ  the equations \eqref{0101fMnew} to describe the motion of stratified compressible viscous fluids.
We introduce the subscript  $+$ and $_-$  to the notations of the known physical parameters, pressure functions and other
unknown functions in \eqref{0101fMnew} for  the upper and lower fluids, respectively. Consequently,
the equations governing the motion of the stratified compressible viscous fluids in a slab domain, subject to
a uniform gravitational field, can be formulated as follows:
 \begin{equation}\label{0101f1} \begin{cases}
\partial_t \rho_\pm+\mm{div}(\rho_\pm{  v}_\pm)=0& \mbox{ in } \Omega_\pm(t), \quad t>0,\\
\rho_\pm \partial_t v_\pm+\rho_\pm v_\pm\cdot\nabla v_\pm+\mm{div}\mathcal{S}_\pm  =-\rho_\pm g\mathbf{e}^3&\mbox{ in }
\Omega_\pm(t),\quad t>0,\end{cases}
\end{equation}
where $\mathcal{S}_\pm$ are defined by \eqref{Cauchy2219} with $(\mu_\pm,\varsigma_\pm,v_\pm,  P_\pm)$ in place of $(\mu,\varsigma,v,  P)$.

For two viscous fluids meeting at a free boundary, the  conventional  assumptions are that the velocity is continuous
across the interface and that the jump in the normal stress is proportional to the mean curvature of the free surface multiplied
by the normal to the surface \cite{XLZPZZFGAR}. This requires us to enforce the following (interfacial) jump conditions
\begin{equation}\label{201612262131}
     v_+|_{\Sigma(t)}- v_-|_{\Sigma(t)}   =0\mbox{ and }
(  \mathcal{S}  _+\nu |_{\Sigma(t)} -  \mathcal{S}_- \nu |_{\Sigma(t)} ) = \vartheta \mathcal{C} \nu \mbox{ on }\Sigma(t).
\end{equation}
 In these conditions, the coefficient of the interfacial surface tension is denoted  by $\vartheta\geqslant 0$,
the normal vector to $\Sigma(t)$ by $\nu$, and  twice the mean curvature  of the internal surface $\Sigma(t)$
by $\mathcal{C}$ \cite{GYTI2}, i.e.,
$$
\mathcal{C}:=\frac{\Delta_{\mm{h}}d+(\partial_1 d)^2\partial_2^2
d+(\partial_2d)^2\partial_1^2d-2\partial_1d\partial_2d\partial_1\partial_2d}
{(1+(\partial_1d)^2+(\partial_2d)^2)^{3/2}}.
$$
Additionally, we also enforce the condition that the fluid velocity vanishes at the fixed boundaries,
 implemented through the boundary conditions:
\begin{equation}\label{201612262117}
v_\pm=\mathbf{0}  \mbox{ on }\Sigma_\pm.
\end{equation}

Moreover,  under the first jump condition in \eqref{201612262131}, the internal surface function is determined
by $v_+$ (or $v_-$), i.e., for each $t>0$,
\begin{equation}\label{201612262217}
d_t+v_{1,+}(x_{\mm{h}},d) \partial_1d+v_{2,+}(x_{\mm{h}},d) \partial_2d=v_{3,+}(x_{\mm{h}},d) \mbox{ in }\mathbb{R}^2,
\end{equation}
where $v_{i,+}$ denotes the $i$-th component of $v_+$ for $1 \leqslant i\leqslant 3$.
Finally, we impose the initial data of $(\rho, v, d)$:
\begin{equation}
\label{201612262216}
(\rho, v )|_{t=0}:=( {\rho}^0,v^0 ) \mbox{ in } \Omega\!\!\!\!\!-\setminus \Sigma(0)\mbox{ and }d|_{t=0}=d^0  \mbox{ on } \mathbb{R}^2,
\end{equation}
 where we denote  $\Omega\!\!\!\!-:=\mathbb{R}^2\times \{h_-,h_+\}$ and $\Sigma(0):=\{x_3=d(x_\mm{h},0)\}$. Then, the system  \eqref{0101f1}--\eqref{201612262216} constitutes an initial-boundary value problem for stratified compressible
 viscous fluids with/without internal surface tension. For the sake of simplicity,  we refer to the problem defined by
 \eqref{0101f1}--\eqref{201612262216} as the stratified compressible viscous fluid (SCVF) model.

The above SCVF model has been used to investigate the evolution of compressible RT instability. To achieve this objective,
it is necessary to further construct a RT equilibrium of the SCVF model. For this purpose, we choose
 a constant $\bar{d}\in (h_-,h_+)$, and density profiles $\bar{\rho}_\pm$, which are \emph{differentiable functions in $\Omega_\pm$},
 dependent solely on $x_3$, and fulfill the hydrostatic relations:
 \begin{equation}
\label{201611051547}
\begin{cases}
\nabla P_\pm(\bar{\rho}_\pm)  =-\bar{\rho}_\pm g \mathbf{e}^3 &\mbox{ in } \Omega_\pm,\\[1mm]
 \llbracket  P(\bar{\rho})  \rrbracket  \mathbf{e}^3=0  &\mbox{ on }\Sigma,
  \end{cases}
\end{equation}
the non-vacuum condition
\begin{equation}
\inf_{x\in \Omega_\pm}\{\bar{\rho}_\pm(x_3)\}>0,
\end{equation}
and the RT (jump) condition
 \begin{equation}
\label{201612291257} \llbracket   \bar{\rho}  \rrbracket >0\mbox{ on }\Sigma,\end{equation}
where we have denoted
\begin{align}
&\Omega_+:=\mathbb{R}^2\times\{\bar{d}<x_3<h_+\},\
\Omega_-:=\mathbb{R}^2\times\{h_-<x_3<\bar{d}\}, \ \Sigma:=\mathbb{R}^2\times\{\bar{d}\},\nonumber \\
& \llbracket  P(\bar{\rho})  \rrbracket :=  P_+(\bar{\rho}_+)|_{\Sigma}  -P_-(\bar{\rho}_-) |_{\Sigma}
\;\;\mbox{ and }\;\;\llbracket   \bar{\rho}  \rrbracket :=   \bar{\rho}_+ |_{\Sigma}  - \bar{\rho}_- |_{\Sigma} . \label{202307302006}
\end{align}

Now, let us further define
$$\bar{\rho} := \bar{\rho}_+ \mbox{ for }x\in {\Omega_+ }\;\mbox{ and }\; \bar{\rho} _- \mbox{ for }x\in {\Omega_- }. $$
 This leads to an RT equilibrium solution $(\rho,v )=(\bar{\rho}, \mathbf{0})$ with $d=\bar{d}$ of the SCVF model.
We mention that such an equilibrium solution $(\bar{\rho},\mathbf{0} )$, satisfying \eqref{201611051547}--\eqref{201612291257},
indeed exists, see \cite{GYTI1}. \emph{In addition, we assume without loss of generality that $\bar{d}=0$ in this article.}
If $\bar{d}$ is non-zero, we can adjust the $x_3$ coordinate to make $\bar{d}=0$. Consequently,
 $h_-<0$,  and $d$ can be referred to as the displacement function, representing the deviation of the interface point
  from the plane $\Sigma$. It has been established since 1953 that the presence of interfacial surface tension $\vartheta \mathcal{C}\nu$
 can slow down the growth rate of RT instability \cite{1953Effects}.
 In particular, the question of whether  interfacial surface tension inhibits the RT instability
 (i.e., the interfacial surface tension prevents the heavier fluid from  sinking into  the lower fluid)
 reduces to  the verification of whether the RT equilibrium solution to the SCVF model is stable.
 The aim of this article is to show the instability of the RT equilibrium solution to the SCVF model
   in the context of non-periodic motion.

\subsection{Reformulation of the SCVF model}\label{2022307272010}
Since the movement of the free interface $\Sigma(t)$ and the subsequent change of the
domains $\Omega_\pm(t)$ in Eulerian coordinates will result in several  mathematical difficulties,
we shall switch our analysis to Lagrangian coordinates, so that the interface and the domains stay fixed in time. To this end, we take $\Omega_+$ and $\Omega_-$ to be the fixed Lagrangian domains, and assume that there exist
invertible mappings
\begin{equation*}\label{0113}
\zeta_\pm^0:\Omega_\pm\rightarrow \Omega_\pm(0),
\end{equation*}
such that $\det \nabla\zeta_\pm^0\neq 0$, and
\begin{equation}
\label{05261240}
\Sigma(0)=\zeta_\pm^0(\Sigma),\
\Sigma_+=\zeta_+^0(\Sigma_+)\mbox{ and }
\Sigma_-=\zeta_-^0(\Sigma_-),
\end{equation}
 where $\mm{det}$ denotes the  determinant operator.
 The first condition in \eqref{05261240} means that the initial interface
$\Sigma(0)$ is parameterized by the mapping $\zeta_\pm^0$
restricted to $\Sigma$, while  the subsequent  two conditions in \eqref{05261240} indicate that
$\zeta_\pm^0$ map the fixed upper and lower boundaries into themselves.
Define the flow maps $\zeta_\pm$ as the solutions to the initial value problems:
\begin{equation*}
            \begin{cases}
\partial_t \zeta_\pm(y,t)=v_\pm(\zeta_\pm(y,t),t)&\mbox{in }\Omega_\pm,
\\
\zeta_\pm(y,0)=\zeta_\pm^0(y)&\mbox{in }\Omega_\pm.
                  \end{cases}
\end{equation*}
We denote the Eulerian coordinates by $(x,t)$ with $x=\zeta(y,t)$, while the fixed variables $(y,t)\in\Omega\times\mathbb{R}^+_0$
stand for the Lagrangian coordinates, where $\mathbb{R}^+_0:=[0,\infty)$.

In order to switch back and forth  between Lagrangian and  Eulerian coordinates, we temporarily assume that
$\zeta_\pm(\cdot ,t)$ are invertible and $\Omega_{\pm}(t)=\zeta_{\pm}(\Omega_{\pm},t)$. And since $v_\pm$ and
$\zeta_\pm^0$ are all continuous across $\Sigma$, one has
$\Sigma(t)=\zeta_\pm(\Sigma,t)$, i.e.,
\begin{equation}
\label{201701011211}
 \llbracket   \zeta  \rrbracket =0\;\mbox{ on }\Sigma.
\end{equation} In other words, the Eulerian domains of upper and lower fluids correspond to  the images of $\Omega_\pm$
under the mappings $\zeta_\pm$, and the free interface is the image of $\Sigma$ under the mappings $\zeta_\pm(t,\cdot)$.
 Moreover, due to  the non-slip boundary condition $v_\pm|_{\Sigma_\pm}=0$, one has
\begin{equation*} y=\zeta_\pm(y, t)\;\mbox{ on }\Sigma_\pm.
\end{equation*}
From now on, we define $\zeta:=\zeta_+$ for $y\in {\Omega_+ }$ and $\zeta_-$ for $y\in {\Omega_- }$,
 and subsequently introduce the quantity $\eta:=\zeta-y$.

Next, we introduce some notations involving $\eta$. Define $\mathcal{A}:=(\ml{A}_{ij})_{3\times 3}$ through
$\ml{A}^{\top }=(\nabla (\eta+y))^{-1}:=(\partial_j (\eta+y)_i)^{-1}_{3\times 3}$, and the differential operators $\nabla_{\ml{A}}$ and $\mm{div}_\ml{A}$ as follows:
\begin{align}
&\nabla_{\ml{A}}w:=(\nabla_{\ml{A}}w_1,\nabla_{\ml{A}}w_2,\nabla_{\ml{A}}w_3)^{\top },\quad \nabla_{\ml{A}}w_i:=(\ml{A}_{1k}\partial_kw_i,
\ml{A}_{2k}\partial_kw_i,\ml{A}_{3k}\partial_kw_i)^{\top },\nonumber \\
\label{062209345}
&\mm{div}_{\ml{A}}(f_1,f_2,f_3)^{\top }=(\mm{div}_{\ml{A}}f_1,\mm{div}_{\ml{A}}f_2,\mm{div}_{\ml{A}}f_3)^{\top },
\quad \mm{div}_{\ml{A}}f_i:=\ml{A}_{lk}\partial_k f_{il},\\
& \Delta_{\mathcal{A}}w:= (\Delta_{\mathcal{A}}w_1,\Delta_{\mathcal{A}}w_2,\Delta_{\mathcal{A}}w_3)^{\top }\;\;\mbox{ and }\;\;
\Delta_{\mathcal{A}}w_i:=\mm{div}_{\ml{A}}\nabla_{\ml{A}}w_i\nonumber
\end{align}
for vector functions $w:=(w_1,w_2,w_3)^{\top }$ and $f_i:=(f_{i1},f_{i2},f_{i3})^{\top }$,
where  the Einstein summation convention over repeated indices is used,  and $\partial_i$ denotes the partial
derivative with respect to the $i$-th component of the variable $y$.

Finally, we further introduce some properties of $\mathcal{A}$ \cite{JFJSNS}.
\begin{enumerate}[\quad\  (1)]
  \item In view of the definition of $\mathcal{A}$, one can deduce the following two important properties:
\begin{equation}
\label{AklJ=0}
\partial_j (J\mathcal{A}_{ij})=0, \end{equation}
\begin{equation}\label{AklJdeta}
  \partial_i(\eta+y)_k\ml{A}_{kj}= \ml{A}_{ik}\partial_k(\eta+y)_j=\delta_{ij},
  \end{equation}
  where $\delta_{ij}=0$ for $i\neq j$, $\delta_{ij}=1$ for $i=j$, $(\eta+y)_j$ denotes the $j$-th component of the vector $\eta+y$, and
  \begin{align}
  \label{2022306241210}
  J:=\det (\nabla (\eta+y)).
   \end{align} The relation
  \eqref{AklJ=0} is often called the geometric identity.
  \item  It holds after a straightforward calculation that
\begin{align}
\label{20230626}
J\mathcal{A}\mathbf{e}^3 =\partial_1(\eta+y)\times \partial_2(\eta+y),
\end{align}
 which, together with \eqref{201701011211}, implies
\begin{equation}
\label{06051441}
 \llbracket  J\mathcal{A}\mathbf{e}^3  \rrbracket =0.\end{equation}
  \item Let
\begin{align}
\label{05291021n}
 \mathbf{n}= \frac{J\mathcal{A}\mathbf{e}^3}{|J\mathcal{A}\mathbf{e}^3|}.
\end{align} By virtue of \eqref{06051441}, $\mathbf{n}|_{y_3=0}$ is the unit normal to $\Sigma(t)=\zeta(\Sigma,t)$.
\end{enumerate}

Let
$$
 1_{\Omega_\pm} =1_{\Omega_\pm}(y) =
\begin{cases}
1&\mbox{for }y\in \Omega_\pm ,\\
0 &\mbox{for }y\in  \Omega_\mp.
\end{cases}$$
We denote
$$P^{(n)}(\psi)=\begin{cases}
P_+^{(n)}(\psi_+(y))&\mbox{for }y\in {\Omega_+},\\
P_-^{(n)}(\psi_-(y))&\mbox{for }y\in {\Omega_-},
\end{cases}
$$
where $P^{(n)}_\pm(\cdot)$ are the $n$-order derivatives of $P_\pm$.

 Now, let us set the Lagrangian unknowns
\begin{equation*}
(\varrho,u)(y,t)=(\rho,v)(\eta(y,t)+y,t) \;\; \mbox{ for } (y,t)\in \Omega \times\mathbb{R}^+
\end{equation*}
and re-define
 \begin{align}
& \label{202307212034} \mu:=\mu_+1_{\Omega_+} +\mu_-1_{\Omega_-}\;\mbox{ and }\; \varsigma:=\varsigma_+1_{\Omega_+} +\varsigma_-1_{\Omega_-},
\end{align}
the SCVF model can thus be written as the following initial-boundary value problem with an interface for $(\eta,\varrho,u )$
in Lagrangian coordinates \cite{GYTI2}:\begin{equation}\label{201611041430M}
\begin{cases}
\eta_t=u &\mbox{in } \Omega,\\[1mm]
\varrho_t+\varrho\mm{div}_{\mathcal{A}}u=0 &\mbox{in } \Omega, \\[1mm]
\varrho u_t+\mm{div}_{\ml{A}}(P(\varrho)\mathbb{I}-\mathbb{S}_{\mathcal{A}}(u)) =-\varrho g \mathbf{e}^3&\mbox{in }  \Omega,\\[1mm]  \llbracket  P(\varrho)\mathbb{I}-\mathbb{S}_{\mathcal{A}} (u)  \rrbracket  \mathbf{n} = \vartheta \mathcal{H} \mathbf{n} ,\
 \llbracket  \eta   \rrbracket = \llbracket  u   \rrbracket =0 &\mbox{on }\Sigma,\\
  (\eta,u)=0 &\mbox{on }\partial\Omega\!\!\!\!\!-,\\
(\eta, \varrho,u)|_{t=0}=(\eta^0,\varrho^0, u^0 ),  &\mbox{in }  \Omega .
\end{cases}\end{equation}
Here the jump notations $\llbracket \cdot \rrbracket $  in \eqref{201611041430M} are defined by \eqref{202307302006}, and $\partial\Omega\!\!\!\!\!-:=\Sigma_+\cup\Sigma_-$, $\Omega:=\Omega_+\cup\Omega_-$, and
\begin{align}
&H^{\mm{n}}:=|\partial_1 \zeta|^2\partial_2^2\zeta-2(\partial_1\zeta\cdot \partial_2\zeta)\partial_1\partial_2\zeta+
|\partial_2\zeta|^2\partial_1^2\zeta,   \label{20sdfs1611041430M} \\
&H^{\mm{d}}:= |\partial_1\zeta|^2|\partial_2\zeta|^2-|\partial_1\zeta\cdot \partial_2\zeta|^2,   \   \mathcal{H}:= {H^{\mm{n}}}\cdot \mathbf{n}/{H^{\mm{d}}}  ,  \label{20161asdfa1041430M}  \\
&\mathbb{S}_{\mathcal{A}} (u) := \mu \mathbb{D}_{\mathcal{A}} u+\left(\varsigma -{2\mu}/{3}\right)\mm{div}_{\mathcal{A}} u \mathbb{I},
\;\mbox{ and } \label{x2022306090918}
\\ &\mathbb{D}_{\mathcal{A}} u:=\nabla_{\mathcal{A}} u+\nabla_{\mathcal{A}} u^{\top } . \label{2022306090918}
\end{align}

Since our aim is to construct  unstable solutions of  the  problem \eqref{201611041430M},
it is better to simplify the problem  as much as possible.
Thus, our next goal is to eliminate $\varrho$ in \eqref{201611041430M} by expressing it in terms of $\eta$ as in \cite{JFJSNS,TZWYJGw}.

 It follows from \eqref{201611041430M}$_1$ that
\begin{equation}\label{Jtdimau}
J_t=J\mm{div}_{\mathcal{A}}u,             \end{equation}
which, together with \eqref{201611041430M}$_2$, yields  the mass conservation of differential version
\begin{equation}\label{0122sd}\partial_t(\varrho J)=0.   \end{equation}
 Hence, we deduce from \eqref{0122sd} that
$\varrho^0J^0=\varrho J $,
which implies $\varrho = \bar{\rho}J^{-1}$, provided the initial data $(\eta^0,\varrho^0)$ satisfies
\begin{align}
&\label{abjlj0i}  \varrho^0 =\bar{\rho}/J^0  .
\end{align}

Consequently, under the assumption  \eqref{abjlj0i}, we derive the following (compressible) RT problem from the initial-boundary value problem \eqref{201611041430M}:
\begin{equation}\label{n0101nnnM}\begin{cases}
\eta_t=u &\mbox{in } \Omega,\\[1mm]
\bar{\rho}J^{-1} u_t+\mm{div}_{\ml{A}}(P(\bar{\rho}J^{-1})\mathbb{I}-\mathbb{S}_{\mathcal{A}}(u)) =-\bar{\rho}J^{-1} g \mathbf{e}^3&\mbox{in }  \Omega,\\[1mm]
 \llbracket P(\bar{\rho}J^{-1})\mathbb{I}-\mathbb{S}_{\mathcal{A}}( u)    \rrbracket  \mathcal{A}\mathbf{e}^3 = \vartheta \mathcal{H} \mathcal{A}\mathbf{e}^3 &\mbox{on }\Sigma,\\  \llbracket \eta   \rrbracket =\llbracket  u  \rrbracket =0 &\mbox{on }\Sigma,\\
(\eta,u)=0  &\mbox{on }\partial\Omega\!\!\!\! -,\\
 (\eta,  u)|_{t=0}=(\eta^0,  u^0)  &\mbox{in }  \Omega,
\end{cases}\end{equation}
where, by virtue of \eqref{202307212034},
\begin{align}
\mu 1_{\Omega_\pm}=\mu_\pm,\ \varsigma 1_{\Omega_\pm}=\varsigma_\pm\;
\mbox{ and }\; P(\bar{\rho}J^{-1})1_{\Omega_\pm}=P_\pm(\bar{\rho}_\pm /\det (\nabla (\eta 1_{\Omega_\pm} +y)) ).
\end{align}
\emph{It should be noted that our constructed unstable solutions to  the above RT problem automatically
implies the instability of the problem \eqref{201611041430M} with $\varrho=\bar{\rho}/J$ and $J\mathcal{A}\mathbf{e}^3\neq 0$, see Remark \ref{202230800212201}.}

Next, we further deduce the perturbation forms of \eqref{n0101nnnM}$_2$ and \eqref{n0101nnnM}$_3$.
To begin with, we should rewrite the pressure $P(\bar{\rho}J^{-1})$ as some perturbation form around $\bar{P}$.
By a simple computation, we have
\begin{equation} \label{201611050816}
 P(\bar{\rho}J^{-1})=\bar{P}
-{P}'(\bar{\rho})\bar{\rho}\mm{div}\eta+ R_P,
\end{equation}
where $\bar{P}:=P(\bar{\rho})$ and
\begin{align}\label{2022307251842}
R_P := {P}'(\bar{\rho})\bar{\rho}(J^{-1}-1+\mm{div}\eta)
+\int_{0}^{\bar{\rho}(J^{-1}-1)}(\bar{\rho}
(J^{-1}-1)-z)\frac{\mm{d}^2}{\mm{d}z^2} P (\bar{\rho}+z)\mm{d}z.
\end{align}
Let $\tilde{\bar{P}}:=P(\tilde{\bar{\rho}})$ and $\tilde{\bar{\rho}}:=\bar{\rho}(\eta_3+y_3)$. Thus, the hydrostatic relation in
\eqref{201611051547}$_1$ in Lagrangian coordinates reads as
$\nabla_{\mathcal{A}}\tilde{\bar{P}}   =  -\tilde{\bar{\rho}} {g} \mathbf{e}^3.$
Using \eqref{201611051547}$_1$, one has
\begin{align}\mm{div}_{\mathcal{A}}(\bar{P}\mathbb{I})=\nabla_{\mathcal{A}}\bar{P}= -\tilde{\bar{\rho}} {g}\mathbf{e}^3-
\nabla_{\mathcal{A}}(\tilde{\bar{P}}-\bar{P})= g\mm{div}_{\mathcal{A}}(\bar{\rho}\eta_3 \mathbb{I} ) -\tilde{\bar{\rho}} {g}\mathbf{e}^3+\mathbf{N}_P,
\label{202306081302}
\end{align}
where
$$\mathbf{N}_P :=\nabla_{\mathcal{A}}\left( \int_{0}^{\eta_3}
(z-\eta_3)\frac{\mm{d}^2}{\mm{d}z^2} \bar{P}(y_3+z)\mm{d}z\right). $$

In addition, we represent the  gravity term as follows.
\begin{equation}
\label{gravit20160309n}
-\bar{\rho}J^{-1}  g\mathbf{e}^3=-(\bar{\rho}(J^{-1}-1)+\bar{\rho}-\tilde{\bar{\rho}}+\tilde{\bar{\rho}})
{g}\mathbf{e}^3 =g \mm{div}(\bar{\rho}\eta)\mathbf{e}^3- \tilde{\bar{\rho}}
{g}\mathbf{e}^3  +\mathbf{N}_g,   \end{equation}
where
 $$\mathbf{N}_g :=g\left(  \int_{0}^{\eta_3} (\eta_3-z) \frac{\mm{d}^2}{\mm{d}z^2}\bar{\rho}( y_3+z) \mm{d}z
 -\bar{\rho} ( J^{-1}-1+\mm{div}\eta)\right)  \mathbf{e}^3.$$
Thus, by \eqref{201611051547}$_2$, \eqref{201611050816}, \eqref{202306081302} and \eqref{gravit20160309n}, one can
transform \eqref{n0101nnnM}$_2$ and \eqref{n0101nnnM}$_3$ into the following equivalent forms:
\begin{equation}\label{201611040926M} \begin{cases}
\bar{\rho}J^{-1}  u_t-\mm{div}_{\ml{A}} (   {P}'(\bar{\rho})\bar{\rho}\mm{div}\eta \mathbb{I} +\mathbb{S}_{\mathcal{A}}(u))=g\bar{\rho}(\mm{div}\eta \mathbf{e}^3-\nabla \eta_3 ) + \mathbf{N}^1 &\mbox{ in }  \Omega,\\[1mm]
  \llbracket   {P}'(\bar{\rho})\bar{\rho}\mm{div}\eta\mathbb{I} +\mathbb{S}_{\mathcal{A}}(u)\rrbracket  J\mathcal{A}\mathbf{e}^3 +\vartheta \Delta_{\mm{h}}\eta_3\mathbf{e}^3  = {\mathbf{N}^2} &\mbox{ on }\Sigma,
\end{cases}\end{equation}
 where  $\Delta_{\mm{h}}:=\partial_1^2+\partial_2^2$, and
 \begin{align}
  {\mathbf{N}^1} :=   \mathbf{N}_g-{\mathbf{N}}_P-\nabla_{\mathcal{A}} {R}_P -g\nabla_{\tilde{\mathcal{A}}}(\bar{\rho}\eta_3)  \mbox{ and } {\mathbf{N}^2} := \llbracket  R_P \rrbracket J  \mathcal{A}\mathbf{e}^3+  \vartheta (   \Delta_{\mm{h}}\eta_3 \mathbf{e}^3- \mathcal{H}J  \mathcal{A}\mathbf{e}^3 ).
  \label{2022307260853}
   \end{align}
 We will use \eqref{201611040926M} to derive the \emph{a priori} estimates for the temporal derivative of $u$.

Now we further rewrite \eqref{201611040926M} as non-homogeneous perturbation forms.
Letting $\mathbf{f}$, $\mathbf{0}\not =\mathbf{r}\in\mathbb{R}^3$, we define $\Pi_{\mathbf{r}}\mathbf{f}:=\mathbf{f}-|\mathbf{r}|^{-2}(\mathbf{f}\cdot \mathbf{r})\mathbf{r}  $. It should be noted that $\Pi_{\mathbf{r}}\mathbf{f}=0$, if only if the vector $\mathbf{f}$ parallels to $\mathbf{r}$.
Applying the operator $\Pi_{\mathbf{n}}$ to the jump condition \eqref{201611040926M}$_2$, and then using the fact that ${\mathbf{n}}$ parallels to $J \mathcal{A} \mathbf{e}^3 $, we get
$$ \Pi_{\mathbf{n}} \llbracket  (  {P}'(\bar{\rho})\bar{\rho}\mm{div}\eta\mathbb{I}+ \mathbb{S}_{\mathcal{A}} (u) ) J \mathcal{A} \mathbf{e}^3    \rrbracket  =0\mbox{ on }\Sigma.$$
We rewrite the above identity as a nonhomogeneous form,  on the left hand side of which the terms  are linear:
\begin{align}
\Pi_{\mathbf{e}^3}  \llbracket \Upsilon(\eta,u) \mathbf{e}^3   \rrbracket  =\mathbf{N}^4\mbox{ on }\Sigma,
\label{201806211105}
\end{align}
where
\begin{align}
& \Upsilon(\eta,u):=  {P}'(\bar{\rho})\bar{\rho}\mm{div}\eta\mathbb{I}+\mathbb{S}   (u)  ,
\label{2022306161627}
\\
&\mathbf{N}^4:=\mathbf{N}^4(\eta,u):=\llbracket ( \Upsilon  (\eta,u) \mathbf{e}^3)\cdot\tilde{\mathbf{n}}\mathbf{n} + ( \Upsilon  (\eta,u)\mathbf{e}^3 ) \cdot \mathbf{e}^3\tilde{\mathbf{n}}\rrbracket\nonumber \\
  &\qquad \qquad\qquad \qquad-
\Pi_{ \mathbf{n} } \llbracket  \Upsilon  (\eta,u)   (J{\mathcal{A}}\mathbf{e}^3-\mathbf{e}^3)+
  \mathbb{S}_{\tilde{\mathcal{A}}}(u )J\mathcal{A}\mathbf{e}^3)\rrbracket,  \label{2saf022306161627} \\
&\tilde{\mathbf{n}}:=\tilde{\mathbf{n}}-\mathbf{e}^3 ,\   \tilde{\mathcal{A}}:=\mathcal{A}-\mathbb{I}\mbox{ and }\mathbb{S}_{\tilde{\mathcal{A}}}(u)
 \mbox{ is  defined by \eqref{x2022306090918} with }\tilde{\mathcal{A}}\mbox{ in place of }\mathcal{A}.
  \end{align}

Taking scalar product of \eqref{n0101nnnM}$_3$ and $J\mathbf{n}/|J \mathcal{A}\mathbf{e}^3|$, we have
$$  \llbracket  P(\bar{\rho}J^{-1})  - (\mathbb{S}_{\mathcal{A}}( u ) \mathbf{n} )  \cdot\mathbf{n}   \rrbracket =\vartheta \mathcal{H},$$
which can be rewritten as a nonhomogeneous form by  \eqref{201611051547}$_2$ and \eqref{201611050816} \begin{equation}
\label{201806202203}
  \llbracket
{P}'(\bar{\rho})\bar{\rho}\mm{div}\eta+2\mu \partial_3u_3 +\left(\varsigma-{2\mu}/{3}\right)\mm{div}u  \rrbracket +\vartheta\Delta_{\mm{h}}\eta_3= \mathcal{N},\end{equation}
where
\begin{align}
&\mathcal{N}:= \mathcal{N}(\eta,u):= \mathcal{N}^\eta+\llbracket R_P +\mathcal{N}^u \rrbracket  ,\label{202306161628} \\
&\mathcal{N}^\eta(\eta): =\vartheta( \Delta_{\mm{h}}\eta_3-\mathcal{H})\nonumber \\
&\qquad   \ \ \  =\vartheta(H^{\mm{n}} \cdot \mathbf{n}(H^{\mm{d}}-1)/H^{\mm{d}}-H^{\mm{n}}\cdot \tilde{ \mathbf{n}}-H^{\mm{n}}\cdot \mathbf{e}^3 +\Delta_{\mm{h}}\eta_3), \label{202005151045} \\
& \mathcal{N}^u:= \mathcal{N}^u(\eta,u):=- \mathbb{S}_{\tilde{\mathcal{A}}}(u )    \mathbf{n} \cdot  \mathbf{n}-
    (\mathbb{S}(u)\tilde{ \mathbf{n}})\cdot \mathbf{n} -(\mathbb{S}(u)  \mathbf{e}^3)\cdot \tilde{ \mathbf{n}},\nonumber \\
&\mbox{and } \mathbb{S}(u)\mbox{ is defined by \eqref{x2022306090918} with }\mathbb{I}\mbox{ in place of }\mathcal{A}.\nonumber
\end{align}

Thanks to \eqref{201806211105} and \eqref{201806202203}, we further
rewrite  \eqref{201611040926M} as
the following nonhomogeneous form:
\begin{equation}\label{n0101nn1928M} \begin{cases}
  \bar{\rho}u_t +  g\bar{\rho}(\nabla \eta_3- \mm{div}\eta \mathbf{e}^3)-\mm{div} \Upsilon(\eta,u)    = {\mathbf{N}}^3  &\mbox{ in }\Omega, \\
 \llbracket\Upsilon(\eta,u) \mathbf{e}^3  \rrbracket +\vartheta\Delta_{\mm{h}}\eta_3\mathbf{e}^3  =(\mathbf{N}_1^4,\mathbf{N}_2^4,\mathcal{N})^{\top }&  \mbox{ on }\Sigma,
\end{cases}\end{equation}
 where
  \begin{align}
   & {\mathbf{N}^3} :={\mathbf{N}^1} +\mm{div}_{\tilde{\ml{A}}}\mathbb{S}_{ {\mathcal{A}}}(u) +\mm{div}  \mathbb{S}_{\tilde{\mathcal{A}}}(u) +\nabla_{\tilde{\mathcal{A}}}({P}'(\bar{\rho})\bar{\rho}\mm{div}\eta ) +{\bar{\rho}}(1-J^{-1} ) u_t,
\end{align}
and $\mm{div}_{\tilde{\mathcal{A}}}$ is defined by \eqref{062209345} with $\tilde{\mathcal{A}}$ in place of $\mathcal{A}$.   The above nonhomogeneous forms in \eqref{n0101nn1928M} are very useful to derive the \emph{a priori} estimates of tangential derivatives of $(\eta,u)$, in particular, a highest-order  boundary estimate of $ u_3  $ at interface for $\vartheta>0$.

Thanks to \eqref{n0101nn1928M}, we can easily write out the linearized  RT problem   corresponding  to \eqref{n0101nnnM} under the small perturbations:
\begin{equation}\label{linearizedxx} \begin{cases}
\eta_t=u &\mbox{ in } \Omega,\\[1mm]
\bar{\rho} u_t = g\bar{\rho}(\mm{div}\eta \mathbf{e}^3- \nabla \eta_3) +\mm{div}  \Upsilon(\eta,u)&\mbox{ in }  \Omega,\\[1mm]
 \llbracket   \eta  \rrbracket  =\llbracket   u  \rrbracket  =0,\  \llbracket     \Upsilon(\eta,u)  \mathbf{e}^3\rrbracket   =- \vartheta \Delta_{\mm{h}}\eta_3  \mathbf{e}^3&\mbox{ on }\Sigma,\\
(\eta, u)=0 &\mbox{ on }\partial\Omega\!\!\!\!\!-,\\
 (\eta,  u)|_{t=0}= (\eta^0,  u^0 ) &\mbox{ in }  \Omega.
\end{cases}
\end{equation}
Moreover, the corresponding spectrum problem of \eqref{linearizedxx} reads as follows.
\begin{equation}\label{spectrumequations} \begin{cases}
\lambda \sigma =w &\mbox{ in } \Omega,\\[1mm]
\lambda\bar{\rho}w= g\bar{\rho}(\mm{div}\sigma \mathbf{e}^3-\nabla \sigma_3 )+\mm{div}   \Upsilon(\sigma,w)&\mbox{ in }  \Omega,\\[1mm]
  \llbracket \sigma  \rrbracket  = \llbracket  w   \rrbracket  =0,\
   \llbracket  \Upsilon(\sigma,w)  \mathbf{e}^3\rrbracket=- \vartheta \Delta_{\mm{h}}\sigma_3 \mathbf{e}^3&\mbox{ on }\Sigma,\\
(\sigma, w)=0 &\mbox{ on }\partial\Omega\!\!\!\!\!\! - \! .
\end{cases}\end{equation}

 The linearized RT problem \eqref{linearizedxx} has the advantage of being convenient in mathematical analysis in order to have an insight into the physical  mechanism  of the instability, and how the internal surface tension affects the RT instability. In particular, an instability solution of the  linearized  RT problem plays a role of standing point to further investigate the nonlinear instability.
We mention that directly linearizing the   problem \eqref{201611041430M} yields the following linear problem  \cite{GYTI2}:
\begin{equation}\label{saf} \begin{cases}
\eta_t=u &\mbox{ in } \Omega,\\[1mm]
\varrho_t=-\bar{\rho}\mm{div}u&\mbox{ in } \Omega,\\[1mm]
\bar{\rho} u_t =\mm{div} (\mathbb{S}(u)-P'(\bar{\rho})\varrho\mathbb{I }) -g(\varrho\mathbf{e}^3+\bar{\rho}\nabla \eta_3 ) &\mbox{ in }  \Omega,\\[1mm]
 \llbracket   \eta  \rrbracket  =\llbracket   u  \rrbracket  =0,\  \llbracket   \mathbb{S}(u)-P'(\bar{\rho})\varrho\mathbb{I }  \rrbracket   \mathbf{e}^3=- \vartheta \Delta_{\mm{h}}\eta_3  \mathbf{e}^3&\mbox{ on }\Sigma,\\
(\eta, u)=0 &\mbox{ on }\partial\Omega\!\!\!\!\!-,\\
 (\eta,  u)|_{t=0}= (\eta_0,  u_0 ) &\mbox{ in }  \Omega,
\end{cases}
\end{equation}
which is called the \emph{linearized (non-periodic) m-RT problem}. Here we add the letter $m$ in the name of linearized m-RT problem to emphasize that it includes the linearized mass equation \eqref{saf}$_2$.

\subsection{Relevant progress and our main result}\label{202306071229}

Let  $\mathbb{T}$ be a usual 1-torus,
\begin{align}
\Omega_{L_1,L_2}=2\pi L_1\mathbb{T}\times 2\pi L_2\mathbb{T} \times ((h_-,0)\cup (0,h_+))\mbox{ and }\Sigma_{L_1,L_2} = 2\pi L_1\mathbb{T}\times 2\pi L_2\mathbb{T}\times \{0\} ,
\label{2022307222447}
 \end{align}then  the linearized RT problem \eqref{linearizedxx}
with $\Omega_{L_1,L_2}$ resp. $\Sigma_{L_1,L_2}$ in place of $\Omega $ resp. $\Sigma $ is called the \emph{linearized  periodic RT problem}.
Then the solution  $(\eta,u)$ of the linearized periodic RT problem satisfies \eqref{linearizedxx} and $$(\eta(y,t),u(y,t))= (\eta(x,t),u(x,t))|_{(x_1,x_2)=(y_1+2\pi m L_1,y_2+2\pi n L_2)}  $$
 for any   integers $m$ and $n$. Similarly the linearized m-RT problem \eqref{saf}
with $\Omega_{L_1,L_2}$ resp. $\Sigma_{L_1,L_2}$ in place of $\Omega $ resp.  $\Sigma $ is called the \emph{linearized   periodic m-RT problem}. By applying a so-called ``normal mode" method to the linearized periodic m-RT problem, Guo--Tice obtained both the linear stability and instability results \cite{GYTI2}:
 \begin{itemize}
   \item   the linearized periodic m-RT problem is stable, if the  Rayleigh number satisfies
\begin{equation}
\label{20223006082126}
R:= \frac{\vartheta}{g\max\{L_1^2,L_2^2\}\llbracket {\bar{\rho}}\rrbracket}>1.
\end{equation}
   \item   the linearized periodic m-RT problem is unstable, if the instability condition is satisfied
\begin{equation}
\label{20223saf006082126}R<1 .
  \end{equation}
 \end{itemize}

In addition, they also \emph{formally} verified that the linearized (non-periodic) m-RT problem \eqref{saf} is unstable, see \cite[Theorem 2.4]{GYTI2} (Here we emphasize the word ``\emph{formally}", since both the integral expressions  in (2.24) and (2.25) in \cite[Theorem 2.4]{GYTI2} are formal). Obviously, the instability of the linearized m-RT problem \eqref{saf}   can be expected from the periodic case, since the domain $\Omega$ can be regarded the limit case of $\Omega_{L_1,L_2}$, as $L_1$, $L_2\to \infty$. We mention that all Guo--Tice's results of the linearized m-RT problem also hold for the   linearized periodic/non-periodic RT problem.

Guo--Tice's linear stability result  roughly presents  that    the interfacial  surface tension can inhibit the occurrence of the RT instability (under the case of small perturbation with horizontally periodic motion). However the rigorous proof of \emph{the phenomenon of  inhibition of RT instability by the interfacial  surface tension} is relatively difficult in Lagrangian coordinates. Therefore Jang--Tice--Wang used the flattening transformation introduced by Beal in \cite{beale1984large} to rewrite the motion equations of stratified  compressible viscous fluids, and then  succeeded in the mathematical verification of the inhibition phenomena under the  sharp  stability condition \eqref{20223006082126}  by further overcoming some difficulties arising from the compressibility, in particular, the stability solutions under the  inhibition  of surface tension enjoy the  exponential decay-in-time \cite{JJTIWYJ}. Such conclusion also was obtained for the corresponding incompressible case by Wang--Tice--Kim \cite{WYJTIKCT}.
Recently, Wilke further gave a similar result for the stratified incompressible viscous fluids in a bounded  cylindrical domain $G\times (h_-,h_+)$  under the sharp stability condition
\begin{equation}
\label{20223006082126x}
R_s:=\frac{\lambda\vartheta}{ g (\rho_+^{\mm{i}}-\rho_-^{\mm{i}})}>1,
\end{equation}
where $G\subset \mathbb{R}^2$, $\rho_\pm^{\mm{i}}$ are the constant densities of  upper and lower incompressible fluids, respectively, and the constant $\lambda$ depends on the geometric structure of $G$; in particular, if the cylindrical domain is a cylinder, the expression of $R_s$ in \eqref{20223006082126x} is given by
$R_s=\vartheta\lambda^*/ g r^2(\rho_+^{\mm{i}}-\rho_-^{\mm{i}} )$, where $r $ is the radius of the cylinder, $ {\lambda^*}=(\mathfrak{B}'_{1,1})^2$ and $\mathfrak{B}'_{1,1}$  is the first zero of the derivative $\mathfrak{B}_1'$ of the Bessel function $\mathfrak{B}_1$ \cite{wilke2022rayleigh}.

Thanks to Guo--Tice's linear instability result, Jang--Tice--Wang \cite{jang2016compressiblexx} further established the nonlinear RT instability result that the RT equilibrium solution $(\bar{\rho},\mathbf{0})$ is unstable to the  \emph{periodic} SCVF model with the instability condition \eqref{20223saf006082126} based on   Guo--Hallstrom--Spirn's bootstrap instability method in \cite{GYHCSDDC}, also see the result for the corresponding incompressible case \cite{wang2011viscous}. However  Guo--Tice's method for the construction of instability solutions of the  linearized periodic RT problem  can not be directly applied to   the corresponding non-periodic case \eqref{linearizedxx}. Next let us take the case $\vartheta=0$ to roughly explain the reason of failure.

Let us rewrite the  two abstract forms of the equations in \eqref{linearizedxx} and \eqref{spectrumequations} with $\vartheta=0$:
\begin{align}\frac{\mm{d}V}{\mm{d}t}=\mathcal{L}({V})
\label{2022300609212sdaf6}
\end{align}
and\begin{align}\lambda U=\mathcal{L}(U),
\label{20223006092126}
\end{align}
where ${V} =(\eta^{\top}, u^{\top})^{\top}$, $U=(\sigma^{\top}, w^{\top})^{\top}$ and $\mathcal{L}$ denotes a linear differential operator. For a function $f\in L^2 $, we define the horizontal Fourier transform of $f$ in
  $\mathbb{R}^2$ via
\begin{equation}\label{hftx}
 \hat{f}(\xi,y_3) := \mathcal{F}_{y_{\mm{h}}\to \xi}(f(y_{\mm{h}},y_3))=\frac{1}{2\pi} \int_{\mathbb{R}^{2}} f(y_{\mm{h}},y_3) e^{-\mm{i}\xi \cdot y_{\mm{h}} }\mm{d}y_{\mm{h}}
\end{equation}
for $\xi \in \mathbb{R}^2$, where $y_{\mm{h}}:=(y_1,y_2)$ and $\xi\cdot y_{\mm{h}}:=\xi_1 y_1+ \xi_2 y_2$. Applying the horizontal Fourier transform to \eqref{20223006092126} yields
\begin{equation}\label{hft}
\lambda(\xi) \hat{U}(\xi,y_3) = \mathbb{P}\left(\xi,\frac{\mm{d}}{\mm{d}y_3}\right)\hat{U}(\xi,y_3),
\end{equation}
where $ \mathbb{P}\left(\xi,\frac{\mm{d}}{\mm{d}y_3}\right):=\left(\mathcal{P}_{ij}\left(\xi_1,\xi_2,\frac{\mm{d}}{\mm{d}y_3}\right)\right)_{6\times 6}$, all entry $\mathcal{P}_{ij}\left(\xi_1,\xi_2,\frac{\mm{d}}{\mm{d}y_3}\right)$ are at most two order complex polynomials with respect to $\xi_1$, $\xi_2$ and the derivative operator $\frac{\mm{d}}{\mm{d}y_3}$ applying some component of $\hat{U}(\xi,y_3)$, where $1\leqslant i$, $j\leqslant 6$,   $\left(\frac{\mm{d}}{\mm{d}y_3}\right)^l:=\frac{\mm{d}^l}{\mm{d}y_3^l}$  and the coefficients of
the complex polynomials depends on $g$, $\mu$, $\varsigma$, $\bar{\rho}$,  $P'(\bar{\rho})$ and $(P'(\bar{\rho})\bar{\rho})'$.  Exploiting a modified variational method introduced by Guo--Tice in \cite{GYTI2}, we can construct  solutions $\hat{U}(\xi,y_3)$ and $\lambda(\xi)$ to \eqref{hft} (with some boundary conditions, which can be obtained by applying the horizontal Fourier transform to the boundary conditions \eqref{spectrumequations}$_3$ and \eqref{spectrumequations}$_4$) for any given   frequency  $\xi$, where $ \lambda(\xi)>0$ is a bounded and continuous function on $\xi$. Formally $V=\mathcal{F}_{\xi\to y_{\mm{h}} }^{-1}( e^{\lambda(\xi)t}\hat{U}(\xi,y_3))$ is a growing mode  solution of \eqref{2022300609212sdaf6}, since $\|V\|_{L^2(\Omega )}= \| e^{\lambda(\xi)t}\hat{U}(\xi,y_3)\|_{L^2(\Omega )}$ and $\lambda(\xi)>0$.

The above idea seems to provide a road to construct a growing mode solution. Unfortunately, it is difficult to prove that $\hat{U}(\xi,y_3)$ is a measurable function on the variable $(\xi,y_3)$, and thus we can not  apply Fourier inverse transform to further construct a  growing mode solution, in particular, enjoying an almost largest growth rate (i.e. close to the largest growth rate of all linear solutions roughly speaking) to the linearized RT problem \eqref{linearizedxx}. It should be noted that solutions $\hat{U}(\xi,y_3)$ and $\lambda(\xi)$ satisfying \eqref{hft} can be also obtained for the case $\vartheta>0$, where  $\xi$ should satisfy $|\xi|^2< g\llbracket \bar{\rho}\rrbracket/\vartheta$.
Based on this fact and \emph{an assumption of  $\hat{U}(\xi,y_3)$ being measurable},  Guo--Tice formally obtained the conclusion that the linearized m-RT problem \eqref{saf} is unstable (see \cite[Theorem 2.4]{GYTI2}), and their's linear result roughly presents that  RT instability always occur in the SCVF model with non-periodic motion.  It is worth to be noted that Pr$\mathrm{\ddot{u}}$ess--Simonett  had verified such conclusion in the case of  stratified incompressible viscous fluids in $\mathbb{R}^3$  by the spectrum analysis of the linearized problem,  the maximal regularity theory of type $L^p$ and  the Henry's instability theorem \cite{PJSGOI5}. Later Wilke also used Pr$\mathrm{\ddot{u}}$ess--Simonett's method to prove the RT instability in a cylindrical domain  under the instability condition $R_s<1$ \cite{wilke2022rayleigh}. \emph{However it is still not clear whether Pr$\mathrm{\ddot{u}}$ess--Simonett's method can be further applied to the corresponding compressible case to obtain the solutions with instability of $L^2$-norm}.

In this article, we still use Guo--Hallstrom--Spirn's bootstrap instability method in \cite{GYHCSDDC} to investigate the  instability of the   (non-periodic) RT problem \eqref{n0101nnnM}, and thus face the difficulty of the  construction of a growing mode solution  for the linearized RT problem \eqref{linearizedxx} with an almost largest growth rate for the linearized problem discussed previously.
However we do not try to fix the construction problem, and develop a new alternative scheme to avoid it, i.e. roughly speaking, using  a growing mode solution  with an almost largest growth rate of the periodic RT problem  to approximately construct the   RT instability solutions for the RT problem \eqref{n0101nnnM}. Thanks to the new approximate scheme, we can rigorously prove the instability of the RT problem, and thus extend Guo--Tice's linear conclusion to the nonlinear case. The brief proof frame will be further mentioned after Theorem \ref{thm:0202}. Our RT instability result presents that  the RT equilibrium solution $(\bar{\rho},\mathbf{0})$ is   always unstable to the SCVF model with non-periodic motion, see Remark \ref{202300808400}.

Before further stating our RT instability result, we will introduce some simplified notations in this article.

(1) Basic notations.

$\partial_\mm{h}^\alpha$ denotes $\partial_{1}^{\alpha_1}\partial_{2}^{\alpha_2}$
with $\alpha=(\alpha_1,\alpha_2)$. $\partial_\mm{h}^i$ represents $\partial_{\mm{h}}^{\alpha}$ for any $|\alpha|:=\alpha_1+\alpha_2=i\geqslant 0$. $\mathbb{R}^+:=(0,\infty)$, $\int :=\int_\Omega$, $\mathbf{f}_{\mm{h}}:=(\mathbf{f}_1,\mathbf{f}_2)^{\top}$, $\nabla_{\mm{h}}f:=(\partial_1f,\partial_2f)^{\top}$   and $\mm{div}_{\mm{h}}\mathbf{f}_{\mm{h}}=\partial_1\mathbf{f}_1+\partial_2\mathbf{f}_2$, where $\mathbf{f}=(\mathbf{f}_1,\mathbf{f}_2,\mathbf{f}_3)^{\top}$.   For the sake of the simplicity, we denote $\sqrt{\sum_{i=1}^n\|w_i\|_X^2}$ by $\|(w_1,\cdots,w_n)\|_X$, where $\|\cdot\|_X$ represents a norm or a semi-norm, and $w_i$ is a scalar function or a vector function for $1\leqslant i\leqslant n$.
We define  the fractional differential operator
\begin{equation}\mathfrak{D}_{\mf{h}}^{3/2}w:= (w(y+\mf{h}  )-w(y))/|\mf{h}|^{3/2}\mbox{ for }\mf{h}\in \mathbb{R}^2\times \{0\}.
\label{2022306231307}
 \end{equation}
$a\lesssim b$ means that $a\leqslant cb$ for some constant $c>0$, where the positive constants  $c$  may
depend on the domain $\Omega$, and other known physical parameters/functions in the RT problem \eqref{n0101nnnM}  such as
$g$, $\vartheta$, $\mu$, $\varsigma$, $\bar{\rho}$, $h_\pm$ and $P_\pm(\tau)$, and may vary from line to line. Similarly the fixed constants $c_i$,  $\tilde{c}_j$ and $\delta_k$ also depend on the domain $\Omega$, and the other known physical parameters/functions in the RT problem for $1\leqslant i\leqslant 8$,  $1\leqslant j\leqslant 6$ and $1\leqslant k\leqslant 4$.

(2) Simplified notations of the Sobolev spaces and norms.
  \begin{align}&
L^p:=L^p (\Omega)=W^{0,p}(\Omega),\
{H}^k:=W^{k,2}(\Omega ),\ {H}^1_0:=W^{1,2}_0(\Omega\!\!\!\!- \, ),\
{H}^s:=W^{s,2}(\mathbb{R}^2 ), \nonumber \\
&\|\cdot \|_k :=\|\cdot \|_{H^k},\ |\cdot|_{y_3=0}|_{s} := \|\cdot  \|_{H^{s}(\mathbb{R}^2)},   \ \| \cdot\|_{i,k}^2:=\sum_{\alpha_1+\alpha_2=i} \|\partial_{1}^{\alpha_1}\partial_{2}^{\alpha_2}\cdot\|_{k}^2,\nonumber\\
&\|\cdot\|^2_{\underline{i},k}:=\sum_{j=0}^i\|\cdot\|_{j,k}^2 ,\ H^{3+k,1/2}_{0,*}:=\{\xi\in H_0^1 \cap H^{3+k}~|~ \|\xi(y)\|_3\leqslant \iota\}, \nonumber
\end{align}
where $1< p\leqslant  \infty$,  $s$ is a real number, $i$, $k$ are non-negative integers, $\iota$ is the constant in Lemma \ref{201809012320}  and the definition  of $W^{s,2}$ can be found in  \cite[Section 1.3.5.10]{NASII04} (the equivalent definition of the norm $|\cdot|_{1/2}$ by Fourier transform can be found in \cite[Theorem 7.63]{RAADSBS}). It should be noted that our constructed nonlinear instability solution $\eta$ belongs to $ H^{3,1/2}_{0,*}$ for each $t\geqslant 0$, and thus $\zeta  $ is invertible by Lemma \ref{201809012320}.

(3) Definitions of functionals.
  \begin{align}
&\mathcal{U}(w):= \int \left((\varsigma-{2\mu}/{3})|\mm{div}w|^2+
 { \mu} |\mathbb{D} w|^2/2\right)\mm{d}y,\nonumber \\
& \mathcal{I}(w)  :=\left\|\sqrt{P'(\bar{\rho})\bar{\rho}}\left(\frac{gw_3} {P'(\bar{\rho})} -
 \mm{div}w \right)\right\|_0^2 + \vartheta|\nabla_{\mm{h}} w_3|_0^2 , \label{2022307201157} \\
 & \mathcal{E}  :=\|\eta \|_3^2 + \| u \|_{2 }^2  + \| u_t \|_{0}^2 ,\ \mathcal{D}:=\|\eta\|_3^2
 +\| u\|_{3 }^2 +\|u_t\|_{1}^2 .\nonumber  \end{align}

(4) Cut-off functions: Let $\chi(r)\in C_0^\infty(-2,2)$ be a given smooth function such that
  \begin{equation}
0\leqslant \chi(r)\leqslant 1\mbox{ and }\chi(r)=\begin{cases}
  1&\mbox{for }|r|\leqslant 1;
  \\
  0&\mbox{for }|r|\geqslant 2.
  \end{cases}\nonumber
  \end{equation}
Then we further define  the smooth function  $\chi_n(r)\in C_0^\infty(\mathbb{R})$  as follows: for given $n\geqslant 1$,
\begin{equation} \chi_n(r)=\begin{cases}
 1&\mbox{for }|r|\leqslant n-1;
 \\  \chi(n+1-|r|)&\mbox{for } n-1<|r|<n;  \\
 0&\mbox{for }|r|\geqslant n.
\end{cases}
\label{202306saf081959}
\end{equation}
In addition, we define  $\chi_{n,n}(y_1,y_2):=\chi_n(y_1)\chi_n(y_2) $.

Now we state our instability result, which physically shows that the RT instability always occurs in the SCVF mode with non-periodic motion.
\begin{thm}\label{thm:0202}
 We assume that
\begin{enumerate}[(1)]
\item $P_\pm(\tau)\in C^4(\mathbb{R})$ are positive, and strictly increasing with respect to $\tau$;
\item $\mu_\pm>0$ and   $\varsigma_\pm>0$ are constants;
\item $\bar{\rho}_- \in W^{3,\infty}(h_-,0)$, $\bar{\rho}_+\in W^{3,\infty}(0,h_+)$ and they satisfy \eqref{201611051547}--\eqref{201612291257}.
\end{enumerate}
The zero solution is unstable to the RT problem \eqref{n0101nnnM} in   Hadamard sense, that is, there are positive constants $\epsilon$, $c_i$ for $1\leqslant i\leqslant 8$, and functions $\tilde{\eta}^0$, $\tilde{u}^0\in H_0^1( \Omega_{c_1,c_2}\!\!\!\!\!\!\!\!\!\!\!\!\!\!-\quad\ \,) \cap H^3(\Omega_{c_1,c_2} )$,
such that, for any $\delta\in (0,c_3]$ and for any $n\geqslant \max\{\delta^{-2},c_4\}$,
\begin{itemize}
  \item $\|\chi_{n,n}(\tilde{\eta}^0,\tilde{u}^0)/n\|_3\leqslant c_5$,
  \item there exists a function $u^\mm{r}\in H_0^1\cap H^2$ depending on $\delta$ and satisfying
\begin{align*}\| u^{\mm{r}} \|_2\leqslant c_6,
\end{align*}
  \item the RT  problem   with the initial data
\begin{align}
\label{2023306082002}
(\eta^0, u^0):=\delta n^{-1}\chi_{n,n}(\tilde{\eta}^0,\tilde{u}^0)  + \delta^2(0,u^\mm{r})\in (H_0^1\cap H^3)\times (H_0^1\cap H^2)
\end{align}
admits a unique strong solution $(\eta,u)$ belongs to $C^0([0,T],H_{0,*}^{3,1/2} \times H^2)$ and  satisfies
\begin{align}
\label{202307171300}
 \|\omega_{\mm{h}}(T^\delta)\|_{0},\ \|\omega_3(T^\delta)\|_{0},\  |{\omega}_3(T^\delta)|_{0}\geqslant\epsilon
\end{align}
for some escape time $T^\delta:= {c_7}^{-1}\mm{ln}(\epsilon c_8/\delta)\in (0,T)$, where   $\omega=\eta$ or $u$.
\end{itemize}
\end{thm}
\begin{rem}
\label{202230800212201}
It should be noted that $ \Omega_{c_1,c_2}\!\!\!\!\!\!\!\!\!\!\!\!\!\!-\quad\ \,=2\pi c_1\times 2\pi c_2\times (h_-,h_+)$ and $\Omega_{c_1,c_2}$ in Theorem \ref{thm:0202} is defined as $\Omega_{L_1,L_2}$  with $c_i$ in place of $L_i$ for $i=1$, $2$, see \eqref{2022307222447} for the definition of $\Omega_{L_1,L_2}$.  In addition,  let $\varrho=\bar{\rho}J^{-1} $, then $\eta$, $\varrho$ and $u$ are just  strong solutions of the problem \eqref{201611041430M}; moreover $\|\varrho-\bar{\rho}\|_3\lesssim \|\eta\|_2$ and $ \varrho-\bar{\rho} \in C^0([0,T],H^2)$ due to  $ \eta\in C^0([0,T],H_{0,*}^{3,1/2})$.
\end{rem}
\begin{rem}
\label{202300808400}
By Lemma \ref{201809012320}, $\xi:= \eta+y$  satisfies the lower bound \eqref{201803121601xx21082109} and the homeomorphism properties \eqref{201803121601xx2108}--\eqref{2018031adsadfa21601xx} for sufficiently small $\delta$, therefore we can write out the corresponding instability result in Eulerian coordinates, please refer to \cite[Theorem 1.2]{JFJSZWC} or Appendix \ref{secsd:09}; in particular, we have the interface instability by   \eqref{201803121601xx2108},  \eqref{202223020155} and \eqref{2020307saf302204}:
\begin{align}|d(T^\delta)|_{0}=  \int_{\mathbb{R}^2} |\eta_3(y_{\mm{h}},0,T^\delta)|^2\det \nabla_{y_{\mm{h} }} \zeta_{\mm{h}}(y_{\mm{h}},0,T^\delta) |\mm{d}y_{\mm{h}} \geqslant   |\eta_3( T^\delta)|_0^2  /2 \geqslant\epsilon/2,
\end{align}
which presents that the RT  equilibrium solution $(\bar{\rho},\mathbf{0})$ with $d=0$ is unstable to the SCVF model with non-periodic motion.
\end{rem}

The proof of Theorem \ref{thm:0202} is based on a so-called bootstrap instability method, which has its origin in \cite{GYSWIC,GYSWICNonlinea}. Later, various versions of bootstrap approaches were presented by many authors, see \cite{FSSWVMNA,GYHCSDDC} for examples. In the spirit of bootstrap instability method  in \cite[Lemma 1.1]{GYHCSDDC}, the proof procedure for our problem will be divided into the following five steps.

Firstly, we will construct unstable solution $e^{c_7 t}(\tilde{\eta}^0,\tilde{u}^0)$ for the linearized periodic  RT problem, see Proposition \ref{growingmodesolneriodic} in Section \ref{2022306010101651}.  Secondly, we will establish a Gronwall-type energy inequality for the solutions of the (non-periodic) RT problem \eqref{n0101nnnM} as in \citep[Lemma 4.10]{JFJSZWC}, see Proposition \ref{pro:0401nd} in Section \ref{2022306101653}. Thirdly,  we want to use $\delta e^{c_7 t}\chi_{n,n}(\tilde{\eta}^0,\tilde{u}^0)/n$ to be the approximate solution of the  RT problem,  therefore we will exploit  the stratified elliptic regularity theory to modify the initial data of the approximate solution $\delta\chi_{n,n} (\tilde{\eta}^0,\tilde{u}^0)/n$ as in   \citep[Proposition 8]{JFJSNS}, so that the obtained modified initial data satisfy the initial compatibility jump condition of the RT  problem, i.e.,
\begin{align}
\label{2022212081221}
 \llbracket P(\bar{\rho}J^{-1})\mathbb{I}-\mathbb{S}_{\mathcal{A}}( u)    \rrbracket  \mathcal{A}\mathbf{e}^3 = \vartheta \mathcal{H} \mathcal{A}\mathbf{e}^3 \mbox{ on }\Sigma\mbox{ for }t=0,
 \end{align}
 and at the same time, are close to $\delta\chi_{n,n} (\tilde{\eta}^0,\tilde{u}^0)/n$, see Proposition \ref{lem:modfied} in Section \ref{2022306101652}. Fourthly, we deduce the error estimates between the approximate solutions  and the solutions of the RT problem as in \cite[Lemma 4.1]{JFJSWWWN}, see Proposition \ref{lem:0401} in Section \ref{20230721}.
Finally, we show the existence of escape times as in \cite[Lemma 1.1]{GYHCSDDC}  and thus obtain Theorem \ref{thm:0202} in Section \ref{sec:030845}.

Now we further mention the first step, which includes new ideas.  Since it seems to be  difficulty to directly construct the unstable solutions for the linearized (non-periodic) RT problem, we first will construct the unstable solution, denoted by $e^{c_7 t}(\tilde{\eta}^0,\tilde{u}^0)$, for the linearized periodic RT  problem by following the argument of \citep[Theorem 2.2]{GYTI2}, see the first two conclusions  in Proposition \ref{growingmodesolneriodic}. Obviously $\delta e^{c_7 t}(\tilde{\eta}^0,\tilde{u}^0)$ are also the   linearized periodic RT  problem, but can not directly used to be the approximate solutions of the nonlinear RT problem \eqref{linearizedxx} due to the periodicity. Therefore we will use a  cut-off function $\chi_{n,n}$ to cut off the obtained  linear unstable periodic solutions. Thanks  to the structure of variable separation with respective to $y_{\mm{h}}$ and $y_3$ of the linear solution  $(\tilde{\eta}^0,\tilde{u}^0)$, we can derive   the
 uniform cut-off estimates with respect to $n$ for $\chi_{n,n}(\tilde{\eta}^0,\tilde{u}^0)/n $   (see the third conclusion in Proposition \ref{growingmodesolneriodic}) and thus excitedly find that  $\delta e^{c_7  t}\chi_{n,n}(\tilde{\eta}^0,\tilde{u}^0)/n$ can be regarded as the approximate solution  of the nonlinear unstable  solution if $n$ is sufficiently large for given $\delta$. We mention that our proof for Theorem \ref{thm:0202} can be also applied to the corresponding incompressible case, which will be recorded in a forthcoming paper.  In addition, our approximate method can be also applied to further construct to other instability solutions with non-periodic motion, such as Parker instability \cite{JFJSSETEFP} and  thermal instability \cite{JFJSOUI}   in a slab domain. In addition,  Wang verified that the sufficiently strong non-horizontal magnetic fields can inhibit RT instability in stratified incompressible viscous non-resistive  magnetohydrodynamic fluids in a slab domain $\mathbb{R}^2\times (h_-,h_+)$ \cite{WYJ2019ARMA}. Applying  our approximate method, we can further prove that the magnetic fields can not inhibit RT instability in a slab domain, if the field intensity is too weak.

Finally we remark that there exist a huge of existence results of instability solutions,
 which are  periodic with respect to the horizontal spacial variables,  for  the problems of  flow instabilities, see the solutions of RT instability  in \cite{HHJGY,JFJSZWC},  the solutions  of Parker instability in \cite{JFJSSETEFP} and the solutions of thermal convection instability in \cite{GYYQH,JFJSOUI} for examples.
  However, by further using our new approximate method, all results of (nonlinear) instability solutions with periodic motion aforementioned can be extended to the corresponding non-periodic cases. In addition, it is pointed out that Wang verified the fact that the sufficiently strong non-horizontal magnetic fields can inhibit RT instability in stratified incompressible viscous non-resistive  magnetohydrodynamic fluids in a slab domain $\mathbb{R}^2\times (h_-,h_+)$ \cite{WYJ2019ARMA}. However, applying  our approximate method, we can further prove that the magnetic fields can not inhibit RT instability in a slab domain, if the field intensity is too weak.

 \section{Linear instability}\label{2022306010101651}

In this section, we wish to construct an unstable solution,  which has growing $H^4$-norm, to the linearized periodic RT  problem:
\begin{equation}\label{linearized} \begin{cases}
\eta_t=u &\mbox{in } \Omega_{L_1,L_2},\\[1mm]
\bar{\rho} u_t = g\bar{\rho}(\mm{div}\eta \mathbf{e}^3- \nabla \eta_3)+\mm{div}   \Upsilon(\eta,u)  &\mbox{in }\Omega_{L_1,L_2},\\[1mm]
\llbracket   u  \rrbracket  =\llbracket   \eta  \rrbracket  =0,\  \llbracket     \Upsilon(\eta,u)\mathbf{e}^3\rrbracket   =- \vartheta \Delta_{\mm{h}}\eta \mathbf{e}^3&\mbox{on }\Sigma_{L_1,L_2},\\
(\eta, u)=0 &\mbox{on }\partial\Omega\!\!\!\!\!- ,\\
(\eta,  u)|_{t=0}= (\eta^0,  u^0 ) &\mbox{in }\Omega_{L_1,L_2},
\end{cases}
\end{equation}
see \eqref{2022307222447} for the definitions of $\Omega_{L_1,L_2}$  and $\Sigma_{L_1,L_2}$.
We will
construct such growing solution via some synthesis as in  \cite{GYTI2} by first constructing
a growing mode for any but fixed spatial frequency.  We mention that the construction of the linear instability solutions follows Guo--Tice's method in \cite{GYTI2}, however some new conclusions, which are useful in the construction of nonlinear (non-periodic) instability solutions,  will be further  supplemented.
\subsection{Growing modes}
To start with, we make a growing mode ansatz of solutions, i.e., for some $\lambda>0$,
\begin{equation*}
\eta({y},t)=\sigma ( {y})e^{\lambda t}\mbox{ and }
u( {y},t)=w(y)e^{\lambda t} .
\end{equation*}
Substituting this ansatz into \eqref{linearized}, and then eliminating
$\sigma$ by using the first  equation, we arrive at the
boundary value problem
$w$:
 \begin{equation}
 \label{201604061413}        \begin{cases}
\lambda^2 \bar{\rho}w= g\bar{\rho}(\mm{div}w \mathbf{e}^3-\nabla w_3 )+\mm{div}   \Upsilon(w,\lambda w)
&\mbox{ in } \Omega_{L_1,L_2},\\
  \llbracket  w  \rrbracket =0,\  \llbracket   \Upsilon(w,\lambda w)\mathbf{e}^3 \rrbracket=- \vartheta \Delta_{\mm{h}} w_3 \mathbf{e}^3&\mbox{ on }\Sigma_{L_1,L_2} ,
\\ {w} =0  &\mbox{ on }\partial\Omega\!\!\!\!\!- . \end{cases}
\end{equation}

 Now we formally define the new unknowns $(\varphi, \theta, \psi):(h_-,h_+) \rightarrow \mathbb{R}^3$ according to
\begin{align}
 &w_1(y)  =  -\mm{i} \varphi(y_3) e^{\mm{i} \xi\cdot y_{\mm{h}} },\  w_2(y) =    - \mm{i} \theta(y_3) e^{\mm{i}  \xi\cdot y_{\mm{h}} }\mbox{ and }  w_3(y) = \psi(y_3)e^{\mm{i} \xi\cdot y_{\mm{h}} }, \nonumber
\end{align}where $\xi \cdot y_{\mm{h}} = \xi_1 y_1+ \xi_2 y_2 $ for $\xi \in L^{-1}_1 \mathbb{Z}
\times L^{-1}_2 \mathbb{Z} $.
It is easy to check that \cite{GYTI2}
\begin{equation}\nonumber
 \mm{div}{w} = (\xi_1 \varphi + \xi_2 \theta + \psi')e^{\mm{i} x_{\mm{h}}\cdot \xi}
\end{equation}
and
\begin{equation}\nonumber
 \mathbb{D} w = \begin{pmatrix}
       2\xi_1 \varphi & \xi_1 \theta + \xi_2 \varphi & \mm{i}(\xi_1 \psi - \varphi' ) \\
       \xi_1 \theta + \xi_2 \varphi & 2\xi_2 \theta & \mm{i}(\xi_2 \psi - \theta')  \\
      \mm{i}(\xi_1 \psi - \varphi' ) & \mm{i}(\xi_2 \psi - \theta') & 2\psi'
      \end{pmatrix} e^{\mm{i} x_{\mm{h}}\cdot \xi}.
\end{equation}

For each fixed $\xi$, and for the new unknowns $\varphi(x_3), \theta(x_3), \psi(x_3)$ and $\lambda$ we arrive at the following system of ODEs:
\begin{align}
&  \left( \lambda \left( \lambda   \bar{\rho}  +   \mu \abs{\xi}^2 + \xi_1^2 \left( \varsigma +  \mu/3\right)   \right)  +\xi_1^2 P'(\bar{\rho}) \bar{\rho} \right) {\varphi} - \lambda \mu \varphi''  \nonumber \\
&=  \xi_1 \left( g\bar{\rho}  {\psi}-\left( \lambda \left(  \varsigma +  \mu/3\right)    +
P' (\bar{\rho})\bar{\rho}\right) {\psi}'  \right)
 -\xi_1 \xi_2 (\lambda \left(  \varsigma +  \mu/3 \right) + P'(\bar{\rho}) \bar{\rho} ) {\theta} , \label{2022305251953}
\end{align}
\begin{align}
&\left( \lambda \left( \lambda \bar{\rho}  + \mu \abs{\xi}^2 + \xi_2^2 \left(   \varsigma +   \mu/3  \right)   \right)   + \xi_2^2 P'(\bar{\rho}) \bar{\rho}  \right) {\theta} - \lambda\mu \theta'' \nonumber  \\
&=  \xi_2\left(g\bar{\rho} {\psi}-(\lambda\left(\varsigma + \mu/3\right)+ P'(\bar{\rho})\bar{\rho}) {\psi}' \right)
 -\xi_1 \xi_2  \left( \lambda(\varsigma +  \mu/3 )  +  P'(\bar{\rho}) \bar{\rho} \right) {\varphi} \label{20223052519531}
\end{align}
and
\begin{align}\label{w3equation}
 &( \lambda^2  \bar{\rho}+\lambda \mu \abs{\xi}^2   ){\psi} -\left(\left( \lambda\left(  4  \mu/3 +  \varsigma\right) +  P'(\bar{\rho}) \bar{\rho} \right) {\psi}'\right)'\nonumber \\
&= \left(\left(\lambda \left( \varsigma +  \mu/3 \right)  +   P'(\bar{\rho}) \bar{\rho} \right)\left( \xi_1  {\varphi} + \xi_2  {\theta} \right)\right)' +  g \bar{\rho} (\xi_1 \varphi  + \xi_2  \theta  ).
\end{align}

The first jump condition in \eqref{201604061413}$_2$  yields jump conditions for the new unknowns:
\begin{equation}
\label{202306141513}
  \llbracket{\varphi}  \rrbracket = \llbracket{\theta}  \rrbracket =  \llbracket{\psi} \rrbracket = 0.
\end{equation}
The second jump condition in \eqref{201604061413}$_2$  becomes
\begin{align}
&  \llbracket  \left(\lambda\left(\varsigma - 2 \mu /3   \right)
+   P'(\bar{\rho}) \bar{\rho}  \right)(\xi_1  {\varphi} +  \xi_2  {\theta} +{ \psi}')  \mathbf{e}^3 \nonumber \\
& + \lambda\mu(  \mm{i} (\xi_1 \psi-\varphi')\mathbf{e}^1+ \mm{i} (\xi_2 \psi-\theta')\mathbf{e}^2 +
  2 \psi'\mathbf{e}^3 ) \rrbracket
    = \vartheta \abs{\xi}^2  {\psi} \mathbf{e}^3, \nonumber
\end{align}
 which implies that
\begin{equation}
\llbracket {  \mu (\varphi'-\xi_1  \psi)} \rrbracket = \llbracket { \mu (\theta'-\xi_2 \psi)} \rrbracket =0 \label{2022306141533}
\end{equation}
and that
\begin{align}
&\llbracket\left( {\lambda(  \varsigma+   \mu/3  ) }
+  P'(\bar{\rho})\bar{\rho} \right) ( \xi_1  {\varphi} + \xi_2  {\theta}+{\psi}'  ) \nonumber  \\
&+ { \lambda \mu \left(\psi' - \xi_1 \varphi - \xi_2 \theta    \right)  }\rrbracket
  = \vartheta \abs{\xi}^2  {\psi}. \label{2023112233}
\end{align}
By \eqref{201604061413}$_3$, the boundary conditions
\begin{equation}
 \varphi(h_-) =  \varphi(h_+) = \theta(h_-) =  \theta(h_+)  =\psi(h_-) =  \psi(h_+)=0
 \label{202305282005}
\end{equation}
must also hold.
We will prove the the existence of real-value solutions $\varphi$, $\theta$, $\psi$ to the boundary value problem \eqref{2022305251953}--\eqref{202305282005} in next section, and thus further get real-value solutions $w$ and $\lambda$ of the original problem \eqref{201604061413}.

Next we provide the energy structure for  the boundary value problem of \eqref{2022305251953}--\eqref{202305282005}.
Taking the inner product of \eqref{2022305251953}, \eqref{20223052519531}, resp.  \eqref{w3equation} and  $\varphi$, $\theta$, resp. $\psi$ in $L^2(h_-,h_+)$, adding the three resulting identities together, and then making use of the integration by parts and the boundary conditions \eqref{2022306141533}--\eqref{202305282005}, we derive the following  energy identity for  \eqref{2022305251953}--\eqref{202305282005}:
\begin{equation}
\label{202305291431xx}
  \lambda^2  {\mathcal{J}(\varphi,\theta,\psi)} = F(\varphi,\theta,\psi):= {E(\varphi,\theta,\psi)- \lambda D(\varphi,\theta,\psi)},
\end{equation}
where we have defined that
\begin{align}
 &\mathcal{J}(\varphi,\theta, \psi): = \int_{h_-}^{h_+}\bar{\rho}  (|\varphi|^2 + |\theta|^2+ |\psi|^2)\mm{d}y_3,  \label{202230614sfd2048} \\
&  E(\varphi,\theta,\psi):= \int_{h_-}^{h_+}
 \left(  2 g\bar{\rho}(\xi_1\varphi+ \xi_2\theta ) \psi-P'(\bar{\rho})\bar{\rho}(\xi_1\varphi+\xi_2\theta+\psi' )^2 \right)\mm{d}y_3- {\vartheta |\xi|^2}{\psi}^2(0)   \label{2022306142048}
 \end{align}
 and
\begin{align}
 D(\varphi,\theta,\psi):=&
\int_{h_-}^{h_+} \varsigma | \xi_1\varphi+ \xi_2\theta+ \psi'  |^2\mm{d}y_3 + \int_{h_-}^{h+} \mu (|\xi_2\varphi-\xi_1\theta|^2+|\varphi'-\xi_1\psi|^2\nonumber  \\
&+|\theta'-\xi_2\psi|^2+ |\xi_1\varphi+\xi_2\theta- \psi'|^2+|
\xi_1\varphi+\xi_2\theta+ \psi'|^2/3) \mm{d}y_3.
\label{202202307312119}
\end{align}

In addition, by  a simple computation and the hydrostatic relations in
 \eqref{201611051547}$_1$, we have
\begin{align*}
&2 g\bar{\rho}(\xi_1\varphi+ \xi_2\theta ) \psi-P'(\bar{\rho})\bar{\rho}\left(\xi_1\varphi+\xi_2\theta+\psi' \right)^2 \\
&=\frac{g^2{\bar{\rho}}}{{P'(\bar{\rho}}) }
 \psi^2 -2 g\bar{\rho}\psi\psi'- P'(\bar{\rho})\bar{\rho}\left(\frac{g }{ {P'(\bar{\rho}}) }
\psi-\xi_1\varphi-\xi_2\theta-\psi' \right)^2
 \end{align*}
 and
$$\int_{h_-}^{h_+}\left(\frac{g^2{\bar{\rho}}}{{P'(\bar{\rho}}) }
 \psi^2 -2 g\bar{\rho}\psi\psi'\right)\mm{d}y_3
=    g\llbracket \bar{\rho}\rrbracket  {\psi}^2(0). $$
Thanks to the above two identities, $E$ can be rewritten as follows:
\begin{align}
E(\varphi,\theta,\psi)=&(g\llbracket \bar{\rho}\rrbracket - {\vartheta |\xi|^2})|{\psi}(0)|^2 -\int_{h_-}^{h_+}P'(\bar{\rho})\bar{\rho}\left|\frac{g }{ {P'(\bar{\rho}}) }
\psi-\xi_1\varphi-\xi_2\theta-\psi' \right|^2\mm{d}y_3.\label{2022306031524}
\end{align}
\emph{In Section \ref{20230721}, we will extend both the expressions of \eqref{202230614sfd2048},   \eqref{202202307312119} and \eqref{2022306031524}  to the case that $\varphi$, $\theta$ and  $\psi$ are complex functions. Under such case, the notation $|\cdot| $ represents the modulus of a complex function.}

\subsection{Solutions to \eqref{2022305251953}--\eqref{202305282005}   via modified variational methods}\label{mod_section}
To obtain the solution of \eqref{2022305251953}--\eqref{202305282005}, we will first consider  the following initial value problem, which is obtained from \eqref{2022305251953}--\eqref{202305282005} by modifying $\lambda$:
\begin{equation}
\begin{cases}\left(\alpha(s)   \bar{\rho}  +   s(\mu \abs{\xi}^2 + \xi_1^2 \left( \varsigma +  \mu/3\right) )    +\xi_1^2 P'(\bar{\rho}) \bar{\rho} \right) {\varphi}-
s \mu \varphi''   \\
= \xi_1 \left(g\bar{\rho}  {\psi} -\left( s \left(  \varsigma +  \mu/3\right)    +
P' (\bar{\rho})\bar{\rho}\right) {\psi}' \right)
 -\xi_1 \xi_2 (s \left(  \varsigma +  \mu/3 \right) + P'(\bar{\rho}) \bar{\rho} ) {\theta} , \\
\left(   \alpha(s)  \bar{\rho}  + s(\mu \abs{\xi}^2 + \xi_2^2 \left(   \varsigma +   \mu/3  \right)  )   + \xi_2^2 P'(\bar{\rho}) \bar{\rho}  \right) {\theta}  - s\mu \theta''  \\
=  \xi_2\left(g\bar{\rho} {\psi} -(s\left(\varsigma + \mu/3\right)+ P'(\bar{\rho})\bar{\rho}) {\psi}' \right)
 -\xi_1 \xi_2  \left( s(\varsigma +  \mu/3 )  +  P'(\bar{\rho}) \bar{\rho} \right) {\varphi} , \\
(\alpha(s)   \bar{\rho}+s\mu \abs{\xi}^2   ){\psi} -\left(\left( s\left(  4  \mu/3 +  \varsigma\right) +  P'(\bar{\rho}) \bar{\rho} \right) {\psi}'\right)'  \\
= \left(\left(s \left( \varsigma +  \mu/3 \right)  +   P'(\bar{\rho}) \bar{\rho} \right)\left( \xi_1  {\varphi} + \xi_2  {\theta} \right)\right)' +  g \bar{\rho} (\xi_1 \varphi  + \xi_2  \theta  ),\\
  \llbracket{\varphi}  \rrbracket = \llbracket{\theta}  \rrbracket =  \llbracket{\psi} \rrbracket = 0,
  \\
\llbracket {  \mu (\varphi'-\xi_1  \psi)} \rrbracket = \llbracket { \mu (\theta'-\xi_2 \psi)} \rrbracket =0,\\
\llbracket\left( {s(  \varsigma+   \mu/3  ) }
+  P'(\bar{\rho})\bar{\rho} \right) ( \xi_1  {\varphi} + \xi_2  {\theta}+{\psi}'  )
+ { s \mu \left(\psi' - \xi_1 \varphi - \xi_2 \theta    \right)  }\rrbracket
  = \vartheta \abs{\xi}^2  {\psi},\\
 \varphi(h_-) =  \varphi(h_+) = \theta(h_-) =  \theta(h_+)  =\psi(h_-) =  \psi(h_+)=0.
 \end{cases}\label{202305291213}
\end{equation}
In view of \eqref{202305291431xx}, the above modified problem
   has the following  variational structure:
\begin{equation}
\label{202305291431}
\alpha(s):=\sup_{\phi,\chi,\omega\in H_0^1(h_-,h_+)}\frac{F(\phi,\chi,\omega;s) }{\mathcal{J}(\phi,\chi,\omega)},
\end{equation}
see the  definition of  $F(\varphi,\theta,\psi;s)$ in \eqref{202305291431xx}.  Since both $F$ and $\mathcal{J}$ are homogeneous of degree $2$,   it holds that
\begin{equation}
\label{2023052safda91431}
\alpha(s)=\sup_{(\phi,\chi,\omega )\in \mathcal{A}}F(\phi,\chi,\omega;s),
\end{equation}
where $\mathcal{A}:=\{(\phi,\chi,\omega)\in (H^1_0(h_-,h_+))^3~|~\mathcal{J}(\phi,\chi,\omega)=1\}$.
Thus we easily  construct solutions to the modified problem \eqref{202305291213} by  utilizing a standard variational method.

\begin{lem}\label{minexist}
Under the assumptions of Theorem \ref{thm:0202}, for any given $s\in \mathbb{R}^+$ and any given $\xi\in \mathbb{R}^2$, $F(\phi,\chi,\omega;s)$ has a maximizer $( {\varphi}, {\theta}, {\psi})$ in  $\mathcal{A}$. Moreover, the maximizer $( {\varphi}, {\theta}, {\psi})$ belongs to $H^4((h_-,0)\cup (0,h_+))$ and satisfies \eqref{202305291213} with $\alpha(s) = F( {\varphi}, {\theta}, {\psi}; s)$.
\end{lem}
\begin{pf} In what follows, we denote $F(\cdot,\cdot,\cdot;s)$ by $F(\cdot,\cdot,\cdot)$  for the sake of the simplicity.

(1) Recalling the definition of $F(\varphi,\theta,\psi)$, we can derive that\begin{align}
F(\varphi,\theta,\psi) =& g  \int_{h_-}^{h_+} \bar{\rho}\left(\xi_1(\varphi ^2 + \psi^2)+\xi_2 (\theta ^2 + \psi^2)  \right) \mm{d}y_3-  {\vartheta |{\xi}|^2} \psi^2(0)
 -sD(\varphi,\theta,\psi )\nonumber \\
&-   \int_{h_-}^{h_+}(  P'(\bar{\rho}) \bar{\rho}  (\xi_1 \varphi+  \xi_2 \theta +  \psi')^2
+  g \bar{\rho}(\xi_1( \varphi-\psi)^2+ \xi_2( \theta-\psi)^2))\mm{d}y_3
\nonumber \\
  \leqslant & g(|{\xi}_1|+ |\xi_2|),\label{Elowerbound}
\end{align} which shows that $F$ is bounded above on $\mathcal{A}$. Therefore there exists a maximizing sequence $\{(\varphi_n,\theta_n,\psi_n)\}_{n=1}^\infty\subset\mathcal{A}$. Obviously $\{\psi_n(0)\}_{n=1}^\infty$ is bounded in $\mathbb{R}$, and $\{\varphi_n\}_{n=1}^\infty$, $\{\theta_n\}_{n=1}^\infty$ and $\{\psi_n\}_{n=1}^\infty$ are bounded in $H^1_0(h_-,h_+)$, so up to the extraction of a subsequence   $(\varphi_n,\theta_n, \psi_n) \rightharpoonup (\varphi,\theta, \psi)$ weakly in $(H_0^1(h_-,h_+))^3$, and $(\varphi_n,\theta_n,\psi_n) \rightarrow (\varphi,\theta,\psi)$ strongly in $(L^2(h_-,h_+))^3$. In addition, the compact embedding $H_0^1(h_-,h_+)  \hookrightarrow C^0(h_-,h_+)$ implies that $\psi_n(0) \rightarrow \psi(0)$ as well (referring to section 1.3.5.8 in \cite{NASII04}).   Because of the quadratic structure of all the terms in the integrals defining $F(\varphi,\theta,\psi)  $, weak lower semi-continuity (referring to section 1.4.2.7 in \cite{NASII04}) and strong $L^2$-convergence imply that
\begin{align}
 \alpha(s)\geqslant F( {\varphi}, {\theta}, {\psi})
 & \geqslant \lim_{n\rightarrow \infty} \left(\int_{h_-}^{h_+}
  2 g\bar{\rho}(\xi_1\varphi_n+ \xi_2\theta_n ) \psi_n\mm{d}y_3- {\vartheta |\xi|^2}{\psi}^2_n(0)\right)   \nonumber \\
& \quad -\limsup_{n\rightarrow \infty} \left(\int_{h_-}^{h_+}  P'(\bar{\rho})\bar{\rho}(\xi_1\varphi_n+\xi_2\theta_n+\psi'_n )^2\mm{d}y_3 + s D(\varphi_n,\theta_n,\psi_n)\right)\nonumber \\
 &\geqslant \limsup_{n\rightarrow \infty} F(\varphi_n,\theta_n,\psi_n) = \alpha(s). \nonumber
\end{align}
That $( {\varphi}, {\theta}, {\psi})\in\mathcal{A}$ follows from the strong $L^2$-convergence. Hence $F(\varphi,\theta,\psi )$ has a maximizer $( {\varphi}, {\theta}, {\psi})$ in  $\mathcal{A}$.

(2) Next we further show that the maximizer $(\varphi,\theta,\psi )$ constructed above satisfies Euler--Langrange equations, which are equivalent to \eqref{202305291213}.
Let $\tilde{\varphi}$, $\tilde{\theta}$, $\tilde{\psi}\in H_0^1(h_-,h_+)$ be given. Then we define
\begin{equation}
f(t):=F(\varphi+t\tilde{\varphi}  ,\theta+t\tilde{\theta} ,\psi + t\tilde{\psi} )-\alpha \mathcal{J}(\varphi+t\tilde{\varphi},\theta+t\tilde{\theta},\psi+t\tilde{\psi}).
\end{equation}
It is easy to see that  $f\in C^\infty(\mathbb{R})$, $f(t)\geqslant 0$ for any $t\in \mathbb{R}$ and $f(0)=0$. Therefore $f'(0)= 0$, which implies that
\begin{align}
  &\int_{h_-}^{h_+}
 \left(   g\bar{\rho}((\xi_1\varphi + \xi_2\theta ) \tilde{\psi}+(\xi_1\tilde{\varphi} + \xi_2\tilde{\theta }) \psi)-P'(\bar{\rho})\bar{\rho}\left(\xi_1\varphi+\xi_2\theta+\psi' \right) \left(\xi_1\tilde{\varphi}+\xi_2\tilde{\theta}+\tilde{\psi}' \right)  \right)\mm{d}y_3 \nonumber \\
  &- {\vartheta |\xi|^2}{\psi}(0)\tilde{\psi}(0)-s
\int_{h_-}^{h_+} \varsigma ( \xi_1\varphi+ \xi_2\theta+ \psi'  ) ( \xi_1\tilde{\varphi}+ \xi_2\tilde{\theta}+ \tilde{\psi}'  )\mm{d}y_3\nonumber \\
& -s \int_{h_-}^{h+} \mu ( (\xi_2\varphi-\xi_1\theta)(\xi_2\tilde{\varphi}-\xi_1\tilde{\theta})+ (\varphi'-\xi_1\psi)(\tilde{\varphi}'-\xi_1\tilde{\psi}) +(\theta'-\xi_2\psi)(\tilde{\theta}'-
\xi_2\tilde{\psi}) \nonumber \\
&+ (\xi_1\varphi+\xi_2\theta- \psi')(\xi_1\tilde{\varphi}+\xi_2\tilde{\theta}- \tilde{\psi}')+ (
\xi_1\varphi+\xi_2\theta+ \psi') (
\xi_1\tilde{\varphi}+\xi_2\tilde{\theta}+ \tilde{\psi}') /3)\mm{d}y_3 \nonumber \\
 &= \alpha \int_{h_-}^{h_+}\bar{\rho}(\varphi \tilde{\varphi} + \theta \tilde{\theta}+ \psi\tilde{\psi})\mm{d}y_3 .
\end{align}

Since $\tilde{\varphi}$, $\tilde{\theta}$ and $\tilde{\psi}$ are independent, the above identity gives rise to the three weak forms:
\begin{align}\label{ele1}
 &\int_{h_-}^{h_+}
 \left(   g\xi_1\bar{\rho} \psi - \xi_1 P'(\bar{\rho})\bar{\rho}\left(\xi_1\varphi+\xi_2\theta+\psi' \right) -s\varsigma \xi_1 ( \xi_1\varphi+ \xi_2\theta+ \psi'  )\right)\tilde{\varphi} \mm{d}y_3\nonumber \\
  & - s  \int_{h_-}^{h+} \mu \left(\frac{1}{3}\left((\xi_1^2+3|\xi|^2)\varphi+\xi_1\xi_2\theta-2\xi_1 \psi'\right) \tilde{\varphi}  +(\varphi'-\xi_1\psi)\tilde{\varphi}' \right)\mm{d}y_3  = \alpha \int_{h_-}^{h_+}\bar{\rho}  \varphi  \tilde{\varphi} \mm{d}y_3 ,
 \end{align}
\begin{align}\label{ele2}
 &\int_{h_-}^{h_+}
 \left(   g\xi_2\bar{\rho} \theta -\xi_2 P'(\bar{\rho})\bar{\rho}\left(\xi_1\varphi+\xi_2\theta+\psi' \right)-s
  \varsigma\xi_2 ( \xi_1\varphi+ \xi_2\theta+ \psi'  )\right)  \tilde{\theta} \mm{d}y_3\nonumber \\
 &-s  \int_{h_-}^{h+} \mu\left(\frac{1}{3} \left( \xi_1\xi_2\varphi+(\xi_2^2+3|\xi|^2)\theta -2\xi_2  \psi'\right)\tilde{\theta}   +(\theta'-\xi_2\psi) \tilde{\theta}'\right) \mm{d}y_3= \alpha \int_{h_-}^{h_+}\bar{\rho}  \theta  \tilde{\theta}\mm{d}y_3
\end{align}
and
\begin{align}\label{ele2x}
 &\int_{h_-}^{h_+}
 \left(   g\bar{\rho}(\xi_1\varphi + \xi_2\theta ) \tilde{\psi} -P'(\bar{\rho})\bar{\rho}\left(\xi_1\varphi+\xi_2\theta+\psi' \right) \tilde{\psi}' \right) \mm{d}y_3- {\vartheta |\xi|^2}{\psi}(0)\tilde{\psi}(0) \nonumber \\
  &+s\int_{h_-}^{h+} \mu \left(   \left( \left(\xi_1\varphi   + \xi_2\theta\right)'-|\xi|^2\psi\right)
\tilde{\psi} +(2\left(\xi_1\varphi  + \xi_2\theta\right)
\tilde{\psi}' - 4\psi' \tilde{\psi}' )/3 \right) \mm{d}y_3 \nonumber \\
&-s
\int_{h_-}^{h_+} \varsigma( \xi_1\varphi+ \xi_2\theta+ \psi'  ) \tilde{\psi}'\mm{d}y_3  = \alpha \int_{h_-}^{h_+}\bar{\rho}   \psi  \tilde{\psi} \mm{d}y_3 .
\end{align}

By making variations with $\tilde{\varphi}$, $\tilde{\theta}$ and $\tilde{\psi}$ compactly supported in either $(h_-,0)$ or $(0,h_+)$, we find that $\varphi$, $\theta$ and $\psi$ satisfy the equations \eqref{202305291213}$_1$--\eqref{202305291213}$_3$ in a weak sense in $(h_-,0)$ and $(0,h_+)$.  Standard bootstrapping arguments then show that $(\varphi,\theta,\psi)$ are in $H^4((h_-,0)\cup(0,h_+))$.  This implies that the equations are also  satisfied on $(h_-,0)$ and $(0,h_+)$.  Since $(\varphi,\theta,\psi)\in (H^4((h_-,0)\cup(0,h_+)))^3$, the traces of the functions and their derivatives are well-defined at the endpoints $y_3 = 0$, $h_\pm$.

To show that the jump conditions are satisfied, we make variations with respect to arbitrary $\tilde{\varphi}$, $\tilde{\theta}$,  $\tilde{\psi} \in C_0^\infty(h_-,h_+)$.  Integrating  the terms in \eqref{ele1} with derivatives of $\tilde{\varphi}$ by parts and using that $\varphi$ solves \eqref{202305291213}$_1$ on $(h_-,0)\cup(0,h_+)$, we find that
\begin{equation}
\llbracket \mu(\varphi'-{\xi}_1 \psi)\rrbracket \tilde{\varphi}(0) = 0.
\end{equation}
Since $\tilde{\varphi}(0)$ may be chosen arbitrarily, we deduce $
\llbracket {  \mu (\varphi'-\xi_1  \psi)} \rrbracket =  0$ in \eqref{202305291213}$_5$. Similarly, we also have $\llbracket { \mu (\theta'-\xi_2 \psi)} \rrbracket =0$. Therefore \eqref{202305291213}$_5$ holds. In addition, performing a similar integration by parts in \eqref{ele2x} yields the jump condition \eqref{202305291213}$_6$. Finally, the conditions \eqref{202305291213}$_4$ and \eqref{202305291213}$_7$ are satisfied trivially since $\varphi$, $\theta$, $\psi \in H_0^1(h_-,h_+) \hookrightarrow C_0^{0,1/2}[h_-,h_+] $. \hfill$\Box$
\end{pf}

The  next lemma establishes the behaviors of  $\alpha(s)$ with respect to $s$.
\begin{lem}\label{eigenlip}
Under the assumptions of Theorem \ref{thm:0202}, for any given  $\xi\in \mathbb{R}^2$,  let $\alpha$: $\mathbb{R}^{+}  \rightarrow \mathbb{R}$ be given by \eqref{2023052safda91431}.
Then the following assertions hold.
\begin{enumerate}[(1)]
 \item$ \alpha\in C^{0,1}_{\mm{loc}}(\mathbb{R}^+)$, and in particular $\alpha\in C^{0}(\mathbb{R}^+)$.
\item There exists a positive constant $b_1= b_1(\xi,\mu ,\varsigma , h_\pm)$ so that
\begin{equation}\label{el00}
\alpha(s) \leqslant  g(|{\xi}_1|+|\xi_2|) - s b_1.
\end{equation}
\item $\alpha(s)$ is strictly decreasing.
\item Let
$$|{\xi}|_{\mm{c}} =
\begin{cases}
\sqrt{g \llbracket{\bar{\rho}}\rrbracket/\vartheta}&\mbox{if }\vartheta>0;
\\
\infty&\mbox{if }\vartheta=0.
\end{cases}$$
If $\xi$ satisfies
\begin{equation}
\label{202307212216}
|{\xi}| \in (0,   |{\xi}|_{\mm{c}}),
\end{equation}
then there exists a $s_0$, depending on $\xi$, $g$, $\vartheta$, $\mu$, $\rho $,  $h_\pm$ and $P_\pm$, such that $\alpha(s)>0$ for any $s\in (0,s_0)$.
\end{enumerate}
\end{lem}
\begin{pf}
(1) Let $Q = [a,b]  \in \mathbb{R}^+$ be  a compact interval. By Lemma \ref{minexist}, for each $s \in \mathbb{R}^+$, there exists   $(\varphi_s, \theta_s, \psi_s)\in\mathcal{A}$ so that
\begin{equation}
 F(\varphi_s,\theta_s, \psi_s; s) = \sup_{(\phi,\chi,\omega)\in \mathcal{A}} F(\phi,\chi,\omega; s) = \alpha(s).
\end{equation}
We deduce from the non-negativity of $D$, the maximality of $(\varphi_s,\theta_s, \psi_s)$, and the equality in \eqref{Elowerbound} that
\begin{equation}
 F(\tilde{\varphi},\tilde{\theta},\tilde{\psi};b) \leqslant  F(\tilde{\varphi},\tilde{\theta},\tilde{\psi};s) \leqslant   F(\varphi_s, \theta_s,\psi_s;s) \leqslant g(|\xi_1|+|\xi_2|) - s D(\varphi_s, \theta_s, \psi_s)
\end{equation}
for all $ (\tilde{\varphi},\tilde{\theta},\tilde{\psi})\in \mathcal{A}  $ and all $s \in Q$.  This implies that there exists a constant $0<\tilde{b}_1: =(g(|\xi_1|+|\xi_2|)-  F(\tilde{\varphi},\tilde{\theta},\tilde{\psi};b))/a   < \infty$, where we fix the functions $\tilde{\varphi}$, $\tilde{\theta}$ and $\tilde{\psi}$, so that, for any $s\in Q$,
\begin{equation}\label{el2}
 \sup_{s \in Q} D(\varphi_s,\theta_s,\psi_s)   \leqslant \tilde{b}_1.
\end{equation}

Let $s_i   \in Q$ for $i=1,2$.  Then we may bound
\begin{align}
\alpha(s_1)
=&   F(\varphi_{s_1},\theta_{s_1}, \psi_{s_1};s_2) + ( s_2-s_1) D(\varphi_{s_1}, \theta_{s_1}, \psi_{s_1}) \nonumber    \\
\leqslant  &\alpha(s_2) + |{s_1 - s_2}| D(\varphi_{s_1},\theta_{s_1},  \psi_{s_1}) ,
\end{align}
which, together with \eqref{el2}, yields
\begin{equation}
\alpha (s_1)
\leqslant  \alpha(s_2) + \tilde{b}_1|{s_1- s_2}|.
\end{equation}
Reversing the role of the indices $1$ and $2$ in the derivation of the above inequality gives the same bound with the indices switched.  Therefore we deduce that
\begin{equation}
 |{\alpha(s_1) - \alpha(s_2)}| \leqslant \tilde{b}_1 |{s_1-s_2}|,
\end{equation}
which proves the first assertion.

(2) To prove \eqref{el00} we note that the equality in \eqref{Elowerbound} and the non-negativity of $D$ imply that
\begin{equation}
\alpha(s) \leqslant g (|\xi_1|+|\xi_2|)- s \inf_{(\varphi,\theta,\psi)\in\mathcal{A}}    D(\varphi,\theta,\psi),
\end{equation}
where $(\varphi,\theta,\psi)\in \mathcal{A}$.
It is a simple matter to see that this infimum, which we call the constant $b_1$, is positive.

(3) To prove the third assertion, note that if  $0 < s_1 < s_2< \infty$,  then
\begin{equation}
\alpha(s_1) = F(\varphi_{s_1},\theta_{s_1}, \psi_{s_1};s_1) \geqslant F(\varphi_{s_2}, \theta_{s_2}, \psi_{s_2};s_1) \geqslant F(\varphi_{s_2}, \theta_{s_2}, \psi_{s_2}; s_2) = \alpha(s_2).
\label{202023020712}
\end{equation}
This shows that $\alpha$ is non-increasing in $s$.  Now suppose  by way of contradiction that $\alpha(s_1) = \alpha(s_2)$, then the previous relation \eqref{202023020712} implies that
\begin{equation}
 s_1 D(\varphi_{s_2}, \theta_s,\psi_{s_2}) = s_2 D(\varphi_{s_2}, \theta_{s_2},\psi_{s_2}),
\end{equation}
which means that $D(\varphi_{s_2},\theta_{s_2}, \psi_{s_2})=0$.  This in turn forces $\varphi_{s_2} = \theta_{s_2} = \psi_{s_2} = 0$, which contradicts the fact that $(\varphi_{s_2}, \theta_{s_2},\psi_{s_2}) \in  \mathcal{A}$.  Hence equality $\alpha(s_1) = \alpha(s_2)$ cannot be achieved, and $\alpha$ is strictly decreasing in $s$.

(4) Finally, we prove the fourth assertion. Obviously, it suffices to show that, under the condition \eqref{202307212216},
\begin{equation}
 \sup_{(\phi,\chi,\omega)\in (H_0^1(h_-,h_+))^3} \frac{F(\phi,\chi,\omega)}{\mathcal{J}(\phi,\chi,\omega)} > 0,
 \label{202306031216}
\end{equation}
but since $\mathcal{J}$ is positive definite, we may reduce to constructing a  triple $(\varphi,\theta, \psi)\in  H_0^1(h_-,h_+) $ such that $F(\varphi,\theta,\psi) >0$. Next we divide the construction into two cases.

 (a) If $\xi_1\neq0$, we let $\theta=0$ and $\varphi = -\psi'/\xi_1$. To obtain \eqref{202306031216},  we must then construct $\psi \in H_0^2(h_-,h_+)$ so that
\begin{align}
\tilde{F}(\psi) := & F(-\psi'/\xi_1,0,\psi) \nonumber \\
= &- \int_{h_-}^{h_+}\left(2  g \bar{\rho}\psi \psi' +s {\mu} \left( \left( \frac{\psi''}{\xi_1} + \xi_1 \psi  \right)^2 +( \xi_2\psi)^2 + (4 +(\xi_2/\xi_1 )^2)(\psi')^2 \right) \right)\mm{d}y_3\nonumber \\
&- {\vartheta |\xi|^2}{\psi}^2(0)>0 , \label{2022306031253x}
\end{align}
where, by the   integration by parts  and  \eqref{201611051547}$_1$,
\begin{align}\label{ni1}
-2\int_{h_-}^{h_+}g\bar{\rho}\psi \psi' \mm{d}y_3=& g \llbracket\bar{\rho}  \rrbracket \psi^2(0)+ \int_{h_-}^{h_+}g \bar{\rho}' \psi^2 \mm{d}y_3
 =   {g \llbracket\bar{\rho} \rrbracket \psi^2(0)}- {g^2} \int_{h_-}^{h_+} \frac{\bar{\rho}}{P'(\bar{\rho})} \psi^2 \mm{d}y_3 .
\end{align}

Now we define the test function $\psi_\beta \in H_0^2(h_-,h_+)$ for $\beta \geqslant 5$  according to
\begin{equation*}
 \psi_\beta(y_3) =
\begin{cases}
  \displaystyle \left(1-\frac{y_3^2}{h^2_+} \right)^{\beta/2} &\mbox{for } y_3 \in [0,h_+); \\
   \displaystyle\left(1-\frac{y_3^2}{h_-^2} \right)^{\beta/2} &\mbox{for }  y_3\in (h_-,0).
\end{cases}
\end{equation*}
By Lebesgue's dominated convergence theorem, it is easy to see that
\begin{align}
 \int_{h_-}^{h_+} (\psi_\beta)^2 \mm{d}y_3 =   o (\beta).
 \label{2022307302250}
\end{align}
where   $ o (\beta)$ represents a quantity that vanishes as $\beta \rightarrow \infty$.

In addition,
\begin{equation}
\int_{h_-}^{h_+}\mu  \left(  \left( \frac{\psi''_\beta}{\xi_1} + \xi_1 \psi_\beta  \right)^2  +( \xi_2\psi_\beta)^2 + (4+(\xi_2/\xi_1 )^2)(\psi'_\beta)^2  \right) \mm{d}y_3 \leqslant \tilde{b}_2\label{202dsa2307302250}
\end{equation}
for some constant $\tilde{b}_2$ depending on  $\beta$, $\xi$, $\mu$ and $h_\pm$.  Exploiting \eqref{2022306031253x}--\eqref{202dsa2307302250} and the fact $\psi_\beta(0)=1$,  we find that
\begin{equation}
\label{202307212241}
 \tilde{F}(\psi_\beta) \geqslant  {g \llbracket\rho\rrbracket-\vartheta \abs{\xi}^2}  + o(\beta ) - \tilde{b}_2 s,\end{equation}
  where $ o (\beta)$ also represents a quantity that vanishes as $\beta \rightarrow \infty$, and depends on $g$, $\bar{\rho}$, $h_\pm$ and $P_\pm'$.
  Since $\xi$ satisfies the condition \eqref{202307212216}, we may then fix $\beta$ sufficiently large so that the sum of the first three terms is equal the half of the sum of the first two terms.  Then there exists  $s_0>0$ depending on   $\xi$, $g$, $\vartheta$, $\mu$, $\bar{\rho}$, $h_\pm$ and  $P_\pm$,  so that for $s\leqslant s_0$ it holds that
\begin{align}
\label{202306031307}
\tilde{F}(\psi_\beta) > 0,
\end{align} thereby proving the desired result for $\xi_1 \neq0$.

(b) If
\begin{align}\xi_1=0\mbox{ and }\xi_2\neq0,
\label{2023061419851}
 \end{align}we let $\varphi=0$ and $\theta = -\psi'/\xi_2$, and thus
\begin{align}
 &F(0,-\psi'/\xi_2, \psi) \nonumber \\
&= - \int_{h_-}^{h_+}\left(2  g \bar{\rho}\psi \psi' + s{\mu} \left( \left( \frac{\psi''}{\xi_2} + \xi_2 \psi  \right)^2
+(\xi_1 \psi)^2+ (4 +(\xi_1/\xi_2)^2)(\psi')^2\right) \right)\mm{d}y_3\nonumber \\
& \quad   - {\vartheta |\xi|^2}{\psi}^2(0)>0 .
\label{20222306031253x}
\end{align}
Comparting both   the  structures of \eqref{2022306031253x} and \eqref{20222306031253x}, and recalling the derivation of \eqref{202306031307}, we easily see the fourth assertion also holds under the case \eqref{2023061419851}. This completest the proof. \hfill$\Box$
\end{pf}

With Lemmas \ref{minexist}--\ref{eigenlip} in hand, we are in the position to the proof of the existence of the solutions to \eqref{2022305251953}--\eqref{202305282005}.
\begin{pro}\label{coupledsolution}  Let  the assumptions in Theorem \ref{thm:0202} and the condition \eqref{202307212216} be satisfied.
\begin{enumerate}[(1)]
  \item There exist solutions $\varphi (\xi,y_3)$, $\theta (\xi,y_3)$ and $\psi (\xi,y_3)$ with $\lambda(\xi)>0$ to  the boundary value problem of \eqref{2022305251953}--\eqref{202305282005}. Moreover, $ \varphi (\xi,y_3)$, $\theta (\xi,y_3)$,  $\psi (\xi,y_3)\in H^4((h_-,0)\cup (0,h_+))$,
 \begin{equation}
\label{2023052safda91431x}
 \lambda^2=\sup_{(\phi,\chi,\omega)\in \mathcal{A}}F(\phi,\chi,\omega;\lambda)= F(\varphi ,\theta ,\psi;\lambda )>0
\end{equation}
and
 \begin{equation}
\label{20saf23052safda91431x}
 \psi (\xi,0),\  |\varphi(\xi,y_3)|+|\theta(\xi,y_3)| ,\  \psi (\xi,y_3) \neq0 .\end{equation}
  \item Let  $\varphi (\xi,y_3)$, $\theta (\xi,y_3)$ and $\psi (\xi,y_3)$ with $\lambda(\xi)>0$ be  constructed above. Then \begin{align}\lambda(\xi)=\lambda(-{\xi}),
\label{202306151601}
\end{align} and
$\varphi (\xi,y_3)$, $\theta (\xi,y_3)$, $-\psi (\xi,y_3)$ with $\lambda(-\xi)>0$  also are the solutions of the  boundary value problem of \eqref{2022305251953}--\eqref{202305282005}.\end{enumerate}
\end{pro}
\begin{pf}(1) By  Lemma \ref{eigenlip}, there exists a positive constant  $\gamma$ such that
\begin{align}
\label{20223061421563}
0<\alpha(s)\in C^0(0,\gamma]\mbox{ and }\alpha(\gamma)=0;
\end{align}
moreover, $\alpha(s)$ is strictly decreasing.
Now we define
\begin{align}
\beta(s):= \sqrt{\alpha(s)}, \label{2022306151244}
 \end{align}then $\beta\in C^{0}( 0,\gamma]$ is a strictly decreasing and  $\beta(s)>0$ for any $s\in (0,\gamma)$. Obviously  there exists a unique $\lambda$ such that $\beta(\lambda):= \lambda$; in particular, $\lambda^2=\alpha(\lambda)$. Moreover, by Lemma \ref{minexist}, there exists $(\varphi (\xi,\cdot),\theta (\xi,\cdot)$, $\psi(\xi,\cdot))$, which belongs to
 $\mathcal{A}\cap(H^4((h_-,0)\cup (0,h_+)))^3$ and is the solution  to \eqref{202305291213} with $s=\lambda$ and $\alpha(s)=\lambda^2$.  In particular, $(\varphi(\xi,\cdot),\theta(\xi,\cdot)$, $\psi({\xi},\cdot))$ with $\lambda(\xi)$ is also the desired solution to the problem \eqref{2022305251953}--\eqref{202305282005} due to $\lambda>0$, and satisfies
$ \lambda^2= F(\varphi ,\theta ,\psi;\lambda )$. Hecne \eqref{2023052safda91431x} holds.
In addition, by the    definitions of $E$ and $F$, the identity \eqref{2022306031524}    and the fact $F(\varphi ,\theta ,\psi;\lambda )>0$, we easily get  \eqref{20saf23052safda91431x}.

(2)   Recalling the definition of $\alpha$, we can see that $\alpha(\xi,s)=\alpha(-\xi,s)$. Therefore we immediately get \eqref{202306151601} by recalling the construction of $\lambda(\xi)$. Finally, since  $(\varphi (\xi,\cdot), \theta (\xi,\cdot), \psi (\xi,\cdot))$ with $\lambda(\xi)$ is the solution of the  boundary value problem of \eqref{2022305251953}--\eqref{202305282005}, thus we easily observe that  $(\varphi (\xi,\cdot), \theta (\xi,\cdot), -\psi (\xi,\cdot))$ with $\lambda(-\xi)$ is also the solution of  \eqref{2022305251953}--\eqref{202305282005}. The proof of Proposition \ref{coupledsolution} is complete. \hfill $\Box$\end{pf}

\subsection{ Behavior of $\lambda$ with respect to $\xi$}

In this section we will study the behavior of $\lambda$ from Proposition \ref{coupledsolution} in terms of $\xi$.
\begin{lem}\label{lambdacont}
Let  the assumptions in Theorem \ref{thm:0202} and the condition \eqref{202307212216} be satisfied.
   The function $\lambda:(0,|{\xi}|_{\mm{c}}) \rightarrow \mathbb{R}^+$ is continuous
 and satisfies the upper bound
 \begin{align}   {\lambda(\xi)}\leqslant {\frac{ 3h_+ g\llbracket \bar{\rho}\rrbracket}{2\mu_+} } .
 \label{202230731012028} \end{align}
Moreover,  \begin{align}
\label{202308030959} \lim_{|{\xi}|\rightarrow 0^+} \lambda({\xi})=0,
  \end{align}
  and, if $\vartheta>0$, then also
  \begin{equation}
 \lim_{|{\xi}|\rightarrow |{\xi}|_{\mm{c}}^-} \lambda(\xi) = 0.
\end{equation}
\end{lem}
\begin{pf}
(1)   Recalling the definition of $\alpha$, we easily see that $\alpha$ depends on $\xi$ for any given positive constant $s$. Therefore we denote $\alpha$ by $\alpha(\xi)$ or $\alpha(\xi,s)$. To obtain the desired continuity claim of $\lambda(\xi)$, we will first prove that, for any given $s>0$,  $\alpha(\xi)$ is continuous with respect to $\xi$.
To this purpose, we choose a point $\xi^0\in \mathbb{R}^2$ such that   $|\xi^0|\in (0, \varpi)$, where we have defined that
$$ \varpi:= \begin{cases}
 |{\xi}|_{\mm{c}} &\mbox{if }\vartheta>0;\\
2| \xi^0| &\mbox{if }\vartheta=0.
 \end{cases} $$
 Next we prove that  $\alpha(\xi)$ is continuous at $\xi^0$ for given $s>0$.

Without loss of a generality, it suffices to consider the case $\xi^0_1\neq 0$.  We assume that $\xi$ satisfies
\begin{align}
|\xi|\in (0, \varpi)\mbox{ and } \sigma:=|\xi-\xi^0  |\leqslant \min\{1, |\xi^0_1|/2\},
\label{2022307222151}
\end{align}
By virtue of Lemma \ref{minexist}, for any given $s>0$ and any given $\xi$,
there exists a triple
\begin{align}(\psi^\xi,\theta^\xi,\psi^\xi):= (\psi(\xi,y_3),\theta(\xi,y_3),\psi(\xi,y_3) )\in \mathcal{A},
\label{2022308080450}
\end{align} such that
\begin{align}\label{0243}
\alpha(\xi)={E(\psi^\xi,\theta^\xi,\psi^\xi)- s D(\psi^\xi,\theta^\xi,\psi^\xi)} .
\end{align}
In addition, we easily  see from the derivation of \eqref{202307212241} and the definition of $\alpha(\xi)$  that, for some positive constant $\tilde{b}_3:=\tilde{b}_3(\xi^0_1,\varpi ,g,\vartheta,\mu,\bar{\rho},h_\pm,
P_\pm) $,
 $$ \alpha(\xi)\geqslant ({g \llbracket\bar{\rho}\rrbracket-\vartheta \abs{\xi}^2} )/2  - \tilde{b}_3  s.$$
Exploiting \eqref{2022308080450} and Young's inequality, we can derive from  \eqref{0243}  and the above inequality  that, for some positive constant $\tilde{b}_4:=\tilde{b}_4(\xi^0_1,\varpi,s ,g,\vartheta,\mu,\bar{\rho},h_\pm$, $
P_\pm)$,
\begin{equation}\label{0244}\|(\psi^\xi,\theta^\xi,\psi^\xi)\|_{H^1_0(h_-,h_+)}\leqslant  \tilde{b}_4 .\end{equation}

Let $\sigma_i =\xi_i-\xi^0_i$. Then plugging $\xi_i=\xi^0_i+\delta_i $ into the expression of $\alpha(\xi)$ in \eqref{0243}, and then exploiting  \eqref{0244} and the condition $|\sigma_i|\leqslant \sigma\leqslant 1$, we can estimate that, for some positive constant $\tilde{b}_5:=\tilde{b}_5( \xi^0_1,\varpi ,s,g,\vartheta,\mu,\varsigma,\bar{\rho},h_\pm,P_\pm)$,
\begin{equation}\nonumber
 \alpha(\xi )- \alpha(\xi^0 ) \leqslant \tilde{b}_5\sigma \label{20223081111216} \end{equation}
 Similarly, we also have \begin{equation}\nonumber
 \alpha(\xi^0 )- \alpha(\xi ) \leqslant \tilde{b}_6\sigma, \nonumber \end{equation}
 which, together with \eqref{20223081111216},  implies, for any given $s>0$,
\begin{equation}\label{0249nn}
\alpha(\xi,s)\to  \alpha(\xi^0,s)\mbox{ as }
{\xi\rightarrow \xi^0} \mbox{ with }\xi \mbox{ satisfying \eqref{2022307222151}}.\end{equation}

Now we prove the continuity of $\lambda(\xi)$. Let $\beta(s)$ be  defined by \eqref{2022306151244}. Since $\beta(s)$ and its the definition domain $(0,\gamma)$  depend  on $\xi$, we denote them by $\beta(\xi,s)$ and $(0,\gamma_\xi)$, respectively.

Recalling \eqref{0249nn} and the proof of the first assertion in  Proposition \ref{coupledsolution}, we conclude  the following two facts.\begin{itemize}
       \item for any $\varepsilon>0$, there exists a $\delta>0$ such that, for any $\xi$ satisfying $|\xi|\in (0, \varpi)$ and   $|\xi-\xi^0|\leqslant \delta$,  it holds that
 $|\beta(\xi,s_{\xi^0})-
 \lambda(\xi^0,s_{\xi^0})|<\varepsilon$, where  $$(0,\gamma_\xi)\ni s_{\xi^0}=\lambda({\xi^0},s_{\xi^0})=\beta({\xi^0},s_{\xi^0})=\sqrt{\alpha({\xi^0,s_{\xi^0}})}>0.$$ \item for each fixed $|\xi|\in (0 ,|\xi|_{\mm{c}})$, $\beta(\xi, s)$ is continuous  and strictly decreasing with respect to
 $s\in (0,\gamma_\xi)$, and there exists a unique
 $s_{\xi}\in (0,\gamma_\xi)$ satisfying
$\lambda(\xi,s_{\xi})=\beta(\xi,s_{\xi})=s_{\xi}>0$.
     \end{itemize}
Consequently, we immediately infer that
$|\lambda(\xi,s_{\xi})-\lambda(\xi^0,s_{\xi^0})|<\varepsilon$
with  $s_{\xi}=\lambda(\xi,s^{\xi})$ from the above two facts. Therefore, $\lambda(\xi)$
is continuous on $(0,|{\xi}|_{\mm{c}}) $.

(2)  Now we turn to the  derivation of the upper bound \eqref{202230731012028}.
Thanks to \eqref{2022306031524} and \eqref{2023052safda91431x}, there exists  a triple $(\varphi ,\theta ,\psi )\in \mathcal{A}$ such that
\begin{align}\lambda^2  = (g\llbracket \bar{\rho}\rrbracket - {\vartheta |\xi|^2}){\psi}^2(0) -\int_{h_-}^{h_+}P'(\bar{\rho})\bar{\rho}\left(\frac{g }{ {P'(\bar{\rho}}) }
\psi -\xi_1\varphi-\xi_2\theta-\psi' \right)^2\mm{d}y_3 - \lambda D(\varphi,\theta,\psi) .  \label{20223073110014} \end{align}
By Newton--Leibniz formula and H\"older's inequality,  it is easy to see that
\begin{align}|\psi(0)|^2=&\frac{1}{4}\left|\int_0^{h_+}(\psi'-\xi_1\varphi-\xi_2\theta+
\psi'+\xi_1\varphi+\xi_2\theta)\mm{d}y_3\right|^2\nonumber \\
 \leqslant &{ \frac{h_+}{2}
\int_0^{h_+} ((\psi'-\xi_1\varphi-\xi_2\theta)^2+
(\psi'+\xi_1\varphi+\xi_2\theta)^2)\mm{d}y_3 }
 \leqslant  \frac{ 3h_+}{2\mu_+}D(\varphi,\theta,\psi).
 \end{align}
 Putting the above estimate into \eqref{20223073110014}  yields
\begin{align*}\lambda^2  \leqslant  \left(\frac{ 3h_+ g\llbracket \bar{\rho}\rrbracket}{2\mu_+}- \lambda \right)D(\varphi,\theta,\psi), \end{align*}
which, together the fact $D \geqslant 0 $, implies  \eqref{202230731012028}.

(3) Finally we  derive the limits as $|{\xi}| \rightarrow 0$ and $|{\xi}|_{\mm{c}}$.
   By the identity in \eqref{Elowerbound} and the first conclusion in Proposition \ref{coupledsolution},
for any $\xi$ satisfying $|\xi|\in (0,|\xi|_c)$, there exists $(\varphi,\theta,\psi)\in \mathcal{A}$ associated a $\lambda$ such that
$$0<\lambda^2=F(\varphi,\theta,\psi)\leqslant  g  \int_{h_-}^{h_+} \bar{\rho}\left( \xi_1(\varphi ^2 + \psi^2)+\xi_2 (\theta ^2 + \psi^2)  \right)\mm{d}y_3 - {\vartheta|{\xi}|^2} \psi^2(0),$$
which yields \eqref{202308030959} and
\begin{equation}
 \psi^2(0) < g(|\xi_1|+|\xi_2|) /{\vartheta|{\xi}| }^2,
 \nonumber \end{equation}
but by \eqref{2022306031524} we also know that
\begin{equation*}
 \lambda^2({\xi}) \leqslant  ({g  \llbracket\bar{\rho}\rrbracket - \vartheta |{\xi}|^2} )\psi^2 (0).
\end{equation*}
Chaining the above two inequalities together then shows that $\lim_{|{\xi} \rightarrow |{\xi}|_{\mm{c}}^-} \lambda( {\xi})=0$ for $\vartheta>0$. \hfill $\Box$
\end{pf}

\subsection{Solutions to the linearized periodic RT problem \eqref{linearized}}\label{growing_section}

In this  section we will construct a linear real-valued solution to the linearized RT problem
\eqref{linearized} which grows in-time in $H^4$.

\begin{pro} \label{growingmodesolneriodic}
\begin{enumerate}[(1)]  \item
Let
\begin{equation}
\label{2308031344}
\Lambda:= \sup_{|\xi|<|\xi|_{\mm{c}}}{\lambda(\xi)}<\infty.
\end{equation}there exists   a frequency
$\xi^1 $, such that $|\xi^1|\in (0,|\xi|_{\mm{c}})$ and
 \begin{align}&\lambda({\xi^1}) \in (2\Lambda /3,\Lambda]  \label{2202306151933} .
 \end{align}
 \item Let $\xi^2=-\xi^1$,
$c_7  =\lambda({\xi^1})$, $\varphi(\xi^1)$, $\theta(\xi^1)$ and $\psi(\xi^1)$ be the solutions provided by the first assertion in Proposition \ref{coupledsolution}, and $(\varphi(\xi^2),\theta(\xi^2), \psi(\xi^2))  =(\varphi(\xi^1),  \theta(\xi^1),-\psi(\xi^1))$.
  We define
\begin{align} \nonumber
& {w}(\xi^j,y_3) =  -\mm{i} \varphi(\xi^j,y_3) \mathbf{e}^1 - \mm{i} \theta(\xi^j,y_3) \mathbf{e}^2 + \psi(\xi^j,y_3) \mathbf{e}^3\mbox{ for }j=1, \ 2,\\
 &  \tilde{u}^0(y):=  \sum_{j=1}^2 (-1)^{j-1}{w}(\xi^j,y_3)  e^{\mm{i}y_{\mm{h}}\cdot \xi^j}\mbox{ and } \tilde{\eta}^0(y) = \tilde{u}^0(y)/c_7, \label{2022306101947}
\end{align}
Then \begin{align}
& \eta(y,t) =  e^{c_7 t} \tilde{\eta}^0(y)
\mbox{ and } u(y,t) =   e^{c_7 t}\tilde{u}^0(y)
\label{2022306041340}
\end{align}
 are real solutions to  the linearized periodic RT problem \eqref{linearized} with, for $k=1$ and $2$,
\begin{align}
L_k:= \begin{cases}
1/|\xi_k^1| &\mbox{if }\,\xi^1_k\neq 0;\\
1 &\mbox{if }\,\xi^1_k= 0.
\end{cases}
\label{20223072522255}
\end{align}
  \item
    There exist   constants  $b_i$ may depending on $\xi^1$, $h_\pm$, $\varphi(\xi^1,y_3)$, $\theta(\xi^1,y_3)$ and $\psi(\xi^1,y_3)$ such that
\begin{align}
& \sum_{  \beta_1+ \beta_2+ \beta _3\leqslant 4,\ 1\leqslant \beta _1+\beta _2}\|\partial_{1}^{\beta _1}\partial_{2}^{\beta _2} \chi_{n,n} \partial_{3}^{\beta _2}\tilde{u}^0\|_{0}^2/n  \leqslant  b_2 ,\label{2022306102034} \\
&  \| \chi_{n,n} \tilde{u}^0\|_4 /n \to b_3 ,\label{20223sfa06102034} \\
& \|\chi_{n,n} \tilde{u}^0_{\mm{h}}\|_{0}/n\to  b_4 , \label{20223060419451} \\
&\| \chi_{n,n}  \tilde{u}^0_3\|_{0}/n\to b_5 \mbox{ and }
  |\chi_{n,n} \tilde{u}^0_3|_{0}/n\to b_6 , \label{x2022306102034}
 \end{align}
 as $n\to \infty$, where $b_2\geqslant 0$ and $b_i>0$ for $3\leqslant i\leqslant 6$.
\end{enumerate}
\end{pro}
\begin{pf} (1) The first assertion in Proposition \ref{growingmodesolneriodic} is obvious, since $\lambda(\xi)$ is bounded and continuous by Lemma \ref{lambdacont}.

(2) Recalling Proposition \ref{coupledsolution} and the derivation of  \eqref{2022305251953}--\eqref{202305282005} from \eqref{linearized}, it is easy to observe that
$(\eta,u)$ defined by \eqref{2022306041340}  is  a real solution to  the linearized periodic RT problem \eqref{linearized} with $L_k$ be defined by \eqref{20223072522255}.

(3) Recalling the definitions of $\tilde{u}^0_1$ and $\chi_{n,n}$,  for any given multi-index  $(\beta_1,\beta_2,\beta_3)$ satisfying $\beta_1+\beta_2+\beta_3\leqslant 4$ and  $ \beta_1+\beta _2\geqslant 1$, there exists a non-negative constant $\tilde{b}_{7}:=\tilde{b}_{7}( \beta, h_\pm,\varphi(\xi^1,y_3) )$ such that
 \begin{align}
& \|\partial_{1}^{\beta _1}\partial_{2}^{\beta _2}\chi_{n,n}\partial_{3}^{\beta _3}\tilde{u}^0_1\|_{0}^2/n \nonumber \\
&= 4n^{-1} \int_{h_-}^{h_+}\partial_{3}^{\beta _3}\varphi^2(\xi^1,y_3)
\mm{d}y_3\int_{-n}^n\int_{-n}^n  (\partial_{1}^{\beta _1}\chi_{n }(y_1) \partial_{2}^{\beta _2}
\chi_{n }(y_2)\sin (\xi^1\cdot{y}_{\mm{h}}))^2\mathrm{d} {y}_ 1\mathrm{d} {y}_2
\nonumber \\
& \leqslant 16 \int_{h_-}^{h_+}\partial_{3}^{\beta _3}\varphi^2(\xi^1,y_3)
\mm{d}y_3\sup_{|y_1|<n}\{(\partial_{1}^{\beta _1}\chi_{n }(y_1))^2\}
\sup_{|y_2|<n}\{ (\partial_{2}^{\beta _2}
\chi_{n }(y_2))^2\}=:\tilde{b}_{7}
\end{align}
for any $n$.
Similarly we also have $\|\partial_{1}^{\beta _1}\partial_{2}^{\beta _2}\chi_{n,n}\partial_{3}^{\beta _3}\tilde{u}^0_i\|_{0}^2/n \leqslant \tilde{b}_{6+i}$, where $\tilde{b}_{6+i}\geqslant 0$ for $i=2$, $3$. Thus we arrive at \eqref{2022306102034}.

Now we derive \eqref{20223060419451}. From now on, we assume that $\beta _1+\beta _2\geqslant 0$ and make a convention that $0^0=1$. It is easy to see that
\begin{align}
& {n}^{-2}\int_{(-n,n)^2} \cos(2\xi^1\cdot y_{\mm{h}}) \mathrm{d} {y}_{\mm{h}}\nonumber \\
&=\frac{1}{n^2 }\int_{-n}^n \int_{-n}^{n} \left(\cos (2\xi_1^1y_1  ) \cos (2\xi_2^1y_2 )  - \sin (2\xi_1^1 y_1) \sin (2\xi_2^1 y_2)\right)  \mathrm{d} {y}_1 \mathrm{d} {y}_2\to 0, \nonumber
\end{align}
as $ n \to \infty$. We further derive from the above limit that
\begin{align}
& {n^{-2}}\int_{(-n,n)^2}
(\partial_{1}^{\beta _1}\partial_{2}^{\beta _2} \sin(\xi^1\cdot{y}_{\mm{h}}))^2\mathrm{d} {y}_{\mm{h}}\nonumber \\
& = \frac{1}{2n^2}\int_{(-n,n)^2}
|\xi_{1}|^{2\beta _1}|\xi_{2}|^{2\beta _2} \left(1+(-1)^{\beta _1+\beta _2+1}\cos(2\xi^1\cdot y_{\mm{h}})\right)\mathrm{d} {y}_{\mm{h}}
\to  2|\xi_{1}|^{2\beta _1}|\xi_{2}|^{2\beta _2}. \label{2022sda306161501}
\end{align}
Similarly, we also have
\begin{align}
& {n^{-2}}\int_{(-n,n)^2}
(\partial_{1}^{\beta _1}\partial_{2}^{\beta _2} \cos(\xi^1\cdot{y}_{\mm{h}}))^2\mathrm{d} {y}_{\mm{h}}\to 2|\xi_{1}|^{2\beta _1}|\xi_{2}|^{2\beta _2} ,\nonumber  \\
& {n^{-2}}\int_{(-n+1,n-1)^2}
(\partial_{1}^{\beta _1}\partial_{2}^{\beta _2}\Phi)^2\mathrm{d} {y}_{\mm{h}} \to2|\xi_{1}|^{2\beta _1}|\xi_{2}|^{2\beta _2} \mbox{ as }n\to \infty, \label{2022306161501}
\end{align}
where $\Phi=  \sin(\xi^1\cdot {y}_{\mm{h}})$ or $ \cos(\xi^1\cdot {y}_{\mm{h}})$.

Making use of the above limits in \eqref{2022sda306161501} and \eqref{2022306161501} and the fact
\begin{align*}
& \int_{(-n+1,n-1)^2}
(\partial_{1}^{\beta _1}\partial_{2}^{\beta _2} \sin(\xi^1\cdot {y}_{\mm{h}}))^2\mathrm{d} {y}_{\mm{h}} \\
&\leqslant   \int_{\mathbb{R}^2}(\chi_{n,n}\partial_{1}^{\beta _1}\partial_{2}^{\beta _2} \sin(\xi^1\cdot{y}_{\mm{h}}))^2\mathrm{d} {y}_{\mm{h}}\leqslant \int_{(-n,n)^2} (\partial_{1}^{\beta _1}\partial_{2}^{\beta _2} \sin(\xi^1\cdot{y}_{\mm{h}}))^2\mathrm{d} {y}_{\mm{h}},
\end{align*}
we obtain that
\begin{align*}
&   \frac{1}{n^2} \int(\chi_{n,n} \partial_{1}^{\beta _1}\partial_{2}^{\beta _2} \sin(\xi^1\cdot{y}_{\mm{h}}))^2\mathrm{d} {y}_{\mm{h}} \to 2|\xi_{1}|^{2\beta _1}|\xi_{2}|^{2\beta _2} \mbox{ for }n\to \infty ,
\end{align*}
which implies
\begin{align}
 \|\chi_{n,n}\partial_{1}^{\beta _1}\partial_{2}^{\beta _2}\partial_{3}^{\beta _3}\tilde{u}^0_1\|_{0}^2/n^2 =& 4n^{-2} \int_{h_-}^{h_+}(\partial_{3}^{\beta _3}\varphi(\xi^1,y_3))^2
\mm{d}y_3\int_{\mathbb{R}^2}(
\chi_{n,n} \partial_{1}^{\beta _1}\partial_{2}^{\beta _2} \sin(\xi^1\cdot{y}_{\mm{h}}))^2\mathrm{d} {y}_{\mm{h}}\nonumber \\
\to &   8|\xi_{1}|^{2\beta _1}|\xi_{2}|^{2\beta _2} \int_{h_-}^{h_+}(\partial_{3}^{\beta _3}\varphi(\xi^1,y_3))^2
\mm{d}y_3 \mbox{ as }n\to \infty . \label{202306161507}
\end{align}

In view of  \eqref{2022306102034}  and \eqref{202306161507},  there exists a non-negative constant $\tilde{b}_{10}:=\tilde{b}_{10}( \beta,\xi^1,h_\pm,\varphi(\xi^1,y_3) )$ such that
  $$ \| \chi_{n,n} \tilde{u}^0_1 \|_{4}/n\to  \tilde{b}_{10}\geqslant
   \sum_{  \beta_1+ \beta_2+ \beta _3\leqslant 4 } 8|\xi_{1}|^{2\beta _1}|\xi_{2}|^{2\beta _2} \int_{h_-}^{h_+}|\partial_{3}^{\beta _3}\varphi(\xi^1,y_3)|^2
\mm{d}y_3 .
  $$
Similarly, we have
 $$ \| \chi_{n,n} \tilde{u}^0_2 \|_{4}/n\to \tilde{b}_{11}\geqslant
   \sum_{   \beta_1+ \beta_2+ \beta _3\leqslant 4 } 8|\xi_{1}|^{2\beta _1}|\xi_{2}|^{2\beta _2} \int_{h_-}^{h_+}|\partial_{3}^{\beta _3}\theta(\xi^1,y_3)|^2
\mm{d}y_3
  $$
  and $$ \| \chi_{n,n} \tilde{u}^0_3 \|_{4}/n\to \tilde{b}_{12}\geqslant
   \sum_{  \beta_1+ \beta_2+ \beta _3\leqslant 4 } 8|\xi_{1}|^{2\beta _1}|\xi_{2}|^{2\beta _2} \int_{h_-}^{h_+}|\partial_{3}^{\beta _3}\psi(\xi^1,y_3)|^2
\mm{d}y_3 .
  $$
 Moreover $ \tilde{b}_{12}(\tilde{b}_{10}+\tilde{b}_{11}) > 0$ due to  $|\varphi(\xi,y_3)|+|\theta(\xi,y_3)|\neq0$, $ \psi(\xi,y_3) \neq0$ in \eqref{20saf23052safda91431x} and $|\xi|\neq 0$. Thanks to the above two limits, we immediately get \eqref{20223sfa06102034}.

In view of the derivation of \eqref{20223060419451}, we easily observe that  \eqref{20223060419451} and \eqref{x2022306102034} also hold.
This completes the proof of Proposition \ref{growingmodesolneriodic}. \hfill $\Box$
\end{pf}

\section{Gronwall-type energy inequality}\label{2022306101653}

Now we \emph{a prior} derive the Gronwall-type energy inequality for the solution $(\eta,u)$ of the RT  problem. To this end, we assume that $(\eta,u)$ satisfies
\begin{equation}\label{aprpiosesnew}
 {\sup_{0\leqslant  t \leqslant  T}(\|\eta(t)\|_3 +\|u(t)\|_2 )}\leqslant  \delta \mbox{ for some  }T>0,
\end{equation}
where $\delta\in  (0,\iota] $  is sufficiently small constant, and $\iota$ is the constant in Lemma \ref{201809012320}. It should be noted that the smallness of $\delta$ depends on the domain $\Omega$ and other known physical functions/parameters, such as $g$, $\vartheta$,
$\mu$, $\varsigma$, $\bar{\rho}$, $h_\pm$ and $P_\pm (\tau)$, in the  RT problem,
and will be repeatedly used in what follows.
In addition, by virtue of Lemma \ref{201809012320} and \eqref{aprpiosesnew} with sufficiently small $\delta$,  the expressions of $J^{-1}$,  $\mathbf{n}$, $\mathcal{H} $   and $\mathbf{N}_g$ in Section \ref{2022307272010} make sense.

\subsection{Estimates involving  $J$ and $\mathcal{A}$}

In this section, we derive some estimates involving the determinant $J$ and matrix $\mathcal{A}$.

\begin{lem}\label{lem:201612041032}
Under the assumption
 \begin{align}
 \label{2022306272107}
  \sup_{0\leqslant  t \leqslant  T} \|\eta(t)\|_3 \leqslant \delta\mbox{ with sufficiently small }\delta,
  \end{align}the determinant $J$  enjoys the following estimates: for $0\leqslant i\leqslant  2$,
\begin{align}
\label{Jdetemrinat}&1\lesssim \inf_{y\in {\Omega}}J\leqslant  \sup_{y\in  {\Omega}} J \lesssim1,
\\& \label{Jdetemrinatneswn}    \|(J-1, J^{-1}-1)\|_i\lesssim  \| \eta\|_{i+1},\\& \label{06011711}\|J^{-1}-1+\mm{div}\eta\|_i \lesssim \|\eta\|_{3}\|\eta\|_{i+1},\\
& \label{Jtestismtsnn}\|\partial_t(J,J^{-1}) \|_i\lesssim  \| u\|_{i+1} ,\\
 & \label{201612082025}\|J^{-1}_t+\mm{div}u\|_i\lesssim
 \| \eta\|_3\|  u\|_{i+1} .
\end{align}
\end{lem}
\begin{pf} By the definition of
$J$, we find that
\begin{align}
& J-1=\mm{div}\eta+P_2(\nabla \eta)+P_3(\nabla \eta)\mbox{ and }J_t=\mm{div}u+\partial_t(P_2(\nabla \eta)+  P_3(\nabla \eta)),  \label{201810221507}
\end{align}
where $P_i(\nabla \eta)$ denotes the homogeneous polynomial of degree $i$  with  respect to $\partial_j \eta_k$ for $1\leqslant i$, $ j$, $k\leqslant  3$. In addition,
 \begin{equation} J^{-1}-1=(1-J^{-1} )\left(\mm{div}\eta+P_2(\nabla \eta)+P_3(\nabla \eta)\right)-\left(\mm{div}\eta+P_2(\nabla \eta)+P_3(\nabla \eta)\right),\label{20181022150711}\end{equation}
which yields \begin{align}
J^{-1}_t =& -\mm{div}u-J^{-1}_t\left(\mm{div}\eta+P_2(\nabla \eta)+P_3(\nabla \eta)\right)-\partial_t\left(P_2(\nabla \eta)+P_3(\nabla \eta) \right)\nonumber \\
& -(J^{-1}-1)\partial_t\left(\mm{div}\eta + P_2(\nabla \eta)+P_3(\nabla \eta) \right) , \label{201810dfsa22150711}\end{align}
Thus, using  the smallness condition \eqref{2022306272107}, the embedding inequality $H^2\hookrightarrow L^\infty$ in \eqref{esmmdforinfty},   the product estimates (of $H^i$) in \eqref{fgestims} and the fact $ \eta_t=u$, we easily derive the desired estimates \eqref{Jdetemrinat}--\eqref{201612082025} from \eqref{201810221507}--\eqref{201810dfsa22150711} for sufficiently small $\delta$. \hfill$\Box$
\end{pf}

\begin{lem}\label{AestimstfroA}
Under the assumptions of  \eqref{2022306272107}, the  matrix $\mathcal{A}$ enjoys the following estimates:
 for  $0\leqslant  i\leqslant 2$,
\begin{align} & \label{aimdse}
 \sup_{y\in  {\Omega}}  |\mathcal{A} |  \lesssim 1 ,  \\
&\label{prtislsafdsfsfds}\|\tilde{\mathcal{A}}\|_{i} \lesssim \|  \eta\|_{i+1},\\
&\label{prtislsafdsfs}\|  \ml{A}_t\|_i \lesssim  \|  u\|_{i+1}.
\end{align}
\end{lem}
\begin{pf} We can compute out that, for sufficiently small $\delta$,
\begin{equation}
 \mathcal{A}^{\top}= (\nabla \eta+\mathbb{I})^{-1}= \mathbb{I}-\nabla \eta+(\nabla\eta)^2\mathcal{A}^{\top} \nonumber
 \end{equation}
and
\begin{equation}
 \mathcal{A}^{\top}_t=(\nabla\eta)^2\mathcal{A}^{\top}_t+( \nabla\eta\nabla u + \nabla u\nabla \eta )\mathcal{A}^{\top} -\nabla u. \nonumber
 \end{equation}
Recalling $\tilde{\mathcal{A}}=\mathcal{A}-\mathbb{I}$, then, similarly to Lemma \ref{lem:201612041032}, we derive the desired estimates \eqref{aimdse}--\eqref{prtislsafdsfs} from the above two identities for sufficiently small $\delta$. \hfill $\Box$\end{pf}

\begin{lem}
 \label{201702041500}
Under the assumption \eqref{2022306272107}, we  have,
for $0\leqslant  i\leqslant 2 $,
\begin{align}
 & 1\lesssim \inf_{y\in  {\Omega}} |\mathbf{n}| \leqslant \sup_{y\in  {\Omega}} |\mathbf{n} |  \lesssim 1, \label{202307291032} \\
 &  \|(J\mathcal{A}\mathbf{e}^3 - \mathbf{e}^3,\tilde{\mathbf{n} } )\|_i\lesssim \|\eta\|_{1,i}, \label{201702022056}  \\
 &\left \|\partial_t( J\mathcal{A}\mathbf{e}^3,\mathbf{n} )\right\|_i\lesssim \|u\|_{1,i},\label{201612071401xx}
\end{align}
where we have defined that $\tilde{\mathbf{n} }:={\mathbf{n} }-\mathbf{e}^3$.
\end{lem}
\begin{pf}
Recalling \eqref{20230626} and \eqref{05291021n}, we have
\begin{align}
&J\mathcal{A}\mathbf{e}^3  = {\mathbf{e}^3+ \mathbf{e}^1\times
\partial_2\eta+\partial_1\eta\times \mathbf{e}^2 +\partial_1\eta\times \partial_2\eta } , \label{202307231038}\\
&\partial_i| J\mathcal{A}\mathbf{e}^3|=(J\mathcal{A}\mathbf{e}^3)\cdot\partial_i(J\mathcal{A}\mathbf{e}^3) /|J\mathcal{A}\mathbf{e}^3|  , \label{20230sdf7231038}\\
&\tilde{\mathbf{n}}=(J {\mathcal{A}}\mathbf{e}^3-\mathbf{e}^3+(1- |J\mathcal{A}\mathbf{e}^3|) \mathbf{e}^3 )/|J\mathcal{A}\mathbf{e}^3|, \label{202032310581} \\
&\mathbf{n} _t=\partial_t( J\mathcal{A}\mathbf{e}^3) |J\mathcal{A}\mathbf{e}^3|^{-1}-
J\mathcal{A}\mathbf{e}^3|J\mathcal{A}\mathbf{e}^3|^{-3}(J\mathcal{A}\mathbf{e}^3)\cdot\partial_t (J\mathcal{A}\mathbf{e}^3),  \label{20203231058}
\end{align}
where $ \mathbf{e}^1=(1,0,0)^{\top} $ and $ \mathbf{e}^2=(0,1,0)^{\top} $.
Utilizing the product estimates and the embedding inequality $H^2\hookrightarrow L^\infty$, we deduce from \eqref{202307231038} and \eqref{20230sdf7231038}
that, for sufficiently small $\delta$,
\begin{align}
 &\| J\mathcal{A}\mathbf{e}^3 - \mathbf{e}^3\|_i \lesssim \|\eta\|_{1,i},\ \| \partial_t(J\mathcal{A}\mathbf{e}^3)\|_i \lesssim \|u\|_{1,i}, \label{060asd12130xx}\\
&\sup_{y\in  {\Omega}}|J\mathcal{A}\mathbf{e}^3|^{-1}  \lesssim 1,\ \left\|1-1/|J\mathcal{A}\mathbf{e}^3|\right\|_i\lesssim \|\eta\|_{1,i},\label{202307231250}
 \end{align}
and
\begin{align}\|1-|J\mathcal{A}\mathbf{e}^3|\|_i\lesssim \|\eta\|_{1,i} .  \label{20223072301039} \end{align}
In particular, by \eqref{202307231250}, we have
\begin{align}
\|f/|J\mathcal{A}\mathbf{e}^3|\|_i\lesssim \|f\|_i\mbox{ for any }f\in H^i .  \label{202saf307231250}
\end{align}
Thus, using \eqref{060asd12130xx}, \eqref{20223072301039}, \eqref{202saf307231250}, the embedding inequality of $H^2\hookrightarrow L^\infty$ and the product estimates, we can easily derive the  desired estimates in \eqref{202307291032}--\eqref{201612071401xx} from \eqref{202032310581} and \eqref{20203231058}. \hfill$\Box$
\end{pf}

\begin{lem}  Under the assumption \eqref{2022306272107}, we have
\begin{align}
& \|w\|_1^2\lesssim  \mathcal{U}_{\ml{A}}(w)\mbox{ for any }w\in H_0^1, \label{201701212009} \\
&\label{f11392016}
 \|w\|_{1+i}^2\lesssim\| \nabla_{\ml{A}}w\|_{i}^2  \lesssim\|\nabla w\|_i^2 \mbox{ for any }w\in H^{1+i},
\end{align}
where $0 \leqslant i\leqslant 2$, and we have defined that
\begin{align} \mathcal{U}_{\ml{A}}(w) :=\int \mathbb{S}_{\mathcal{A}}( w): \nabla_{\ml{A}}   w\mm{d}y.
\label{2022306192107}
\end{align}
\end{lem}
\begin{pf}
 Noting that
$$
 \mathcal{U}(w) = \mathcal{U}_{\ml{A}}(w) - \int \mathbb{S}( w):
 \nabla_{\tilde{\mathcal{A}}}    w\mm{d}y-\int \mathbb{S}_{\tilde{\mathcal{A}}}( w):
 \nabla_{ {\mathcal{A}}}   w\mm{d}y$$
 and
 $$   \mathcal{U}(w)  =  \frac{1}{2}\|\sqrt{\mu}(\mathbb{D}w-2\mm{div}w\mathbb{I}/3) \|_{0}^2+ \varsigma \|\mm{div}w\|_0^2  ,$$
thus, employing  Korn's inequality \eqref{2022306241238},  product estimate    and  \eqref{prtislsafdsfsfds}, we obtain \eqref{201701212009} for sufficiently small $\delta$.
Similarly to the derivation of \eqref{201701212009},  we easily see  that \eqref{f11392016} holds by further using Poinc\'are's inequality \eqref{2022306232041}.
\hfill $\Box$
\end{pf}

\subsection{Estimates of the nonlinear terms in the RT problem}
Now we proceed to derive some estimates on the nonlinear terms in the RT problem.
We first control the nonlinear terms $\mathbf{N}^3 $ and $(\mathbf{N}_1^4,\mathbf{N}_2^4,\mathcal{N})$ in the nonhomogeneous form \eqref{n0101nn1928M}.
\begin{lem}\label{201612041505}
Under the assumptions of \eqref{2022306272107}  with  $\delta\in (0,\iota]$, it holds that,   for   $   i=0$, $1$,
\begin{align}
 &  \|\mathbf{N}^3  \|_i \lesssim  \| \eta \|_{3} (\|(\eta,u)\|_{2+i}+\|u_t\|_{i}) \label{06011711jumpv},\\
&  |  \mathbf{N}^4 _{\mm{h}}    |_{i+1/2} +| \llbracket  R_P+\mathcal{N}^{u}  \rrbracket |_{i+1/2}
 \lesssim \|\eta\|_3 \|(\eta,u)\|_{2+i} ,\label{201807241645} \\
&|\mathcal{N}^{\eta}|_{y_3=0}   |_{1/2}\lesssim \|\mathcal{N}^{\eta}   \|_{1}\lesssim   \|\eta\|_3\| \eta \|_{2,1}.
    \label{060117sdfa34}
\end{align}
\end{lem}
\begin{pf}
(1) To being with, we bound $\mathbf{N}^3$. Using the product estimates, \eqref{Jdetemrinatneswn}, \eqref{prtislsafdsfsfds}  and  the regularity conditions of $\bar{\rho}$ and $P_\pm$ in Theorem \ref{thm:0202}, we infer that
\begin{align}
\|\mathbf{N}^3 \|_i\lesssim &\|(\mathbf{N}_g,\mathbf{N}_P) \|_i+(1+\|{\tilde{\ml{A}}}\|_2)\| {R}_{P} \|_{1+i}  +\|{\tilde{\ml{A}}}\|_2
 \|( \eta, u)\|_{2+i}  +\|J^{-1}-1\|_2\|u_t\|_i\nonumber\\
 \lesssim &\|\eta\|_3 (\|(\eta,u)\|_{2+i}+\|u_t\|_i ) +\|(\mathbf{N}_g ,\mathbf{N}_P)\|_i+\| {R}_{P} \|_{1+i}  \label{201612061920}.
\end{align}
In addition, if we use \eqref{Jdetemrinatneswn}, \eqref{06011711}, the regularity conditions of $\bar{\rho}$ and $P_\pm$, the homeomorphism property of $\eta+y$ (see \eqref{201803121601xx}),  we have, for sufficiently small $\delta$,
\begin{align}
&\|\mathbf{N}_g \|_i\lesssim \left\|   \int_{0}^{\eta_3}
(\eta_3-z) \frac{\mm{d}^2}{\mm{d}z^2}\bar{\rho}( y_3+z) \mm{d}z \right\|_i+ \|J^{-1}-1+\mm{div}\eta\|_i\lesssim \|\eta\|_3\|\eta\|_{1+i}, \label{2022307315} \\
&\|{\mathbf{N}}_P \|_i\lesssim \left\| \int_{0}^{\eta_3}
(\eta_3-z)\frac{\mm{d}^2}{\mm{d}z^2} \bar{P}(y_3+z)\mm{d}z\right\|_{1+i}\lesssim \|\eta\|_2\|\eta\|_{1+i}, \label{2022sfa307315}\\&
\| {R}_{P} \|_{i+1}\lesssim  \|(J^{-1}-1+\mm{div}\eta)\|_{1+i}
\nonumber \\
&\qquad \qquad\ + \left\|\int_{0}^{\bar{\rho}(J^{-1}-1)}(\bar{\rho}
(J^{-1}-1)-z)\frac{\mm{d}^2}{\mm{d}z^2} P (\bar{\rho}+z)\mm{d}z\right\|_{1+i}\lesssim
\|\eta\|_3   \|\eta\|_{2+i}.
\label{201612071026}
\end{align}
 Inserting the above three estimates into  \eqref{201612061920} immediately yields \eqref{06011711jumpv}.

(2)
Thanks to \eqref{201702022056}, the trace estimate \eqref{2022306232141}  and the product estimate, we easily estimates that
\begin{align}
\label{202307231255}
 |\Pi_{\mathbf{\mathbf{n}}}\mathbf{f}|_{i+1/2}\lesssim \|\mathbf{f}-(\mathbf{f}\cdot \mathbf{n} )\mathbf{n} \|_{1+i}\lesssim \|f\|_{1+i}\mbox{   for any }f\in H^{1+i} \mbox{ with  }  i=0, \ 1.\end{align}
 Making use of  \eqref{201702022056}, \eqref{202307231255}, trace estimate and the product estimate, we easily derive that
\begin{align}|\mathbf{N}^4_{\mm{h}}|_{i+1/2}
\lesssim  &
\|( \Upsilon  (\eta,u) \mathbf{e}^3)\cdot\tilde{\mathbf{n} }\mathbf{n}+ ( \Upsilon  (\eta,u)\mathbf{e}^3 ) \cdot \mathbf{e}^3\tilde{\mathbf{n} } \|_{1+i}\nonumber \\
 &+\|  \Upsilon  (\eta,u)   (J{\mathcal{A}}\mathbf{e}^3-\mathbf{e}^3)+
  \mathbb{S}_{\tilde{\mathcal{A}}}(u )J\mathcal{A}\mathbf{e}^3\|_{1+i}\lesssim \|\eta\|_3\|(\eta,u)\|_{2+i} \nonumber
\end{align}
Similarly, we have
\begin{align}| \llbracket   \mathcal{N}^{u}  \rrbracket |_{i+1/2}\lesssim \|\mathcal{N}^{u}\|_{1+i}\lesssim \|\eta\|_3\|u\|_{2+i}. \nonumber \end{align}
Thus we immediately derive \eqref{201807241645} from \eqref{201612071026}, the above two estimates  and trace estimate.

(3) Recalling the definitions of  $H^{\mm{n}}$ and $H^{\mm{d}}$ in \eqref{20sdfs1611041430M} and \eqref{20161asdfa1041430M}, and then using the product estimates, we
easily get, for sufficiently small $\delta$,
\begin{align}
 & \|H^{\mm{n}}\|_1\lesssim   \|\eta \|_{2,1}, \ \|H^{\mm{n}}\cdot \mathbf{e}^3-\Delta_{\mm{h}}\eta_3\|_1 \lesssim  \|\eta\|_3\|\eta\|_{2,1}, \ \|H^{\mm{d}}-1\|_2\lesssim \|\eta\|_{1,2}\label{0dsafdf6012130n}  ,\\
& \sup_{y\in \Omega}\{1/H^{\mm{d}}\}\lesssim 1\mbox{ and } \|1-1/H^{\mm{d}}\|_2\lesssim \|\eta\|_{1,2} \label{0dsafdf60sa12fsa130n}
\end{align}
Moreover, by \eqref{0dsafdf60sa12fsa130n}, we have
\begin{align} \label{2022307311845}
\|f/H^{\mm{d}}\|_2 \lesssim \|f\|_2\mbox{ for any }f\in H^2.
\end{align}

 Recalling the definition of $\mathcal{N}^\eta$ in \eqref{202005151045}, and then making use of \eqref{201702022056}, \eqref{0dsafdf6012130n},  \eqref{2022307311845}, the product estimate  and the trace estimate to get that
\begin{align}
|\mathcal{N}^{\eta}|_{y_3=0}   |_{1/2}\lesssim & \|\mathcal{N}^{\eta}  \|_{1}\lesssim \|(H^{\mm{n}} \cdot \mathbf{n}(H^{\mm{d}}-1)/H^{\mm{d}}\|_{1}+\|H^{\mm{n}}\cdot \tilde{ \mathbf{n}}\|_{1}+\|H^{\mm{n}}\cdot \mathbf{e}^3 -\Delta_{\mm{h}}\eta_3\|_{1}\nonumber \\
\lesssim &\|H^{\mm{n}}\|_{1}((1+\|\tilde{\mathbf{n}}\|_2) \|(H^{\mm{d}}-1)/H^{\mm{d}}\|_{2}+\|
\tilde{\mathbf{n} }\|_{2})+\|H^{\mm{n}}\cdot \mathbf{e}^3 -\Delta_{\mm{h}}\eta_3\|_{1}\nonumber \\
 \lesssim &  \|\eta\|_3\| \eta \|_{2,1}, \label{2022307241832}
\end{align}
which yields \eqref{060117sdfa34}.
 \hfill $\Box$
\end{pf}

To estimate the temporal derivatives of $u$ in Section  \ref{2022306171604}, we can
apply $\partial_t $ to \eqref{201611040926M} to derive that
\begin{equation}\label{n0101nnnn2026m}\begin{cases}
\bar{\rho}J^{-1}  u_{tt}-\mm{div}_{\ml{A}} ( {P}'(\bar{\rho})\bar{\rho}\mm{div}u) \mathbb{I}  +\partial_t\mathbb{S}_{\mathcal{A}}(u))=g\bar{\rho}(\mm{div}u \mathbf{e}^3-\nabla u_3 )+ \mathbf{N}^5&\mbox{in }  \Omega,\\[1mm]
 \llbracket    {P}'(\bar{\rho})\bar{\rho}\mm{div}u +  \partial_t\mathbb{S}_{\mathcal{A}}(u) \rrbracket  J\mathcal{A}\mathbf{e}^3 +\vartheta \Delta_{\mm{h}}u_3\mathbf{e}^3 =  \mathbf{N}^6,\ \llbracket   u_t  \rrbracket =0&\mbox{on }\Sigma,\\
  u_t=0 &\mbox{on }\partial\Omega\!\!\!\!-,
\end{cases} \end{equation}
where we have defined that
  \begin{align*}
\mathbf{N}^5:=& \mathbf{N}_t^1+ \mm{div}_{ {\ml{A}}_t} ( {P}'(\bar{\rho})\bar{\rho}\mm{div}\eta
\mathbb{I}+\mathbb{S}_{\mathcal{A}}(u) )-\bar{\rho} J^{-1}_t
u_t,\\
 \mathbf{N}^6 :=& \mathbf{N}^2_t - \llbracket{P}'(\bar{\rho})\bar{\rho}\mm{div}\eta \mathbb{I} +\mathbb{S}_{\mathcal{A}}(u)  \rrbracket  \partial_t (J\mathcal{A}\mathbf{e}^3)  .
\end{align*}Then we will establish the following estimates for both  the nonlinear terms $ \mathbf{N}^5$ and $ \mathbf{N}^6$.
\begin{lem}
\label{lem:0933}
Under the assumptions of \eqref{2022306272107} and $\delta\in (0,\iota]$,
\begin{align}
&  |  \mathcal{N}_t^\eta|_{y_3=0} |_{1/2} \lesssim \|  \mathcal{N}_t^\eta \|_{1} \lesssim \|\eta\|_3\| u \|_{3},\label{201806271120dsfs}\\
& \|\mathbf{N}^5\|_0 +|  \mathbf{N}^6  |_{1/2}
\lesssim \|\eta\|_3\|u\|_3+\|u\|_2 (\|u\|_3+\|u_t\|_1). \label{201806291}
\end{align}
\end{lem}
\begin{pf}
(1)
 Noting that
$$\ \|H^{\mm{n}}_t\|_1  \lesssim \|u\|_{1,2},\|H^{\mm{d}}_t\|_2  \lesssim \|u\|_{1,2} \mbox{ and } \|H^{\mm{n}}_t\cdot \mathbf{e}^3-\Delta_{\mm{h}}u_3\|_1 \lesssim  \|\eta\|_{1,2}\|u\|_{1,2},
$$
thus we can use \eqref{201702022056}, \eqref{201612071401xx}, \eqref{0dsafdf6012130n}, \eqref{2022307311845} and the above three estimates to derive that
\begin{align}
| \mathcal{N}^\eta_t|_{y_3=0}|_{1/2}\lesssim \|\mathcal{N}^\eta_t\|_1=  & \|\partial_t (H^{\mm{n}} \cdot \mathbf{n} (H^{\mm{d}}-1)/H^{\mm{d}}-H^{\mm{n}}\cdot \tilde{ \mathbf{n} }\nonumber \\
&-H^{\mm{n}}\cdot \mathbf{e}^3 +\Delta_{\mm{h}}\eta_3)\|_1
\lesssim\| \eta\|_3\|u\|_{1,2},   \nonumber
\end{align}
which yields \eqref{201806271120dsfs}.

(2)
Noting that
\begin{align}& \partial_t(\mathbf{N}_g-{\mathbf{N}}_P) \nonumber \\
&=g\left( u_3 \int_{0}^{\eta_3}
 \frac{\mm{d}^2}{\mm{d}z^2}\bar{\rho}( y_3+z) \mm{d}z - \bar{\rho} ( J^{-1}_t+\mm{div}u)\right)  \mathbf{e}^3\nonumber \\
&\quad +  \nabla_{\mathcal{A}_t}\left( \int_{0}^{\eta_3}
(z-\eta_3)\frac{\mm{d}^2}{\mm{d}z^2} \bar{P}(y_3+z)\mm{d}z\right)-\nabla_{\mathcal{A}}\left(
u_3 \int_{0}^{\eta_3}\frac{\mm{d}^2}{\mm{d}z^2} \bar{P}(y_3+z)\mm{d}z\right), \nonumber
  \end{align}
thus, following the arguments of \eqref{2022307315} and \eqref{2022sfa307315} by further exploiting \eqref{201612082025},  \eqref{aimdse} and  \eqref{prtislsafdsfs}, we  easily get  from the above identity that
\begin{align}
\label{202307261248}
\|  \partial_t(\mathbf{N}_g-{\mathbf{N}}_P) \|_0\lesssim \|\eta\|_3\|u\|_2.
\end{align}

Similarly, making use of \eqref{Jdetemrinatneswn} and \eqref{201612082025}, we have
 \begin{align}
 \|\partial_t R_P\|_1\lesssim \left\|   {P}'(\bar{\rho})\bar{\rho} (J^{-1}_t+\mm{div}u)
+\bar{\rho}J^{-1}_t\int_{0}^{\bar{\rho}(J^{-1}-1)}
\frac{\mm{d}^2}{\mm{d}z^2} P (\bar{\rho}+z)\mm{d}z \right\|_{1} \lesssim \|\eta\|_3\|u\|_2,
\label{20223070245204}
\end{align}
which, together with \eqref{aimdse}, \eqref{prtislsafdsfs} and \eqref{201612071026}, implies
\begin{align}
\label{20232307272053}
\|\partial_t\nabla_{\mathcal{A}}R_P\|_0\lesssim \|\eta\|_3\|u\|_2.
\end{align}
Thanks to  \eqref{prtislsafdsfsfds}, \eqref{prtislsafdsfs},  \eqref{202307261248} and \eqref{20232307272053}, we easily get
\begin{align}
\|\mathbf{N}_t^1\|_0= \|\partial_t( \mathbf{N}_g-{\mathbf{N}}_P-\nabla_{\mathcal{A}} {R}_P-g\nabla_{\tilde{\mathcal{A}}}(\bar{\rho}\eta_3) )\|_0\lesssim \|\eta\|_3\|u\|_2.
\label{20223080800630}
\end{align}
In addition, it is easy see that
\begin{align}
&   \|  \mm{div}_{ {\ml{A}}_t}  ({P}'(\bar{\rho})\bar{\rho}\mm{div}\eta \mathbb{I} +\mathbb{S}_{\mathcal{A}}(u) )-\bar{\rho} J^{-1}_t
u_t  \|_0\nonumber \\
 &\lesssim  \|\eta\|_3\|u\|_{ 2}+\|u\|_2 (\|u\|_3+\|u_t\|_1),\nonumber
 \end{align}which, together with \eqref{20223080800630}, yields
\begin{align}
\|\mathbf{N}^5\|_0
 \lesssim  &  \|\mathbf{N}_t^1\|_0 +\|\partial_t(  \mm{div}_{ {\ml{A}}_t} ( {P}'(\bar{\rho})\bar{\rho}\mm{div}\eta
\mathbb{I}+\mathbb{S}_{\mathcal{A}}(u) )-\bar{\rho} J^{-1}_tu_t)\|_0\nonumber \\
 \lesssim &\|\eta\|_3\|u\|_2+\|u\|_2 (\|u\|_3+\|u_t\|_1). \label{202307241824}
 \end{align}

Exploiting \eqref{202005151045}, \eqref{201702022056}, \eqref{201612071401xx}, \eqref{060117sdfa34},  \eqref{201612071026}, \eqref{201806271120dsfs}, \eqref{20223070245204} and trace estimate,  we  have
 \begin{align} |  \mathbf{N}^2_t |_{1/2}\lesssim &
     |\partial_t(\llbracket  R_P \rrbracket J  \mathcal{A}\mathbf{e}^3+  \vartheta (   \Delta_{\mm{h}}\eta_3 \mathbf{e}^3- \mathcal{H}J  \mathcal{A}\mathbf{e}^3 )) |_{1/2}   \nonumber \\
     \lesssim &  \|\partial_t ( R_P  J  \mathcal{A}\mathbf{e}^3)\|_1+\|  \partial_t(\mathcal{N}^{\eta}  J\mathcal{A}\mathbf{e}^3   -\Delta_{\mm{h}}\eta_3 (J\mathcal{A}\mathbf{e}^3-\mathbf{e}^3 )) \|_1
    \lesssim   \|\eta\|_3\|u\|_3.
     \label{20210880633}
 \end{align}
Similarly, we have
  \begin{align} | \llbracket{P}'(\bar{\rho})\bar{\rho}\mm{div}\eta \mathbb{I} +\mathbb{S}_{\mathcal{A}}(u)  \rrbracket  \partial_t (J\mathcal{A}\mathbf{e}^3)|_{1/2}   \lesssim (\|\eta\|_3+\|u\|_3)\|u\|_2,
 \end{align}
 which, together with \eqref{202307241824} and \eqref{20210880633}, yields  \eqref{201806291}. \hfill$\Box$
\end{pf}

\subsection{Basic estimates for $(\eta,u)$}\label{2022306171604}
In this subsection we  derive the $y_{\mm{h}}$-derivative estimates of $(\eta,u)$ in Lemmas \ref{201612132242nx} and \ref{201612132242nxsfs}, the temporal derivative estimates of $u$ in Lemma \ref{2016sadf12132242nx} and the
$y_{3}$-derivative estimates of $u$ in Lemma \ref{201asfd612132242nx}.

\begin{lem}\label{201612132242nx}Under the assumptions of \eqref{2022306272107}  and $\delta\in (0,\iota]$,
the following estimates hold:
\begin{align}
&\frac{\mm{d}}{\mm{d}t}\left( \int \bar{\rho}  \partial_\mm{h}^i \eta \cdot  \partial_\mm{h}^i u\mm{d}y+   \mathcal{U}(\partial_\mm{h}^i \eta)/2 \right) + \mathcal{I}(\partial_\mm{h}^i \eta) \nonumber \\
&\lesssim| \eta_3|_i^2+ \|   u\|^2_{i,0}
 + \sqrt{\mathcal{E}} \mathcal{D}\mbox{ for }
                                   \begin{cases}
                                                            0\leqslant i\leqslant 1 ;\\
                              ( i,\vartheta)=(2, 0),
                                   \end{cases}
\label{ssebdaiseqinM0846}   \\
& \frac{\mm{d}}{\mm{d}t}\int_{\mathbb{R}^2}
\left(  \int  \bar{\rho}  \mathfrak{D}_{\mf{h}}^{3/2} \partial_{\mm{h}}\eta \cdot   \mathfrak{D}_{\mf{h}}^{3/2} \partial_{\mm{h}} u \mm{d}y +  \mathcal{U}( \mathfrak{D}_{\mf{h}}^{3/2} \partial_{\mm{h}}\eta)/2 \right)\mm{d}\mf{h}   +
\vartheta|\nabla_{\mm{h}}\partial_{\mm{h}} \eta_3|_{1/2}^2\nonumber \\
 &\lesssim |\partial_{\mm{h}} \eta_3|_{1/2}^2+|\nabla_{\mm{h}}\partial_{\mm{h}} \eta_3|_0^2+\|u\|_{1,1}^2 + \sqrt{\mathcal{E}}\mathcal{D}\label{ssebdaisedsaqinM0dsfadsff846asdfadfad},\\
   &    \frac{\mm{d}}{\mm{d}t}  \left(\int\bar{\rho}  \partial_\mm{h}^2 \eta \cdot  \partial_\mm{h}^2 u \mm{d}y + \mathcal{U}(\partial_\mm{h}^2 \eta)/2 \right) + \| \partial^2_{\mm{h}} \mm{div}\eta\|_0^2  \nonumber \\
    & \lesssim    |   \partial^2_{\mm{h}} \eta_3 |^2_{0}+  \|  (  \eta_3 ,  u)\|^2_{2,0}+ |\partial_{\mm{h}}^2 \eta_3|_{1/2}   \left(\sqrt{\| \nabla_{\mm{h}} \mm{div}\eta\|_{\underline{1},0}\|\partial_3 \mm{div} \eta\|_{1,0}}
    \right.\nonumber \\
    &\quad \left.+ \|\nabla_{\mm{h}} \mm{div} \eta\|_{\underline{1},0} +\|\nabla_{\mm{h}} u\|_{\underline{1},1}+ \sqrt{\| \nabla_{\mm{h}} u\|_{\underline{1},1}\|u\|_{1,2}}  \right) + \sqrt{\mathcal{E}} \mathcal{D}
.
\label{ssebdaiseqinM0asdfa846}  \end{align}
\end{lem}
\begin{pf}   To begin with we derive the  estimate  \eqref{ssebdaiseqinM0846}. Applying $\partial_\mm{h}^i$ to \eqref{n0101nnnM}$_4$, \eqref{n0101nnnM}$_5$ and \eqref{n0101nn1928M}, we have
\begin{align}\label{n0101nnnnM}
&\bar{\rho} \partial_\mm{h}^iu_t+ \partial_\mm{h}^i(g\bar{\rho}(\nabla \eta_3- \mm{div}\eta \mathbf{e}^3)-\mm{div}\Upsilon(\eta,u))
 = \partial_\mm{h}^i\mathbf{N}^3  \mbox{ in }  \Omega,\\[1mm]
&\label{201612011052M}  \llbracket \partial_\mm{h}^i u  \rrbracket = \llbracket  \partial_\mm{h}^i\eta  \rrbracket =0,\ \partial_\mm{h}^i (\llbracket \Upsilon(\eta,u) \mathbf{e}^3
   \rrbracket
 +\vartheta\Delta_{\mm{h}}\eta_3 \mathbf{e}^3) = \partial_\mm{h}^{i} (\mathbf{N}_1^4,\mathbf{N}_2^4,\mathcal{N})^{\top }  \mbox{ on }\Sigma,\\
 &\label{201612011053M} \partial_\mm{h}^i\eta=\partial_\mm{h}^iu=0  \mbox{ on } \partial\Omega\!\!\!\!\!- ,
 \end{align}

Taking the inner product of \eqref{n0101nnnnM} and $\partial^i_\mm{h}\eta$ in $L^2$, we have
\begin{align}
 \frac{\mm{d}}{\mm{d}t} \int \bar{\rho} \partial_\mm{h}^i \eta \cdot \partial_\mm{h}^i u \mm{d}y =\int {\bar{\rho} }|\partial_\mm{h}^i u|^2\mm{d}y+\sum_{j=1}^3 I_{j}  ,\label{202306270924}
\end{align}
 where we have defined that
\begin{align*}
& I_1:= \int g\bar{\rho}\partial_\mm{h}^i( \mm{div} \eta \mathbf{e}^3-\nabla\eta_3 )\cdot  \partial_\mm{h}^i\eta\mm{d}y,\\
& I_2:=\int \mm{div}\partial_{\mm{h}}^{i}
\Upsilon(\eta,u) \cdot\partial_\mm{h}^i \eta\mm{d}y\mbox{ and }   I_3:=  \int \partial_\mm{h}^i{\mathbf{N}}^3\cdot  \partial_\mm{h}^i\eta\mm{d}y.
\end{align*}

Integrating by parts, using  \eqref{201612011052M}, the boundary condition of $\partial_\mm{h}^i\eta$ in \eqref{201612011053M}, and the symmetry of $\Upsilon $, we have
\begin{equation}  \label{201611222014}
 I_{1}= g \llbracket \bar{\rho}\rrbracket  |\partial_\mm{h}^i\eta_3|^2_0 +g\int ( \bar{\rho}'|\partial_\mm{h}^i\eta_3|^2 +2   \bar{\rho} \mm{div}\partial_\mm{h}^i\eta\partial_\mm{h}^i\eta_3)\mm{d}y
 \end{equation}
and
\begin{align}
I_{2} = &- \int \partial_\mm{h}^{i} \Upsilon(\eta,u) : \nabla  \partial_{\mm{h}}^i\eta\mm{d}y-\int_\Sigma
  \llbracket    \partial_{\mm{h}}^{i}
\Upsilon(\eta,u)\rrbracket  \mathbf{e}^3 \cdot\partial_{\mm{h}}^i\eta\mm{d}y_{\mm{h}} \nonumber \\
 =& -\int{P}'(\bar{\rho})\bar{\rho}|\mm{div} \partial_\mm{h}^i \eta|^2 \mm{d}y-\vartheta|\nabla_{\mm{h}}\partial_\mm{h}^i\eta_3|_0^2-\frac{1}{2}\frac{\mm{d}}{\mm{d}t}\mathcal{U}(\partial_\mm{h}^i\eta )+I_4, \label{201611222016}
  \end{align}
where we have defined that
\begin{equation*}
I_4:=-\int_\Sigma \partial_\mm{h}^{i} (\mathbf{N}_1^4,\mathbf{N}_2^4,\mathcal{N})^{\top } \cdot\partial_{\mm{h}}^i\eta\mm{d}y_{\mm{h}}.
\end{equation*}

Exploiting \eqref{201611051547}$_1$, we have the relation \begin{align}
 \int(   P'(\bar{\rho})\bar{\rho}|\mm{div}w|^2-g\bar{\rho}'w_3^2  -  2g
\bar{\rho}\mm{div}w w_3 )\mm{d}y  =\left\|\sqrt{P'(\bar{\rho})\bar{\rho}}\left(\frac{gw_3}{P'(\bar{\rho})} -
 \mm{div}w \right)\right\|_0^2 .
\label{2022307021259}
\end{align}
Making use of  \eqref{201611222014}, \eqref{201611222016}, the definition of $\mathcal{I}$ in \eqref{2022307201157} and the  above relation, we obtain
\begin{equation}  \label{estimforhoedsds1stm}
 \frac{\mm{d}}{\mm{d}t}\left(\int \bar{\rho}  \partial_\mm{h}^i \eta \cdot \partial_\mm{h}^i u \mm{d}y+ \mathcal{U}(\partial_\mm{h}^i \eta)/2\right)
+\mathcal{I}(\partial_\mm{h}^i \eta) \leqslant  c( |\eta_3|_{i}^2+  \|  u\|^2_{i,0})  +  I_3+I_4.  \end{equation}

Making use of  \eqref{06011711jumpv}--\eqref{060117sdfa34},  the integration by parts and H\"older's inequality, we get, for $0\leqslant  i\leqslant  1$ and $(\vartheta,i)=(0,2)$,
\begin{align}
I_3+I_4 \lesssim   \sqrt{\mathcal{E}} \mathcal{D}.
\label{201612011325M}
\end{align}
Putting the above estimates into \eqref{estimforhoedsds1stm} yields   \eqref{ssebdaiseqinM0846}.

Now we estimate for \eqref{ssebdaisedsaqinM0dsfadsff846asdfadfad}. Applying the fractional differential operator $\mathfrak{D}_{\mf{h}}^{3/2}$ (see the definition \eqref{2022306231307}) to \eqref{n0101nnnnM} with $i=1$, and then arguing in a way similar to that in the derivation of \eqref{estimforhoedsds1stm}, we get
\begin{align}
&
\frac{\mm{d}}{\mm{d}t}\left( \int \bar{\rho} \mathfrak{D}_{\mf{h}}^{3/2} \partial_{\mm{h}} \eta \cdot\mathfrak{D}_{\mf{h}}^{3/2}  \partial_{\mm{h}}  u\mm{d}y +  \mathcal{U} (\mathfrak{D}_{\mf{h}}^{3/2} \partial_{\mm{h}}\eta )/2 \right) + \mathcal{I}( \mathfrak{D}_{\mf{h}}^{3/2} \partial_{\mm{h}}\eta  )   +\vartheta  |\nabla_{\mm{h}} \partial_{\mm{h}} \eta_3|_{(y_{\mm{h}},y_3)=\mf{h}} |^2  \nonumber \\
    & \lesssim    |\nabla_{\mm{h}} \partial_{\mm{h}} \eta_3 |_{(y_{\mm{h}},y_3)= \mf{h} }|^2  + |\mathfrak{D}_{\mf{h}}^{3/2}  \partial_{\mm{h}} \eta_3|_{0 }^2+\|\mathfrak{D}_{\mf{h}}^{3/2} \partial_{\mm{h}}  u \|_0^2  \nonumber
\\
 &\quad+ \left|\int \mathfrak{D}_{\mf{h}}^{3/2} \partial_{\mm{h}}  \mathbf{N}^3  \cdot \mathfrak{D}_{\mf{h}}^{3/2} \partial_{\mm{h}}  \eta \mm{d}y -\int_{\Sigma}\mathfrak{D}_{\mf{h}}^{3/2}    \partial_{\mm{h}} (\mathbf{N}_1^4,\mathbf{N}_2^4,\mathcal{N})^{\top } \cdot \mathfrak{D}_{\mf{h}}^{3/2}  \partial_{\mm{h}} \eta \mm{d}y_{\mm{h}}\right|
=:I_5(\mf{h}). \label{ssebdaiseqinM0846asdfadfad}
 \end{align}
Making  use of   \eqref{201806291928},  \eqref{06011711jumpv}--\eqref{060117sdfa34}, the dual estimate \eqref{201808121247},  the definition of the norm $|\cdot|_{1/2}$  and a partial integration, we have
\begin{align}
\int_{\mathbb{R}^2}I_5(\mf{h})\mm{d}\mf{h}\lesssim& |\partial_{\mm{h}} \eta_3|_{1/2}^2+  |\nabla_{\mm{h}}\partial_{\mm{h}} \eta_3|_0^2+\|  u\|_{1,1}^2  +\|\eta\|_3 \|\mathbf{N}^3  \|_1 +|\nabla_{\mm{h}}\partial_{\mm{h}}\eta|_{1/2}| (\mathbf{N}_{\mm{h}}^4,\mathcal{N}) |_{1/2} \nonumber \\
\lesssim &|\partial_{\mm{h}} \eta_3|_{1/2}^2+  |\nabla_{\mm{h}}\partial_{\mm{h}} \eta_3|_0^2+\|  u\|_{1,1}^2+\sqrt{\mathcal{E}} \mathcal{D}. \nonumber
\end{align}
Thus, integrating \eqref{ssebdaiseqinM0846asdfadfad} with respect to $\mf{h}$ over $\mathbb{R}^2$, and then utilizing   the above estimate, we obtain  \eqref{ssebdaisedsaqinM0dsfadsff846asdfadfad}.

Finally, we turn to the derivation of \eqref{ssebdaiseqinM0asdfa846}. For $i=2$, we integrate by parts to find that
\begin{equation}
\label{201807022000}
I_2 =I_6-\int {{P}'(\bar{\rho})\bar{\rho}}  |\partial_{\mm{h}}^2  \mm{div}\eta|^2\mm{d}y  -\frac{1}{2}\frac{\mm{d}}{\mm{d}t}\mathcal{U}(\partial_{\mm{h}}^2\eta) ,
\end{equation}
where we have defined that
\begin{align*}
I_6:=&-\int_\Sigma    \llbracket \partial_{\mm{h}}^2({P}'(\bar{\rho})\bar{\rho}\mm{div}\eta+2\mu \partial_3u_3 +\left(\varsigma-{2\mu}/{3}\right)\mm{div}u )\rrbracket \partial_{\mm{h}}^2\eta_3\mm{d}y_{\mm{h}} -\int_\Sigma \partial_{\mm{h}}^2  \mathbf{N}^4_{\mm{h}}\cdot \partial_{\mm{h}}^2\eta_{\mm{h}}  \mm{d}y_{\mm{h}} .\nonumber
\end{align*}
Putting \eqref{201807022000} into \eqref{202306270924} with $i=2$, and then employing \eqref{201611222014} with $i=2$ and \eqref{2022307021259}, we conclude
\begin{align}
 & \frac{\mm{d}}{\mm{d}t}\int  \left( \bar{\rho}  \partial_\mm{h}^2 \eta \cdot  \partial_\mm{h}^2 u  +   \frac{1}{2} \mathcal{U}(\partial_\mm{h}^2 \eta )\right) \mm{d}y
 +  \left\|\sqrt{P'(\bar{\rho})\bar{\rho}}\partial_{\mm{h}}^2 \left(\frac{g\eta _3} {P'(\bar{\rho})} -
 \mm{div}  \eta \right)\right\|^2_0 \nonumber \\
 & \leqslant  c ( | \partial_{\mm{h}}^2\eta_3  |^2_{0}+ \|   u\|^2_{2,0} )+ I_6+I_7,\label{201807022042}
\end{align}
where we have defined that $I_7:=  \int  \partial_\mm{h}^2{\mathbf{N}}^3 \cdot  \partial_\mm{h}^2\eta\mm{d}y $.

Making use of  \eqref{06011711jumpv},  \eqref{201807241645}, the trace estimate,  dual estimate  and  integration by parts, we can  deduce that
\begin{align}
 I_6+I_7
\lesssim &   |\partial_{\mm{h}}^2 \eta_3|_{1/2}    \left(\sqrt{\| \nabla_{\mm{h}} \mm{div}\eta\|_{\underline{1},0}\|\partial_3 \mm{div} \eta\|_{1,0}}
    \right.\nonumber \\
    &\left.+ \|\nabla_{\mm{h}} \mm{div} \eta\|_{\underline{1},0} +\|\nabla_{\mm{h}} u\|_{\underline{1},1}+ \sqrt{\| \nabla_{\mm{h}} u\|_{\underline{1},1}\|u\|_{1,2}}  \right)+\sqrt{\mathcal{E}}\mathcal{D}.\label{2023071610911}
\end{align}
Consequently, inserting the above estimate  into \eqref{201807022042}, and then using Young's inequality, we  obtain \eqref{ssebdaiseqinM0asdfa846}.
\hfill$\Box$
\end{pf}
\begin{lem}\label{201612132242nxsfs} Under the assumptions of  \eqref{2022306272107}  and  $\delta\in (0,\iota]$,
the following estimates hold:
\begin{align}
&
 \frac{\mm{d}}{\mm{d}t}(\|\sqrt{\bar{\rho} } \partial_\mm{h}^i u\|^2_0+
  \mathcal{I}(\partial_{\mm{h}}^i \eta ))
+ c\|\partial_\mm{h}^i   u \|_{1}^2 \nonumber \\
&\lesssim  \sqrt{\mathcal{E} }\mathcal{D}+
   \begin{cases}\| \eta_3\|_{1}^2&\mbox{for }i=0; \\
  \| \eta_3\|_{i-1,1}^2&\mbox{for }i=1\mbox{ and }(\vartheta,i)=(0,2), \end{cases} \label{201702061418}  \\
&\frac{\mm{d}}{\mm{d}t}\int_{\mathbb{R}^2}\left( \|\sqrt{\bar{\rho} }   \mathfrak{D}_{\mf{h}}^{3/2} \partial_{\mm{h}}  u\|^2_0+
  \mathcal{I}( \mathfrak{D}_{\mf{h}}^{3/2} \partial_{\mm{h}}  \eta  )\right)\mm{d}\mathbf{h} \lesssim  \|\eta_3\|_{1,1}^2 + \sqrt{\mathcal{E}}\mathcal{D},
 \label{ssebdaisedsaqinM0dsfadsff846sdfaaasdfadfad}\\&\frac{\mm{d}}{\mm{d}t}  \left( \| \sqrt{\bar{\rho} }\partial_{\mm{h}}^2 u \|^2_0+\left\|\sqrt{P'(\bar{\rho})\bar{\rho}}\partial_{\mm{h}}^2 \left(\frac{g  \eta _3} {P'(\bar{\rho})} -
 \mm{div}  \eta \right) \right\|_0^2 \right)+ c\|\partial_\mm{h}^2 u\|_1^2  \nonumber \\
 &\lesssim  \|\eta_3\|_{1,1}^2+ |\partial_{\mm{h}}^2 u_3|_{1/2} \left(\sqrt{\| \nabla_{\mm{h}} \mm{div}\eta\|_{\underline{1},0}\|\partial_3 \mm{div} \eta\|_{1,0}}
    \right.\nonumber \\
    &\quad \left.+ \|\nabla_{\mm{h}} \mm{div} \eta\|_{\underline{1},0} +\|\nabla_{\mm{h}} u\|_{\underline{1},1}+ \sqrt{\| \nabla_{\mm{h}} u\|_{\underline{1},1}\|u\|_{1,2}}  \right)+ \sqrt{\mathcal{E}} \mathcal{D}.
 \label{2018072033}
 \end{align}
\end{lem}
\begin{pf}
Taking the inner product of \ \eqref{n0101nnnnM}$_1$  and $\partial^i_{\mm{h}}u$ in $L^2$, we obtain
\begin{align}
 & \frac{1}{2}\frac{\mm{d}}{\mm{d}t}\int \bar{\rho} |\partial_{\mm{h}}^i u|^2 \mm{d}y \nonumber\\
 & =  \int g\bar{\rho}\partial_\mm{h}^i( \mm{div} \eta \mathbf{e}^3-\nabla\eta_3 )\cdot  \partial_\mm{h}^i u\mm{d}y+\int \mm{div}\partial_{\mm{h}}^{i}
\Upsilon(\eta,u) \cdot\partial_\mm{h}^i u\mm{d}y   +\int \partial_\mm{h}^i{\mathbf{N}}^3\cdot \partial_\mm{h}^i u  \mm{d}y.
  \label{201807022044}
  \end{align}
   Thus, following the process used for \eqref{ssebdaiseqinM0846}, and applying Korn's inequality, we arrive at
\begin{align}
& \frac{\mm{d}}{\mm{d}t} \left(\|\sqrt{\bar{\rho}} \partial_{\mm{h}}^i u\|^2_0 +
  \mathcal{I}(\partial_{\mm{h}}^i \eta)\right)+c\| \partial_{\mm{h}}^i u\|^2_1\lesssim \left|
\int_\Sigma \partial_{\mm{h}}^i \eta_3 \partial_{\mm{h}}^iu_3\mm{d}y_{\mm{h}}  \right| + \sqrt{\mathcal{E}}\mathcal{D}.
 \label{estimforhoedsds1stnn1524}\end{align}

We can utilize the trace estimate and the dual estimate to get
\begin{align}
\left|
\int_\Sigma \partial_{\mm{h}}^i \eta_3 \partial_{\mm{h}}^iu_3\mm{d}y_{\mm{h}}
\right|\lesssim &                   \begin{cases}
  |\eta_3|_0|u_3|_0\lesssim \|\eta_3\|_1\|u_3\|_1 & \hbox{for }i=0; \\[1mm]
   | \partial_{\mm{h}}^{i-1} \eta_3|_{1/2} |\partial_{\mm{h}}^iu_3|_{1/2}\lesssim
\| \eta_3\|_{i-1,1} \|\partial_{\mm{h}}^i u_3\|_{1}   & \hbox{for }1\leqslant i\leqslant 2.
                   \end{cases} \label{202308031119}
\end{align}
Consequently, plugging the above estimate into \eqref{estimforhoedsds1stnn1524}, and then using Young's inequality,
we obtain  \eqref{201702061418} from \eqref{estimforhoedsds1stnn1524}.

 Similarly to \eqref{ssebdaiseqinM0846asdfadfad} and \eqref{estimforhoedsds1stnn1524}, we can derive from \eqref{n0101nnnnM} with $i=1$  that
\begin{align}
&\frac{1}{2}\frac{\mm{d}}{\mm{d}t}\left( \|\sqrt{ \bar{\rho}}\mathfrak{D}_{\mf{h}}^{3/2} \partial_{\mm{h}}  u\|^2_0+
  \mathcal{I}(\mathfrak{D}_{\mf{h}}^{3/2} \partial_{\mm{h}}  \eta )\right) +    \|\mathfrak{D}_{\mf{h}}^{3/2} \partial_{\mm{h}}u \|_1^2 \nonumber   \\
  & \lesssim |\mathfrak{D}_{\mf{h}}^{3/2} \partial_{\mm{h}}\eta_3|_0 |\mathfrak{D}_{\mf{h}}^{3/2} \partial_{\mm{h}}  u_3|_0
\nonumber
\\
 &\quad+ \left|\int \mathfrak{D}_{\mf{h}}^{3/2} \partial_{\mm{h}}  \mathbf{N}^3 \cdot\mathfrak{D}_{\mf{h}}^{3/2} \partial_{\mm{h}} u \mm{d}y -\int_{\Sigma} \mathfrak{D}_{\mf{h}}^{3/2}    \partial_{\mm{h}} (\mathbf{N}_1^4,\mathbf{N}_2^4,\mathcal{N})^{\top } \cdot\mathfrak{D}_{\mf{h}}^{3/2}  \partial_{\mm{h}} u \mm{d}y_{\mm{h}}\right| . \nonumber
 \end{align}
Thus, similarly to the argument of \eqref{ssebdaisedsaqinM0dsfadsff846asdfadfad},
 integrating the above inequality  with respect to $\mf{h}$ over $\mathbb{R}^2$,  we easily get \eqref{ssebdaisedsaqinM0dsfadsff846sdfaaasdfadfad} by further using trace estimate and  Young's inequality.

  Finally, we derive \eqref{2018072033}. Similarly to \eqref{201807022042}, we can derive from \eqref{n0101nnnnM} with $i=2$
that
 \begin{align}
 & \frac{1}{2}\frac{\mm{d}}{\mm{d}t} \left( \| \sqrt{\bar{\rho}   }\partial_{\mm{h}}^2 u \|^2_0+
  \left\|\sqrt{P'(\bar{\rho})\bar{\rho}}\partial_{\mm{h}}^2 \left(\frac{g   \eta_3} {P'(\bar{\rho})} -
 \mm{div}    \eta \right) \right\|_0^2   \right)+ \mathcal{U}(\partial_\mm{h}^2 u)\nonumber \\
 & =
\int_\Sigma \partial_{\mm{h}}^2 \eta_3 \partial_{\mm{h}}^2u_3\mm{d}y_{\mm{h}}+ I_8, \label{2018070220257}
 \end{align}
 where we have defined that
\begin{align*}
I_8:=&\int     \partial_\mm{h}^2{\mathbf{N}}^3\cdot  \partial_\mm{h}^2u \mm{d}y-\int_\Sigma \partial_{\mm{h}}^2  \mathbf{N}^4_{\mm{h}} \cdot \partial_{\mm{h}}^2u_{\mm{h}}\mm{d}y_{\mm{h}}\\
& -\int_\Sigma    \llbracket \partial_{\mm{h}}^2({P}'(\bar{\rho})\bar{\rho}\mm{div}\eta+2\mu \partial_3u_3 +\left(\varsigma-{2\mu}/{3}\right)\mm{div}u )\rrbracket \partial_{\mm{h}}^2u_3\mm{d}y_{\mm{h}}    .
   \end{align*}
Analogously to \eqref{2023071610911}, we can show
\begin{align*}
I_{8}\lesssim &   |\partial_{\mm{h}}^2 u_3|_{1/2} \left(\sqrt{\| \nabla_{\mm{h}} \mm{div}\eta\|_{\underline{1},0}\|\partial_3 \mm{div} \eta\|_{1,0}}
    \right.\nonumber \\
    &\quad \left.+ \|\nabla_{\mm{h}} \mm{div} \eta\|_{\underline{1},0} +\|\nabla_{\mm{h}} u\|_{\underline{1},1}+ \sqrt{\| \nabla_{\mm{h}} u\|_{\underline{1},1}\|u\|_{1,2}}  \right) +\sqrt{\mathcal{E}}\mathcal{D}.
\end{align*}
Putting the above estimate into \eqref{2018070220257} and then using \eqref{202308031119} with $i=2$ and Korn's inequality, we get \eqref{2018072033}.
\hfill$\Box$
 \end{pf}
\begin{lem}\label{2016sadf12132242nx}Under the assumptions of \eqref{aprpiosesnew} and $\delta\in (0,\iota]$,
the following estimates hold.
\begin{align}
  &\label{Lem:0301m0832}\frac{\mm{d}}{\mm{d}t}\left(\|\sqrt{\bar{\rho}   }  u_t\|_{0}^2
  + \mathcal{I}( u)\right)
 +c\|  u_t \|^2_{1}  \lesssim \| u_3\|_1^2+
 \sqrt{\mathcal{E} } \mathcal{D} ,\\
  &\label{201702071610nb}
\|u_t\|_{0}\lesssim  \|(\eta,u)\|_2.  \end{align}
\end{lem}
\begin{pf} Taking the inner product of  \eqref{n0101nnnn2026m}$_1$ and $J  u_t$ in $L^2$,  we have
\begin{align}
&\frac{1}{2}\frac{\mm{d}}{\mm{d}t}\int \bar{\rho}   | u_t|^2\mm{d}y=
g\int\bar{\rho} (\mm{div} u \mathbf{e}^3-
\nabla  u_3 )\cdot  u_t\mm{d}y\nonumber\\ &+\int J \mm{div}_{\ml{A}} ({P}'(\bar{\rho})\bar{\rho}\mm{div}u\mathbb{I}+\partial_t
\mathbb{S}_{\mathcal{A}}(u) ) \cdot u_t\mm{d}y+  \int J \mathbf{N}^5\cdot u_t\mm{d}y
\nonumber \\
& +
g\int\bar{\rho}(J-1)(\mm{div} u \mathbf{e}^3-
\nabla  u_3 )\cdot  u_t\mm{d}y= :\sum_{j=9}^{12}I_j .\label{0425m}
\end{align}

Similarly to \eqref{201611222014}, $I_9 $  can be written as
$$
I_9 = \frac{g}{2}\frac{\mm{d}}{\mm{d}t}\left(\int( \bar{\rho}'| u_3|^2+ 2 \bar{\rho} u_3 \mm{div}u)\mm{d}y\right)+g\llbracket  \bar{\rho}   \rrbracket  \int_\Sigma u_3 \partial_t u_3\mm{d}y_{\mm{h}}.
 $$Following the argument of deriving \eqref{201611222016}, we apply \eqref{AklJ=0}, \eqref{06051441}, \eqref{n0101nnnn2026m}$_2$ and \eqref{n0101nnnn2026m}$_3$ to get
\begin{align*}
I_{10} =& -\int J(   {P}'(\bar{\rho})\bar{\rho}\mm{div}u\mathbb{I} +\partial_t\mathbb{S}_{\mathcal{A}}(u) ): \nabla_{\mathcal{A}}  u_t  \mm{d}y\nonumber \\
&-\int_{\mathbb{R}^2}\llbracket (  {P}'(\bar{\rho})\bar{\rho}\mm{div}u \mathbb{I} +\partial_t\mathbb{S}_{\mathcal{A}}(u) ) J \mathcal{A} \mathbf{e}^3\cdot   u_t \rrbracket
\mm{d}y_{\mm{h}}\\
= &I_{13} - \frac{1}{2}\frac{\mm{d}}{\mm{d}t}(
  \|P'(\bar{\rho})\bar{\rho}\mm{div}u\|_0^2 + \vartheta |\nabla_{\mm{h}}u_3|_0^2 )- \mathcal{U}_{\sqrt{J}\mathcal{A}}(u_t)   ,
\end{align*}
where  $\mathcal{U}_{\sqrt{J}\mathcal{A}}(u_t)$ is defined by \eqref{2022306192107} with $\sqrt{J}\mathcal{A}$, resp. $u_t$ in place $w$ resp. $\mathcal{A}$, and
$$  \begin{aligned}
I_{13} :=&\int\left(
 (1-J)
 {P}'(\bar{\rho})\bar{\rho}\mm{div}u \mathbb{I} -J\mathbb{S}_{\ml{A}_t}
( u) \right): \nabla_{\mathcal{A}}u_t  \mm{d}y\\
&- \int  {P}'(\bar{\rho})\bar{\rho} \mm{div} u  \mm{div}_{\tilde{\ml{A}}}  u_t\mm{d}y
-\int_\Sigma   \mathbf{N} ^6 \cdot  u_t \mm{d}y_{\mm{h}}.
\end{aligned}$$
Inserting the above new expressions of $I_{9}$ and $I_{10}$ into \eqref{0425m},  we arrive at
\begin{equation} \label{060817561757m}
\frac{1}{2}\frac{\mm{d}}{\mm{d}t}\left(\|\sqrt{\bar{\rho} }  u_t\|^2_0 +
\mathcal{I}(u)  \right)+  \mathcal{U}_{\sqrt{J}\mathcal{A}}(u_t) =g\llbracket \bar{\rho}  \rrbracket \int_{\Sigma} u_3\partial_t u_3\mm{d}y_{\mm{h}}+ \sum_{j=11}^{13}I_{j}.
\end{equation}

Making use of  \eqref{Jdetemrinatneswn}, \eqref{201806291} and the product estimate, we have
\begin{equation}\nonumber
 I_{11} + I_{12}\lesssim(\|\mathbf{N}^5 \|_0 +\|J-1\|_2(\|  u\|_1 +\|\mathbf{N}^5 \|_0))\|  u_t\|_{0} \lesssim
\sqrt{\mathcal{E} } \mathcal{D} .
\end{equation}
In addition, making use of \eqref{Jdetemrinatneswn},  \eqref{prtislsafdsfsfds}, \eqref{prtislsafdsfs},  \eqref{201806291} and trace estimate, we have
\begin{align}
I_{ 13} \lesssim  &(\|\tilde{\ml{A}}\|_2+(1+\|\tilde{\mathcal{A}}\|_2 )
(\| J-1 \|_2(1+ \|{\mathcal{A}}_t\|_2)\nonumber\\
&+\|{\mathcal{A}}_t\|_2))\| u \|_1\| u_t\|_1
  +|\mathbf{N} ^6 |_0| u_t|_0
 \lesssim
\sqrt{\mathcal{E} } \mathcal{D} . \nonumber
\end{align}
 Plugging the  above two estimates into \eqref{060817561757m},  and then using \eqref{Jdetemrinat},  \eqref{201701212009},  Young's inequality and trace estimate, we arrive at \eqref{Lem:0301m0832}.

Finally, applying $\|\cdot\|_0$ to \eqref{n0101nn1928M}$_1$, and then exploiting \eqref{06011711jumpv} yields that
$$ \| \bar{\rho}  u_t\|_0 \lesssim  \|g\bar{\rho}(\mm{div}\eta \mathbf{e}^3 -\nabla \eta_3 )+\mm{div} \Upsilon(\eta,u)+{\mathbf{N}^3}\|_0 \lesssim \|(\eta,u)\|_2^2+\|\eta\|_3\|u_t\|_0^2,
 $$
which implies \eqref{201702071610nb} for sufficiently small $\delta$.
 \hfill$\Box$
\end{pf}
\begin{lem}\label{201asfd612132242nx}Under the assumptions of  \eqref{2022306272107}  and $\delta\in (0,\iota]$,
the following estimates hold:
\begin{align}
& \|u\|_{ 2}\lesssim
 \|\eta\|_{\underline{1},2} +\|u_t\|_0 ,
   \label{201702071610} \\
&  \|u \|_3  \lesssim \|\eta\|_{3}+ \|  u \|_{\underline{2},1}+\| u_{t}\|_1 ,\label{20161111195sdfafsd2}\\
& \|u \|_{1,2}   \lesssim \|(\mm{div}\eta,\nabla_{\mm{h}}\eta)\|_{ 1,1}+\|\nabla_{\mm{h}}u\|_{\underline{1},1}+\| u_{t}\|_{1}+\sqrt{\mathcal{E}\mathcal{D}}.\label{201808292010}
 \end{align}
\end{lem}
\begin{pf}
 We can rewrite \eqref{n0101nnnM}$_4$, \eqref{n0101nnnM}$_5$ and \eqref{n0101nn1928M} as a stratified Lam\'e problem  with jump conditions:
\begin{equation}\label{n0101nn928m} \begin{cases}
\mu\Delta   u +\tilde{\mu}\nabla \mm{div} u= \mathbf{F}^1 &\mbox{ in }  \Omega, \\[1mm]
  \llbracket  u  \rrbracket =0,\ \llbracket \mathbb{S}(u)\rrbracket  \mathbf{e}^3= \mathbf{F}^2&\mbox{ on }\Sigma, \\
 u=0 &\mbox{ on }\partial\Omega\!\!\!\!\!-,
\end{cases}\end{equation}
where we have defined that
\begin{align*}
  &\mathbf{F}^1 := g\bar{\rho}(\nabla \eta_3-\mm{div}\eta \mathbf{e}^3)-\nabla (P'(\bar{\rho}) \bar{\rho}\mm{div}\eta)+{\bar{\rho}} u_t-{\mathbf{N}^3} ,\\
&\mathbf{F}^2 :=(\mathbf{N}^4_1,\mathbf{N}^4_2,\mathcal{N})^{\top}- \llbracket  P'(\bar{\rho}) \bar{\rho}\mm{div}\eta \mathbf{e}^3 \rrbracket - \vartheta \Delta_{\mm{h}}\eta_3 \mathbf{e}^3\mbox{ and }\tilde{\mu}:=\varsigma+\mu/3.
\end{align*}
 Applying the stratified elliptic estimate \eqref{Ellipticestimate} to \eqref{n0101nn928m},
we find that
  \begin{align}
   \|   u\|_{2} \lesssim  \| \mathbf{F}^1 \|_0+  | \mathbf{F}^2 |_{1/2}\lesssim&
 \|\eta\|_{2}    + | ( \llbracket  P'(\bar{\rho}) \bar{\rho}\mm{div}\eta \mathbf{e}^3 \rrbracket,\Delta_{\mm{h}} \eta_3)|_{1/2}\nonumber \\
 &+\| (u_t,\mathbf{N}^3)\|_0+|(\mathbf{N}^4_1,\mathbf{N}^4_2,\mathcal{N})|_{1/2} ,
 \end{align}
 where $\mathcal{N}=\mathcal{N}^{\eta} +\llbracket  R_P+\mathcal{N}^{u}  \rrbracket   $.
Making use of \eqref{2022306272107},  \eqref{06011711jumpv}--\eqref{060117sdfa34} and trace estimate, we further get  \eqref{201702071610} from the above estimate.

 We rewrite \eqref{n0101nn928m}$_1$ as  a stratified Lam\'e problem with Dirichlet boundary conditions:
\begin{equation}\label{n0101nn928n} \begin{cases}
\mu\Delta   u  +\tilde{\mu} \nabla \mm{div} u =\mathbf{F}^1  &\mbox{in }  \Omega , \\[1mm]
u =u |_{\Sigma}& \mbox{on }  \Sigma, \\
 u =0 &\mbox{on }\partial\Omega\!\!\!\!\!- .
\end{cases}\end{equation}
Applying the elliptic estimate \eqref{xfsddfsf201705141252} to the above  Lam\'e  problems, we obtain
\begin{align}
 \| \partial_{\mm{h}}^i u\|_{3-i}  \lesssim \|   \mathbf{F}^1 \|_{i,1-i}+| \partial_{\mm{h}}^i u |_{5/2-i  } ,
 \label{20171111029}
 \end{align}
where $i=0$ and $1$.
In addition, we have
\begin{align*}
\|   \mathbf{F}^1 \|_{i,1-i } \lesssim &  \|(\mm{div}\eta,\nabla_{\mm{h}}\eta)\|_{i,2-i}    +\|( u_t,\mathbf{N}^3)\|_{i,1-i} \\ \lesssim &
\|(\mm{div}\eta,\nabla_{\mm{h}}\eta)\|_{i,2-i}+ \|u_t\|_{1} +\|\eta\|_3\|(\eta,u)\|_3.
\end{align*}
Putting the above estimate into \eqref{20171111029} and then using \eqref{201808051835} we further obtain  \eqref{20161111195sdfafsd2} and \eqref{201808292010}.
\hfill $\Box$
\end{pf}

\subsection{Highest-order boundary estimates of $u_3$ at interface for $\vartheta> 0$}\label{202306141110}

In this subsection we further establish a highest-order boundary estimate of $u_3$. We should remark that it is difficult to directly derive
the desired estimate based on the RT  problem \eqref{n0101nnnM}. Motivated by the well-posdeness result of  incompressible stratified viscoelastic fluids \cite{XLZPZZFGAR}, we will break up the nonhomogeneous form of the RT  problem into two subproblems.

To this purpose, we will first  use Lemma \ref{2022307241635}  to construct a function $\mathcal{M}\in C^0(\mathbb{R}^+,H^{1}(\mathbb{R}^2))$ such that
\begin{align}
&\mathcal{M}|_{t=0}=\mathcal{H}^0,\label{2020307272139} \\
&\|\mathcal{M}\|_{L^\infty(\mathbb{R}^+,H^{1}(\mathbb{R}^2))}+\|\mathcal{M}\|_{L^2(\mathbb{R}^+,
H^{3/2}
)}+\|\mathcal{M}_t
\|_{L^2(\mathbb{R}^+,H^{1/2})}\lesssim |\mathcal{H}^0|_{1}, \label{2022307252021}
 \end{align}
 where $\mathcal{H}^0=(\Delta_{\mm{h}}\eta_3-{\mathcal{N}}^\eta  /\vartheta)|_{t=0}=\mathcal{H}|_{t=0}$ (see \eqref{202005151045} for the definition of $\mathcal{H}$).
Let $({\eta}^1,{u}^1):=({\eta},{u})-({\eta}^2,{u}^2)$, where $({\eta}^2,{u}^2)$ is a solution to the linear problem
 \begin{equation}\label{s0106pndfsgnnnxx}
\begin{cases}
 \eta_t^2=  u^2 &\mbox{in }\Omega,  \\[1mm]
 \bar{\rho} u^2_{t}-\mm{div}\Upsilon(\eta^2,u^2)=\widetilde{\mathbf{N}}^3:= g\bar{\rho}(\mm{div}\eta  \mathbf{e}^3-\nabla \eta_3 ) + {\mathbf{N}}^3  &\mbox{in }\Omega, \\[1mm]
  \llbracket  \eta^2  \rrbracket =\llbracket  u^2  \rrbracket =0&\mbox{on }\Sigma,\\
 \llbracket    \Upsilon(\eta^2,u^2) \mathbf{e}^3  \rrbracket =
(\mathbf{N}_1^4,\mathbf{N}_2^4, \llbracket R_P  + \mathcal{N}^u \rrbracket-\vartheta\mathcal{M})^{\top } &\mbox{on }\Sigma, \\[1mm]
(\eta^2,u^2)=0&\mbox{on }\partial\Omega\!\!\!\!\!-,\\[1mm]
(\eta^2,u^2)|_{t=0}=(\eta^0,u^0)&\mbox{on }\Omega.
\end{cases}\end{equation}
It should be noted that, due to \eqref{2020307272139} and \eqref{s0106pndfsgnnnxx}$_6$,
 \begin{align}
 \llbracket    \Upsilon(\eta^2,u^2) \mathbf{e}^3  \rrbracket =
(\mathbf{N}_1^4,\mathbf{N}_2^4, \llbracket R_P  + \mathcal{N}^u \rrbracket-\vartheta\mathcal{M})^{\top } \mbox{on }\Sigma\mbox{ for }t=0,
\label{2022307312050}
\end{align}
 which makes sure  the existence of a unique solution of the above problem \eqref{s0106pndfsgnnnxx}, see the second conclusion in Proposition \ref{pro:0401nd}.
Then, $({\eta}^1,{u}^1)$ satisfies
\begin{equation}\label{s0106pndfsgnnn}
 \begin{cases}
 \eta_t^1=  u^1 &\mbox{in }\Omega,  \\[1mm]
 \bar{\rho} u_{t}^1=\mm{div}\Upsilon(\eta^1,u^1 ) &\mbox{in }\Omega, \\[1mm]
  \llbracket  \eta^1  \rrbracket =\llbracket  u^1 \rrbracket =0&\mbox{on }\Sigma,\\
 \llbracket   (\Upsilon(\eta^1,u^1))  \mathbf{e}^3   \rrbracket
=(\widetilde{\mathcal{N}}^\eta -\vartheta\Delta_{\mm{h}}\eta_3^1+\vartheta\mathcal{M} )\mathbf{e}^3&\mbox{on }\Sigma, \\[1mm]
(\eta^1,u^1)=0&\mbox{on }\partial\Omega\!\!\!\!\!-, \\[1mm]
(\eta^1,u^1)|_{t=0}=0&\mbox{on }\Omega,
\end{cases}
\end{equation}
where we have defined that $\widetilde{\mathcal{N}}^\eta :=\mathcal{N}^\eta-\vartheta\Delta_{\mm{h}}\eta_3^2 $.
Thus, we can use the above two auxiliary problems to derive the following highest-order boundary estimate of $u_3$.
\begin{lem}\label{201807251431}
Under the assumptions of $\vartheta> 0$, \eqref{2022306272107} and $\delta\in (0,\iota]$,
we have
\begin{align}
& \|  \eta^2  \|_{2,1}^2  + \int_0^t  | \nabla_{\mm{h}}^2 u_3|^2_{1/2} \mm{d}\tau \nonumber \\
& \leqslant  c\left(\|\eta^0\|_3^2+\|u^0\|_2^2 + |\mathcal{H}^0|_1^2 +\int_0^t\left(\|\eta\|_{2}^2+  \sqrt{\mathcal{E} }\mathcal{D}\right)\mm{d}\tau\right)
+ c_7  \int_0^t  \|  \eta^2 (t) \|_{2,1}^2 \mm{d}\tau,
\label{201807012235}
\end{align}
where $c_7$ is the constant after \eqref{2202306151933}  in Proposition \ref{growingmodesolneriodic} and  $|\nabla_{\mm{h}}^2 u_3|^2_{1/2}:= \sum_{|\alpha|=2}| \partial_{\mm{h}}^\alpha u_3|^2_{1/2}$.
\end{lem}
\begin{pf}
 (1)  Let $0\leqslant i\leqslant 1$.
Applying $\partial_{\mm{h}}^i\partial_t$ to \eqref{s0106pndfsgnnn} yields
\begin{equation}\label{sadfdxsafdsa}
\begin{cases}
\bar{\rho} \partial_{\mm{h}}^i u_{tt}^1 = \partial_{\mm{h}}^i\mm{div}\Upsilon(u^1,u^1_t)  & \mbox{in }\Omega, \\[1mm]
  \llbracket    \partial_{\mm{h}} u^1  \rrbracket =\llbracket    u^1_t \rrbracket =0&\mbox{on }\Sigma,\\
   \partial_{\mm{h}}^i\llbracket \Upsilon(u^1,u_t^1)   \mathbf{e}^3  \rrbracket
 =\partial_{\mm{h}}^i \partial_t (    \widetilde{\mathcal{N}}^\eta-\vartheta\Delta_{\mm{h}}\eta_3^1 +\vartheta\mathcal{M}  )\mathbf{e}^3 &\mbox{on }\Sigma,\\
  ( \partial_{\mm{h}}u^1, u^1_t)=0&\mbox{on }\partial\Omega\!\!\!\!\!-.
\end{cases}
\end{equation}

Multiplying \eqref{sadfdxsafdsa}$_1$ with $i=0$ by $  u^1_t$ in $L^2$ and then integrating by parts, we infer that
\begin{align}
 \frac{\mm{d}}{\mm{d}t}\left(\|\sqrt{\bar{\rho} } u^1_t\|^2_0+  \|  P'(\bar{\rho}) \bar{\rho}\mm{div}u^1 \|_0^2 + \vartheta|\nabla_{\mm{h}} u^1_3|_0^2 \right)+ \| \mathcal{U} (u^1_t)\|^2_0
  \lesssim |\partial_t   u^1_3|_{ 1/2}|  \partial_t (  \widetilde{\mathcal{N}}^\eta  ,\mathcal{M})|_{-1/2}  . \label{201811041954}
\end{align}
Obviously we  have from \eqref{s0106pndfsgnnn}$_2$ that
$$\| u_t^1\|_0^2\lesssim \|(\eta^1,u^1)\|_2^2  ,$$
which, together with the zero initial data \eqref{s0106pndfsgnnn}$_6$,  gives
\begin{equation}
\label{201811041957}
\| u_t^1|_{t=0}\|_0^2 \lesssim  0.
\end{equation}
Integrating \eqref{201811041954} with respect to $t$, and then using \eqref{201811041957}, Korn's, Young's inequalities and trace estimate, we further have
\begin{align}
  \|  u^1_t \|^2_0 +   |\nabla_{\mm{h}} u^1_3|_0^2 +  \int_0^t \| u^1_\tau\|^2_{1}\mm{d}\tau
  \label{201807011935}\lesssim   \int_0^t  |  \partial_{\tau} (  \widetilde{\mathcal{N}}^\eta ,\mathcal{M})|_{-1/2} ^2 \mm{d}\tau.
\end{align}

We multiply \eqref{sadfdxsafdsa}$_1$ with $i=1$ by $\partial_{\mm{h}}  u^1$ in $L^2$ and then integrate by parts to deduce
 \begin{align*}
 &\frac{\mm{d}}{\mm{d}t}\int\left(\frac{1}{2}\mathcal{U} ( \partial_{\mm{h}} u^1 )-
 \bar{\rho}  u_{t}^1\partial_{\mm{h}}^2 u^1 \right)\mm{d}y
+   \|  P'(\bar{\rho}) \bar{\rho}\partial_{\mm{h}}\mm{div}u^1 \|_0^2 +  \vartheta|\nabla_{\mm{h}}\partial_{\mm{h}} u^1_3 |_0^2\\
&\lesssim  \|   \partial_{\mm{h}}   u_{t}^1\|^2_0+|\nabla_{\mm{h}}\partial_{\mm{h}}u^1_3|_{0} |    \partial_t(\widetilde{\mathcal{N}}^\eta,\mathcal{M} )|_{0}  .
\end{align*}
Thus, integrating the above equality with respect to $t$, and using then \eqref{s0106pndfsgnnn}$_6$, and
Korn's, Young's inequalities, one infers that
\begin{align*}
   \| \partial_{\mm{h}} u^1\|^2_1
+ \int_0^t|\nabla_{\mm{h}}\partial_{\mm{h}} u^1_3 |_0^2 \mm{d}\tau \lesssim \|u_t\|_0^2+\int_0^t\left(\|   \partial_{\mm{h}}   u_{\tau}^1\|^2_0
+ |  \partial_{\tau}( \widetilde{\mathcal{N}}^\eta, \mathcal{M}_{\tau})|_{0}^2  \right)\mm{d}\tau ,
\end{align*}
which, together with \eqref{201807011935}, yields
\begin{align}
 \int_0^t |\nabla_{\mm{h}} \partial_{\mm{h}} u^1_3|_0^2\mm{d}\tau
 \lesssim   \int_0^t   | \partial_\tau (   \widetilde{\mathcal{N}}^\eta , \mathcal{M}_\tau)|_{0}^2   \mm{d}\tau. \label{2001807012105}
\end{align}

Applying the operator $\mathfrak{D}_{\mf{h}}^{3/2}$ to \eqref{sadfdxsafdsa}, and then following the same process as in the derivation of
\eqref{2001807012105},  we obtain
\begin{align}
  \int_0^t  | \mathfrak{D}_{\mf{h}}^{3/2}\nabla_{\mm{h}}\partial_{\mm{h}} u^1_3|_0^2\mm{d}\tau
 \lesssim    \int_0^t   |    \mathfrak{D}_{\mf{h}}^{3/2} \partial_\tau(\widetilde{\mathcal{N}}^\eta,   \mathcal{M})|_{0} ^2    \mm{d}\tau.\nonumber
\end{align}
Integrating the above integral over $\mathbb{R}^2$, and then adding the resulting estimate to \eqref{2001807012105} yields
\begin{align}
 \int_0^t  | \nabla_{\mm{h}}\partial_{\mm{h}}  u^1_3|^2_{1/2} \mm{d}\tau
 \lesssim    \int_0^t   (   | \Delta_{\mm{h}} u_3^2 |_{1/2}^2+|\partial_t( \mathcal{N}^\eta,\mathcal{M})|_{1/2}^2 )  \mm{d}\tau, \nonumber
\end{align}
which, together with \eqref{201806271120dsfs} and \eqref{2022307252021}, implies
\begin{align}
 \int_0^t | \nabla_{\mm{h}}\partial_{\mm{h}}  u^1_3|^2_{1/2}\mm{d}\tau
 \lesssim  |\mathcal{H}^0|_1^2+   \int_0^t   ( | \Delta_{\mm{h}} u_3^2 |_{1/2}^2 + \sqrt{\mathcal{E}}\mathcal{D} )  \mm{d}\tau. \label{20018070saf121sfafsad05}
\end{align}

(3) Analogously to \eqref{estimforhoedsds1stm} with $(i,\vartheta)=(2,0)$, we can derive from \eqref{s0106pndfsgnnnxx} that
\begin{align}
&\frac{\mm{d}}{\mm{d}t}\int  \left(  \bar{\rho} \partial_\mm{h}^2 \eta^2 \cdot  \partial_\mm{h}^2 u^2
 + \mathcal{U}(\partial_\mm{h}^2 \eta^2 )/2\right) \mm{d}y + \|\sqrt{P'(\bar{\rho}) \bar{\rho}}\partial_{\mm{h}}^2\mm{div}\eta^2\|_0^2\nonumber  \\
& \lesssim   \|   u^2\|^2_{2,0}  +\| \partial_{\mm{h}}^2 \eta^2\|_{1}\|\widetilde{\mathbf{N}}^3\|_{1}
  +|\partial_{\mm{h}}^2 \eta^2|_{1/2}|(\mathbf{N}_1^4,\mathbf{N}_2^4, \llbracket R_P+\mathcal{N}^u \rrbracket ,\mathcal{M})|_{3/2}.\nonumber
\end{align}
We further utilize \eqref{06011711jumpv}, \eqref{201807241645} and the trace estimate to have
 \begin{align}
&\int  \left( \bar{\rho} \partial_\mm{h}^2 \eta^2 \cdot  \partial_\mm{h}^2 u^2  +
\mathcal{U}(\partial_{\mm{h}}^i\partial_\mm{h}^2 \eta^2 ) /2\right) \mm{d}y\lesssim
\|\eta^0\|_3^2+\|u^0\|_2^2 \nonumber \\
&+\int_0^t\left( \|   u^2\|^2_{2,0}
    +\|\eta^2\|_{2,1} (\|\eta\|_2+|\mathcal{M}|_{3/2}+\|\eta\|_3\sqrt{\mathcal{D}})\right)\mm{d}\tau. \label{201811050925}
\end{align}
In addition, by the same manner as in the derivation of \eqref{estimforhoedsds1stnn1524} with $(i,\vartheta)=(2,0)$,
we can deduce from \eqref{s0106pndfsgnnnxx} that
 \begin{align}
 &\frac{1}{2}\frac{\mm{d}}{\mm{d}t} (\|\sqrt{\bar{\rho} } \partial_\mm{h}^2 u^2\|^2_0
+   \|\sqrt{P'(\bar{\rho}) \bar{\rho}}\partial_{\mm{h}}^2\mm{div}\eta^2\|_0^2)+ c\|\partial_\mm{h}^2   u^2 \|_{1}^2\nonumber  \\
&\lesssim \| \partial_{\mm{h}}^2 u^2\|_{1}\|\widetilde{\mathbf{N}}^3\|_{1}+|\partial_{\mm{h}}^2 u^2|_{1/2}|(\mathbf{N}_1^4,\mathbf{N}_2^4, \llbracket R_P+ \mathcal{N}^u \rrbracket,\mathcal{M}) |_{3/2}, \nonumber
\end{align}
which further yields that
\begin{align}
 \| \partial_\mm{h}^2 u^2\|^2_0 +  \int_0^t\|\partial_\mm{h}^2  u^2 \|_{1}^2\mm{d}\tau \lesssim \|\eta^0\|_3^2+\|u^0\|_2^2+\int_0^t\left(\| \eta \|_{2}^2
+|\mathcal{M}|_{3/2}^2+ \sqrt{\mathcal{E} }\mathcal{D}\right)\mm{d}\tau. \label{201811050925xyz}
\end{align}
Using Korn's inequality  and trace estimate, we can derive from  \eqref{201811050925} and \eqref{201811050925xyz}  that
 \begin{align}
 &
\|\eta^2\|^2_{2,1} +  \|  u^2\|^2_{2,0} +  \int_0^t|\nabla_{\mm{h}}^2  u^2 |_{1/2}^2\mm{d}\tau  \nonumber  \\
&\lesssim   \|\eta^0\|_3^2+\|u^0\|_2^2 + \int_0^t\bigg( \|\eta\|_2^2+|\mathcal{M}|_{3/2}^2+
  \|\eta^2\|_{2,1}(\|\eta \|_{2} \nonumber \\
&\qquad +|\mathcal{M}|_{3/2}+ \|\eta\|_3\sqrt{\mathcal{D}})+\sqrt{\mathcal{E} }\mathcal{D} \bigg)\mm{d}\tau ,  \label{xx201808310859}
\end{align}
which, together with  \eqref{20018070saf121sfafsad05}, yields
 \begin{align}
 &\|\eta^2\|_{2,1}^2+   \int_0^t    | \nabla_{\mm{h}}^2 u_3|^2_{1/2} \mm{d}\tau  \nonumber \\
 &  \lesssim  \|\eta^0\|_3^2+\|u^0\|_2^2+ |\mathcal{H}^0|_1^2 + \int_0^t\bigg(\|\eta\|_2^2+  |\mathcal{M}|_{3/2}^2\nonumber \\
 &\qquad+
  \|\eta^2\|_{2,1}(\|\eta \|_{2} +|\mathcal{M}|_{3/2}  + \|\eta\|_3\sqrt{\mathcal{D}})+\sqrt{\mathcal{E} }\mathcal{D}  \bigg)\mm{d}\tau  ,  \nonumber
\end{align}
which yields  \eqref{201807012235}  by further using \eqref{2022307252021} and  Young's inequality.
\hfill $\Box$
\end{pf}

\subsection{Gronwall-type energy inequality}\label{2022307291526}
With the estimates of $(\eta,u)$ in Lemmas \ref{201612132242nx}--\ref{201807251431}, we are in a position to \emph{a prior} derive Gronwall-type energy inequality, which couples with the solution $(\eta^2,u^2)$ of the linear problem \eqref{s0106pndfsgnnnxx} for the case $\vartheta> 0$.
\begin{lem}\label{20230726}   Let $(\eta,u)$ be the solution of the RT  problem \eqref{n0101nnnM} and satisfy \eqref{aprpiosesnew} with  $\delta  \in (0, \iota]$, where $ \iota$  is the constant in Lemma \ref{201809012320}, and thus the following definition makes sense
\begin{align}\label{201811191638}
& d(x_{\mm{h}},t):=\zeta_3((\zeta_{\mm{h}})^{-1}(x_{\mm{h}},t),0,t)\in (h_-,h_+)\mbox{ for each }t\in [0,T].
 \end{align} For $\vartheta>0$, we additionally assume that $d^0:=d|_{t=0}\in H^3(\mathbb{\mathbb{R}}^2)$  and $(\eta^2,u^2)$ is
a solution of the linear problem \eqref{s0106pndfsgnnnxx} with $\mathcal{H}^0\in H^1(\mathbb{R}^2)$. There are an energy functional $\tilde{\mathcal{E}}(t)$ of $(\eta(t),u(t))$, and constants $\delta_1\in (0,\iota ]$, $c>0$
such that, for any $\delta\leqslant \delta_1$,    $(\eta,u)$ enjoys the Gronwall-type energy inequality:
\begin{align}
&  \tilde{\mathcal{E}}(t)+ \vartheta\left(\|  \eta^2  \|_{2,1}^2(t) + c^{-1}|d(t)|_3^2 \right)   +c^{-1}\int_0^t\mathcal{D}(\tau)\mm{d}\tau\nonumber \\
&
 \leqslant c_7\int_0^t (\tilde{\mathcal{E}}(\tau)+\vartheta \|  \eta^2 (\tau) \|_{2,1}^2)\mm{d}\tau +  c \left(\|\eta^0\|_3^2+ \|u^0\|_2^2+ \vartheta|d^0|_{3}^2  +\int_0^t \|(\eta,\vartheta u)\|_{0}^2\mm{d}\tau\right),
\label{2016121521430850}
\end{align}
where   $c_7$ is the constant after \eqref{2202306151933}  in Proposition \ref{growingmodesolneriodic},   $\tilde{\mathcal{E}}(t)$ satisfies
\begin{align}
\label{2018008121027}
 c^{-1}\mathcal{E}(t)\leqslant  \tilde{\mathcal{E}}(t)\leqslant c \mathcal{E}(t) \mbox{ for any }t\in [0,T] ,
\end{align}  and
the constants $\delta_1$, $c$ depend on the domain $\Omega$ and parameters/functions in the RT problem. It should be noted that $\|  \eta^2  \|_{2,1}^2$ and $|d|_3^2$  exist only for $\vartheta>0$.
\end{lem}
\begin{pf}  If we make use of trace estimate, \eqref{ssebdaiseqinM0846} for $0\leqslant i\leqslant 1$, \eqref{201702061418} for $0\leqslant i\leqslant 1$ and \eqref{Lem:0301m0832}, we can infer that there is a constant $c$ such that, for    any sufficiently large constant $\tilde{c}_1\geqslant 1$ and any sufficiently small $\delta$,
\begin{align}
\frac{\mm{d}}{\mm{d}t} \mathcal{E}_1+c^{-1} \mathcal{D}_1\leqslant c \tilde{c}_1 \left(\| \eta\|_{2}^2 +  \sqrt{{\mathcal{E}}} {\mathcal{D}} \right) ,
\label{202216072023726}
 \end{align}
  where we have defined that
 \begin{align*}\mathcal{E}_1:=&\sum_{|\alpha|\leqslant 1}
 \left(\int \bar{\rho}   \partial_{\mm{h}}^{\alpha}\eta \cdot  \partial_{\mm{h}}^{\alpha}u\mm{d}y+\tilde{c}_1\|\sqrt{\bar{\rho}} \partial_{\mm{h}}^{\alpha}u\|^2_0 +
 \tilde{c}_1  \mathcal{I}(  \partial_{\mm{h}}^{\alpha}\eta )
+\mathcal{U}(\partial_{\mm{h}}^{\alpha}\eta) /2\right)   +   \|\sqrt{\bar{\rho}}  u_{t}\|_{0}^2 + \mathcal{I}(u)
\end{align*}
and
$$  \mathcal{D}_1:= \|u_t\|_{1}^2+  \tilde{c}_1 \| u \|_{\underline{1},1}^2
  .  $$
Thanks to  Korn's and Young's inequalities, and the trace estimate, we  have, for    any sufficiently large constant $\tilde{c}_1$,
\begin{equation}  \label{201809081915}
 \| \eta\|_{\underline{1},1}^2+\tilde{c}_1 \| u\|_{\underline{1},0}^2+\|u_t\|_0^2\lesssim \mathcal{E}_1 .
\end{equation}

Next we further derive the estimate for the higher-order normal derivatives of $\eta$. To this purpose, we rewrite \eqref{n0101nn928m}$_1$ as follows:
\begin{equation}   \label{Stokesequson1137}
   \mu \partial_3^2 u_\mm{h}
   =\nabla_\mm{h}(g\bar{\rho}\eta_3 -P'(\bar{\rho})\bar{\rho}\mm{div}\eta )-
   \mu  \Delta_\mm{h} u_\mm{h} - \tilde{\mu} \nabla_\mm{h} \mm{div} u
  +\bar{\rho} \partial_t u_\mm{h}-{\mathbf{N}}^3_\mm{h} =:(\mathbf{F}_1^3,\mathbf{F}_2^3)^{\top}    \end{equation}
and
\begin{align} & (( \mu +\tilde{\mu})\partial_3^2 u_3+ P'(\bar{\rho})\bar{\rho}  \partial_3^2\eta_3)
 = g\bar{\rho} \partial_3\eta_3- \mu \Delta_\mm{h} u_3-\tilde{\mu}\partial_3\mm{div}_{\mm{h}}u_{\mm{h}}\nonumber  \\
&  -(P'(\bar{\rho}))'\bar{\rho}\mm{div}\eta
- P'(\bar{\rho})\bar{\rho}   \partial_3\mm{div}_{\mm{h}}\eta_{\mm{h}}   +\bar{\rho} \partial_t u_3- {\mathbf{N}}^3_3  =:\mathbf{F}_{3}^3.  \label{Stokesequson1}
\end{align}

We can deduce from
 \eqref{Stokesequson1137} that, for $0\leqslant  j\leqslant  1$,
 \begin{equation}\label{etah1502}
 \frac{1}{2}\frac{\mm{d}}{\mm{d}t}\|\sqrt{\mu } \partial_3^{2+j} \eta_\mm{h}\|_{\underline{1-j},0}^2 \lesssim \| \partial_3^{2+j} \eta_\mm{h}\|_{\underline{1-j},0}  \|\partial_3^j\mathbf{F}_{\mm{h}}^3 \|_{\underline{1-j},0}.
 \end{equation}
Similarly, we can also get from  \eqref{Stokesequson1} that
 \begin{align}
 &\frac{1}{2}\frac{\mm{d}}{\mm{d}t} \left\|\sqrt{ \mu +\tilde{\mu}   }\partial_3^{2+j} \eta_3\right\|_{\underline{1-j},0}^2+
 c\|\partial_3^{2+j}  \eta_3 \|_{\underline{1-j},0}^2\nonumber\\
 &\lesssim {\|\partial_3^{2+j}\eta_3\|}_{\underline{1-j},0}\|\partial_3^j\mathbf{F}_{3}^3 \|_{\underline{1-j},0}^2+
 \begin{cases}
 0&\mbox{for }j=0;\\
\|\partial_3^2\eta_3 \|_ 0 {\|\partial_3^{3}\eta_3\|}_{0} &\mbox{for }j=1.
 \end{cases}
 \label{201702091459}
 \end{align}

In addition, by \eqref{06011711jumpv},
 $$ \|\partial_3^j\mathbf{F}^3\|_{\underline{1-j},0}^2\lesssim \|\partial_3^j(\eta,u)\|_{\underline{2-j},1}^2+\|(u_t,\mathbf{N})\|_1^2\lesssim \|\partial_3^j(\eta,u)\|_{\underline{2-j},1}^2+\|u_t\|_{1}^2+\sqrt{\mathcal{E}} \mathcal{D},$$
where $\mathbf{F}^3 =(\mathbf{F}^3_1, \mathbf{F}^3_2, \mathbf{F}^3_3)^{\top}$.
Thus we derive from \eqref{etah1502}, \eqref{201702091459} and Young's inequality  that
 \begin{equation}\label{etah1502nn}
 \frac{\mm{d}}{\mm{d}t}\overline{\|\partial_3^2\eta\|}_{\underline{1},0}^2
\leqslant  \tilde{c}_1^{-1}\| \partial_3^{2} \eta \|_{\underline{1},0}^2  +
c\tilde{c}_1 (\| (\eta,u)\|_{\underline{2},1}^2+\|u_t\|_{1}^2+\sqrt{\mathcal{E}} \mathcal{D}).
 \end{equation}
 and
 \begin{equation}\label{etah1502nnn12}
 \frac{\mm{d}}{\mm{d}t}
\overline{\|\partial_3^3\eta\|}_{0}^2\leqslant \tilde{c}_1^{-1}\| \partial_3^3\eta \|_{0}^2  + c\tilde{c}_1(\| \partial_3^2\eta \|_{\underline{1},0}^2  +\|\eta\|_{\underline{2},1}^2+\| u\|_{\underline{1},2}^2+\|u_t\|_{1}^2+\sqrt{\mathcal{E}} \mathcal{D}),
 \end{equation}
 where we have defined that $ \overline{\|\partial_3^3\eta\|}_{0}^2:=\overline{\|\partial_3^3\eta\|}_{\underline{0},0}^2  $ and
 $$\overline{\|\partial_3^{2+j}\eta\|}_{\underline{1-j},0}^2:=\left\| \left(\sqrt{\mu } \partial_3^{2+j} \eta_\mm{h},\sqrt{ \mu +\tilde{\mu}  }\partial_3^{2+j} \eta_3\right)\right\|_{\underline{1-j},0}^2 \mbox{ for }0\leqslant j\leqslant 1.$$
Moreover,
\begin{equation}
\label{201702131611}
{\|\partial_3^3\eta\|}_{ {0}}
 \lesssim \overline{\|\partial_3^3\eta\|}_{ {0}} \mbox{ and } {\|\partial_3^2\eta\|}_{\underline{1},0}
 \lesssim \overline{\|\partial_3^2\eta\|}_{\underline{1},0} .
 \end{equation}Plugging \eqref{201702071610} and \eqref{201808292010} into \eqref{etah1502nnn12}, we further get
 \begin{equation}\label{etah1502nnn}
\begin{aligned}
 \frac{\mm{d}}{\mm{d}t}
\overline{\|\partial_3^3\eta\|}_{0}^2\leqslant  \tilde{c}_1^{-1}\| \partial_3^3\eta \|_{0}  + c\tilde{c}_1( \| \partial_3^2\eta \|_{\underline{1},0}^2  +
 \| (\eta,u)\|_{\underline{2},1}^2+\|u_t\|_{1}^2+\sqrt{\mathcal{E}} \mathcal{D}).
 \end{aligned}
 \end{equation}

Multiplying \eqref{etah1502nn} and \eqref{etah1502nnn} by $\tilde{c}_1^{-2}$ and $\tilde{c}_1^{-4}$ respectively, and then adding the two resulting inequalities, we deduce that
 \begin{align}
 \frac{\mm{d}}{\mm{d}t}(\tilde{c}_1^{-2} \overline{ \|\partial_3^2\eta\|}_{\underline{1},0}^2  +\tilde{c}_1^{-4} \overline{ \|\partial_3^3\eta\|}_{0}^2 )\lesssim  & \tilde{c}_1^{-1} (\|(\eta,u)\|_{\underline{2},1}^2+ \|u_t\|_1^2+\sqrt{\mathcal{E}} \mathcal{D})\nonumber \\
 &+ \tilde{c}_1^{-3}\|\partial_3^2\eta\|_{\underline{1},0}^2+ \tilde{c}_1^{-5}\|\partial_3^3\eta\|_{0}^2 ,
\label{201702171722} \end{align}
Now we prove \eqref{2016121521430850} by two cases.

(1) Case of $ \vartheta =0$.

We can derive from \eqref{ssebdaiseqinM0846},  \eqref{201702061418} with $(\vartheta,i)=(0,2)$, \eqref{202216072023726}  and \eqref{201702171722}
that
\begin{align}
 & \frac{\mm{d}}{\mm{d}t}\mathcal{E}_2+c^{-1} \mathcal{D}_2 \nonumber \\
 &\leqslant c  (    \tilde{c}_1^{-1}( \|\eta\|_{\underline{2},1}^2+\|u_t\|_1^2)+
 \tilde{c}_1^{-3} { \|\partial_3^2\eta\|}_{\underline{1},0}^2+ \tilde{c}_1^{-5} { \|\partial_3^3\eta\|}_{0}^2
\nonumber \\
 &\quad +| \eta_3|_2^2+\|   u\|^2_{\underline{2},1}+ \tilde{c}_1 ( \| \eta\|_{2}^2+    \sqrt{\mathcal{E}} \mathcal{D})) ,  \label{2022307291921}
 \end{align}
where we have defined that
 \begin{align*}\mathcal{E}_2: =&\mathcal{E}_1+\tilde{c}_1^{-2} \overline{ \|\partial_3^2\eta\|}_{\underline{1},0}^2 +  \tilde{c}_1^{-4} \overline{ \|\partial_3^3\eta\|}_{0}^2 \\
& +\sum_{|\alpha|=2}
 \left(\int \bar{\rho}   \partial_{\mm{h}}^{\alpha}\eta \cdot  \partial_{\mm{h}}^{\alpha}u+\tilde{c}_1\|\sqrt{\bar{\rho}   } \partial_{\mm{h}}^{\alpha}  u\|^2_0  +
 \tilde{c}_1  \mathcal{I}(  \partial_{\mm{h}}^{\alpha}\eta ) +
 \mathcal{U}(\partial_{\mm{h}}^{\alpha}\eta)/2   \right)
\end{align*}
and
$$\mathcal{D}_2:= \tilde{c}_1^{-5 }\|\eta\|_3^2+ \tilde{c}_1 \|u \|_{\underline{2},1}^2+\mathcal{D}_1 .$$

Thanks to \eqref{201702071610},  \eqref{201702131611} and Korn's inequality, we see that $\mathcal{E}_2$ satisfies,  for any sufficiently large $\tilde{c}_1$, \begin{equation}
\label{201809081saf915}
 c^{-1}\tilde{c}_1^{-4}\mathcal{E}\leqslant  \mathcal{E}_2   \leqslant  c\tilde{c}_1\mathcal{E}
  \end{equation} and
 \begin{equation}
  \|\eta\|_{\underline{2},1}^2 +\tilde{c}_1^{-2} { \|\partial_3^2\eta\|}_{\underline{1},0}^2+ \tilde{c}_1^{-4} { \|\partial_3^3\eta\|}_{0}^2\leqslant \mathcal{E}_2. \label{202308031243}
\end{equation}
In addition, by \eqref{20161111195sdfafsd2}, Newton--Leibniz formula and H\"older's inequality,   it is obvious that
\begin{equation}  \label{2018safd09081915}
\|u_t\|_1^2+\tilde{c}_1\|u\|^2_{\underline{2},1}+\tilde{c}_1^{-5}\mathcal{D} \lesssim    \mathcal{D}_2
\end{equation}
and
\begin{align}
|\eta_3|_2^2= \int_{h_+}^0\partial_3 |\eta_3|^2_2\mm{d}y_3
\lesssim \|\eta_3\|_{2,0}\|\eta_3\|_{2,1}.
\label{202237261609}
\end{align}
Thus, exploiting \eqref{2018safd09081915}, \eqref{202237261609} and Young's inequality,  it follows from \eqref{2022307291921} that,  sufficiently large $\tilde{c}_1$,
\begin{align}
 & \frac{\mm{d}}{\mm{d}t}\mathcal{E}_2+c^{-1}( \|u_t\|_1^2+\tilde{c}_1\|u\|^2_{\underline{2},1}+\tilde{c}_1^{-5}\mathcal{D}) \nonumber \\
 &\leqslant c  (  \tilde{c}_1^{-1} \|\eta\|_{\underline{2},1}^2 +\tilde{c}_1^{-3} { \|\partial_3^2\eta\|}_{\underline{1},0}^2+ \tilde{c}_1^{-5} { \|\partial_3^3\eta\|}_{0}^2 + \tilde{c}_1(  \| \eta\|_{2}^2+   \sqrt{\mathcal{E}} \mathcal{D})).\nonumber
 \end{align}

Making use of   \eqref{202308031243} and the interpolation inequality \eqref{201807291850}, we further derive from the above estimate that there exist a positive constant $c$ and a sufficiently large $\tilde{c}_1$ such that, for any sufficiently small $\delta$,
\begin{align}
 \frac{\mm{d}}{\mm{d}t}\mathcal{E}_2+c^{-1}\tilde{c}_1^{-5}\mathcal{D}  \leqslant  c_7\mathcal{E}_2+ c ( \|\eta\|_0^2 +\sqrt{\mathcal{E}} \mathcal{D})  .
 \label{2022307261638}
 \end{align}
Now let $\tilde{\mathcal{E}}= \mathcal{E}_2$, where  $\mathcal{E}_2$ satisfies \eqref{201809081saf915}. Finally we immediately obtain
\eqref{2016121521430850} with $\vartheta=0$ by integrating \eqref{2022307261638} over $(0,t)$, and then using \eqref{201702071610nb} and \eqref{201809081saf915} with $t=0$.

(2) Case of  $\vartheta> 0$.

It follows from \eqref{ssebdaisedsaqinM0dsfadsff846asdfadfad}, \eqref{ssebdaiseqinM0asdfa846}, \eqref{ssebdaisedsaqinM0dsfadsff846sdfaaasdfadfad},  \eqref{2018072033}, \eqref{202216072023726} and \eqref{201702171722}  that, for any sufficiently small $\delta$,
 \begin{align}
 \frac{\mm{d}}{\mm{d}t}\mathcal{E}_3+c^{-1}\mathcal{D}_3\leqslant   &c ( \tilde{c}_1^{-1}(\|(\eta,u)\|_{{2},1}^2+\|u_t\|_1^2)+\tilde{c}_1^{-3} { \|\partial_3^2\eta\|}_{\underline{1},0}^2+\tilde{c}_1^{-5} { \|\partial_3^3\eta\|}_{0}^2 \nonumber \\
 &+  \tilde{c}_1^3  (|   \eta |^2_2 +   \|  (\eta, u)\|^2_{2}+ \sqrt{{\mathcal{E}}} {\mathcal{D}} )) \nonumber
\\
&     + c   (|\partial_{\mm{h}}^2  \eta_3 |_{1/2}  +\tilde{c}_1 |\partial_{\mm{h}}^2  u_3   |_{1/2})   \left(\sqrt{\| \eta\|_{\underline{2},1}\|\eta\|_{1,2}}
    \right.\nonumber \\
    &\quad \left.+ \|  \eta\|_{\underline{2},1} +\| u\|_{\underline{2},1}+ \sqrt{\| u\|_{\underline{2},1}\|u\|_{1,2}}  \right) , \label{2022307229194}
 \end{align}
where we have defined that
$$ \begin{aligned}
\mathcal{E}_3:=&  \mathcal{E}_1+\tilde{c}_1^{-2} \overline{ \|\partial_3^2\eta\|}_{\underline{1},0}^2+ \tilde{c}_1^{-4} \overline{ \|\partial_3^3\eta\|}_{0}^2\\
&+ \tilde{c}_1^2\sum_{|\alpha|= 1} \int_{\mathbb{R}^2}\bigg(
\int   \bar{\rho}  \mathfrak{D}_{\mf{h}}^{3/2 }\partial_{\mm{h}}^{\alpha} \eta \cdot  \mathfrak{D}_{\mf{h}}^{3/2 } \partial_{\mm{h}}^{\alpha}u \mm{d}y
 +\mathcal{U}(\mathfrak{O}_{\mf{h}}^{3/2}\partial_{\mm{h}}^{\alpha} \eta)/2 )\\
 &+\tilde{c}_1 ( \|\sqrt{\bar{\rho} }   \mathfrak{D}_{\mf{h}}^{3/2} \partial_{\mm{h}}^{\alpha} u\|^2_0+
  \mathcal{I}( \mathfrak{D}_{\mf{h}}^{3/2} \partial_{\mm{h}}^{\alpha}  \eta  ))\bigg) \mm{d}\mathbf{h}   +  \sum_{|\alpha|=2}\left( \int\bar{\rho}   \partial_{\mm{h}}^{\alpha}  \eta \cdot  \partial_{\mm{h}}^{\alpha}  u \mm{d}y + \mathcal{U}(\partial_{\mm{h}}^{\alpha} \eta )/2 \right) \\
   &+
   \tilde{c}_1 \left(\|\sqrt{\bar{\rho}  } u\|^2_{2,0}+ \left\|\sqrt{P'(\bar{\rho})\bar{\rho}}\partial_{\mm{h}}^2 \left(\frac{g  \eta _3} {P'(\bar{\rho})} -
 \mm{div}  \eta \right) \right\|_{2,0}^2 \right)
\end{aligned}$$
and
 $$
\mathcal{D}_3:=\|    \mm{div}\eta   \|^2_{2, 0}+ \vartheta  \tilde{c}_1^2|\nabla_{\mm{h}}^2\eta_3|_{1/2}^2 +  \tilde{c}_1^{-5}\|\eta\|_3^2+  \tilde{c}_1 \|     u \|^2_{2, 1}+\mathcal{D}_1 ,
 $$

Thanks to \eqref{201702071610},  \eqref{201702131611}, \eqref{201806291928}, trace estimate and Korn's inequality, we see that $\mathcal{E}_3$ satisfies,  for any sufficiently large $\tilde{c}_1$, \begin{equation}  \label{201809081sasaff915}
 c ^{-1}\tilde{c}_1^{-4}\mathcal{E}\leqslant  \mathcal{E}_3   \leqslant  c \tilde{c}_1^{3}\mathcal{E}
\end{equation}
and
 \begin{equation}
  \|\eta\|_{\underline{2},1}^2 +\tilde{c}_1^{-2} { \|\partial_3^2\eta\|}_{\underline{1},0}^2+ \tilde{c}_1^{-4} { \|\partial_3^3\eta\|}_{0}^2\leqslant \mathcal{E}_3. \label{2308031247}
\end{equation}
In addition, by \eqref{20161111195sdfafsd2},
\begin{equation}  \label{2018sasaffd09081915}
\tilde{c}_1^2|\nabla_{\mm{h}}^2\eta_3|_{1/2}^2+\tilde{c}_1\|u\|^2_{\underline{2},1} + \|u_t\|_1^2  +\tilde{c}_1^{-5}\mathcal{D}  \lesssim    \mathcal{D}_3.
\end{equation}

Using \eqref{201808292010}, \eqref{2018sasaffd09081915} and Young's inequality, we further get from \eqref{2022307229194} that, for any sufficiently large $ \tilde{c}_1 $,
 \begin{align*}
 &\frac{\mm{d}}{\mm{d}t}\mathcal{E}_3+c^{-1}( \tilde{c}_1^2|\nabla_{\mm{h}}^2\eta_3|_{1/2}^2+\tilde{c}_1\|u\|^2_{\underline{2},1} +   \|u_t\|_1^2+\tilde{c}_1^{-5}\mathcal{D}  )\nonumber \\
 &\leqslant   c ( \tilde{c}_1^{-1}\|\eta\|_{2,1}^2 + \tilde{c}_1^{-3}\|\eta\|_{1,2}^2+\tilde{c}_1^{-5} { \|\partial_3^3\eta\|}_{0}^2 \nonumber \\
 &\qquad +\tilde{c}_1^4 (|   \eta |^2_2 +  \|  (\eta, u)\|^2_{2}+  |\nabla_{\mm{h}}^2 u_3|_{1/2}^2 +\sqrt{{\mathcal{E}}} {\mathcal{D}} )).
 \end{align*}

 Making use of   \eqref{202237261609}, \eqref{2308031247}, the interpolation inequality and Young's inequality, we derive from the above inequality that there exist a positive constant $c$ such that, for any sufficiently small $\delta$,
\begin{align}
 \frac{\mm{d}}{\mm{d}t}\mathcal{E}_3+c^{-1} \mathcal{D}  \leqslant  c( \|(\eta,u)\|_0^2+|\nabla_{\mm{h}}^2 u_3|_{1/2}^2+\sqrt{{\mathcal{E}}} {\mathcal{D}} ) +c_7\mathcal{E}_3/2,
 \label{202230726dsaf1638}
 \end{align}
which, together with  \eqref{201809081sasaff915} with $t=0$, yields
\begin{align}
 \mathcal{E}_3+c^{-1} \int_0^t\mathcal{D}  \mm{d}\tau\leqslant   c\mathcal{E}|_{t=0}+  \int_0^t( c( \|(\eta,u)\|_0^2+ |\nabla_{\mm{h}}^2 u_3|_{1/2}^2 +\sqrt{{\mathcal{E}}} {\mathcal{D}} )+c_7\mathcal{E}_3/2)\mm{d}\tau.\nonumber
 \end{align}
 Exploiting the interpolation inequality and \eqref{201702071610nb},  we further deduce from  \eqref{201807012235} and the above inequality that, for any sufficiently small $\delta$,
\begin{align}
& \vartheta\|  \eta^2  \|_{2,1}^2 +\tilde{c}_2\mathcal{E}_3+c^{-1}\int_0^t(\mathcal{D} + \vartheta| \nabla_{\mm{h}}^2 u_3|^2_{1/2} )\mm{d}\tau \nonumber \\
& \leqslant c\left(\|\eta^0\|_3^2+\|u^0\|_2^2 + |\mathcal{H}^0|_1^2 +\int_0^t\|(\eta,u)\|_{0}^2\mm{d}\tau\right)+ \int_0^t ( c_7 \vartheta\|  \eta^2 (t) \|_{2,1}^2 +   c_7\tilde{c}_2\mathcal{E}_3)\mm{d}\tau. \label{2020308111350}
\end{align}

In addition, we can deduce that
\begin{align}   &  |\mathcal{H}|_1\lesssim |d|_3,\label{20181114101015}  \\
& |d(t)|_{3}^2\lesssim |d^0|_3^2+\int_0^t \mathcal{D}(\tau)\mm{d}\tau   \label{20181114101015xx}
\mbox{ for any }t\in [0,T]
\end{align}
(also cf. the derivation \eqref{201811141623} and \eqref{201811141623xx} in Section \ref{2022308081027}).
Here we will omit the proof of the above two estimates, since the auxiliary properties of inverse transform in Lagrangian coordinates are required
and will be introduced in Section \ref{201811191634}.
Consequently,  let $\tilde{\mathcal{E}}= \tilde{c}_2\mathcal{E}_2$, we get \eqref{2016121521430850} from \eqref{2020308111350}--\eqref{20181114101015xx}. This completes the proof of Lemma \ref{20230726}.
\hfill$\Box$ \end{pf}

The local existence of strong solutions to the equations of stratified compressible viscous fluids has been
established, see \cite{jang2016compressible} for example. Similarly to  \cite{jang2016compressible},
we can use a Faedo--Galerkin approximation scheme for the linearized problem and an iterative method to obtain a local-in-time existence
result of a unique strong solution $(\eta,u)$ to the RT  problem \eqref{n0101nnnM},
and a unique global-in-time strong solution $(\eta^2,u^2)$   to the linear problem \eqref{s0106pndfsgnnnxx}
for  given $\mathcal{M}$, $\mathbf{N}_1^4$, $\mathbf{N}_2^4$, $\widetilde{\mathbf{N}}^3$,  $R_P$ and $\mathcal{N}^u$ defined by $(\eta,u)$. Moreover, the strong  solutions also satisfy
the \emph{a priori }estimates in Lemma \ref{20230726}. Since the proof is standard in the  well-posedness theory of PDEs,
and hence we omit its details here, and only state the  well-posedness results, in which the solutions enjoy  the Gronwall-type energy inequality.
\begin{pro} \label{pro:0401nd}
\begin{enumerate}
  \item[(1)] Let  $(\eta^0,u^0)\in H^{3,1/2}_{0,*}\times H_0^2$ and $\zeta^0:=\eta^0+y$. There is a sufficiently small $\delta_2\in (0,\iota)$,
  such that if $(\eta^0,u^0)$ satisfies
\begin{align}
&\sqrt{\|\eta^0\|_3^2+\|{u}^0\|_2^2 }\leqslant   \delta_2,\  \mbox{the (necessary) compatibility condition \eqref{2022212081221}},
\end{align}
and additionally
\begin{align}
 d^0:=\zeta_3^0((\zeta_{\mm{h}}^0)^{-1}(x_{\mm{h}}),0)\in H^3(\mathbb{R}^2)\mbox{ if }\vartheta>0 , \label{202308031251}
\end{align}
then there is a local existence time $T_{\mm{loc}}>0$, depending on $\delta_2$, the domain and the known parameters/functions in the RT  problem \eqref{n0101nnnM},
and a unique local-in-time strong  solution
$(\eta, u)\in C^0([0, T_{\mm{loc}}),H^{3,1/2}_{0,*}\times H^2 )$  to the RT  problem, where $(\eta,u)$
enjoys the following regularity
\begin{align*}
  u_t \in C^0([0,T_{\mm{loc}}),L^2    )\mbox{ and }(u,u_t )\in L^2((0,T_{\mm{loc}}),H^3\times H^1 ) .
\end{align*}
Moreover, if $\vartheta> 0$, $|\mathcal{H}|_1 \lesssim |d|_3$, where
 $d:=\zeta_3((\zeta_{\mm{h}})^{-1}(x_{\mm{h}},t),0,t)$ makes sense and belongs to $C^0([0,T_{\mm{loc}}),$ $H^3(\mathbb{R}^2))$.
  \item[(2)] For the case $\vartheta>0$, we further use the above strong  solution $(\eta,u)$ to define $\mathbf{N}_1^4$, $\mathbf{N}_2^4$, $\widetilde{\mathbf{N}}^3$,  $R_P$, $\mathcal{N}^u$ and $\mathcal{M}$ with $\mathcal{M}|_{t=0}= \mathcal{H}|_{t=0}\in H^1(\mathbb{R}^2)$ and the linear problem \eqref{s0106pndfsgnnnxx} possesses a unique global-in-time strong solutions
$(\eta^2, u^2) $, where $(\eta^2, u^2)$ enjoys the regularity as well as $(\eta, u)$.
  \item[(3)] In addition, if $(\eta,u)$ further satisfies
$$\sup_{0\leqslant t\leqslant T}\sqrt{\|\eta(t)\|_3^2 +\|u(t)\|_2^2}\leqslant \delta_1\mbox{ for some }T<T_{\mm{loc}},$$
 then the solution $(\eta, u)$ enjoys the Gronwall-type energy inequality \eqref{2016121521430850} and the equivalent estimate \eqref{2018008121027} on $[0,T]$.
\end{enumerate}
\end{pro}
\begin{rem}By the smallness condition $\|\eta^0\|_3\leqslant   \delta_2< \iota$ and  Lemma \ref{201809012320}, $\zeta^0_{\mm{h}}(y_{\mm{h}},0):\mathbb{R}^2 \to \mathbb{R}^2$
   is a homeomorphism mapping, and thus  $d^0:=\zeta_3((\zeta_{\mm{h}})^{-1}(x_{\mm{h}}),0) $ in \eqref{202308031251} makes sense; moreover $d^0\in (h_-,h_+)$.
In addition,  we do not state the necessary compatibility conditions in the second assertion in Proposition \ref{pro:0401nd} as in the first assertion,
since it is easy to observe that the initial data in the linear problem \eqref{s0106pndfsgnnnxx} automatically satisfy the necessary compatibility condition, i.e. \eqref{2022307312050}.
\end{rem}

\section{Construction of initial data for the nonlinear problem} \label{2022306101652}

 Let $(\tilde{\eta}^0,\tilde{u}^0) $ come from \eqref{2022306101947} in Proposition \ref{growingmodesolneriodic} and $c_7$ be the constant after \eqref{2202306151933}.
By the definition of $\tilde{\eta}^0$ in \eqref{2022306101947} and \eqref{2022306102034}--\eqref{x2022306102034}, there exist positive constants $c_4$, $c_5$ and $c_8$ such that $(\tilde{\eta}^0,\tilde{u}^0) $ enjoys the estimates, for any  $n\geqslant c_4$,
 \begin{align}
& \sum_{  \beta_1+ \beta_2+ \beta _3\leqslant 4,\ 1\leqslant \beta _1+\beta _2}\|\partial_{1}^{\beta _1}\partial_{2}^{\beta _2} \chi_{n,n} \partial_{3}^{\beta _2}(\tilde{\eta}^0,\tilde{u}^0)\|_{0}^2+n^{-1}\|\chi_{n,n}  (\tilde{\eta}^0,\tilde{u}^0)\|_4^2 \leqslant  c_5^2n,
\label{appesimtsofu1857}\\
&2c_8 n\leqslant \min_{\omega=\tilde{\eta}^0,\tilde{u}^0}\{  \|\chi_{n,n} \omega_{\mm{h}}\|_{0}, \|\chi_{n,n} \omega_3\|_{0} ,
|\chi_{n,n}\omega_3|_{0}\}.  \label{2023061632013}
 \end{align}

 From now on, for any given $\delta>0$, the integer $n$ always satisfies
  \begin{align}
 \label{2022308081221}
n\geqslant\max\{c_4, \delta^{-2}\}.
  \end{align}
 Let
\begin{equation}\label{0501}
\left(\eta^\mm{a}, u^\mm{a}\right)=\delta e^{{c_7  t}}\chi_{n,n}  (\tilde{\eta}^0,\tilde{u}^0)/n\in (H_0^1\cap H^4)^2,
\end{equation}
Then the approximate solution $\left(\eta^\mm{a}, u^\mm{a}\right)$ satisfies the estimate by \eqref{appesimtsofu1857} and the form \eqref{0501}:  for any $0\leqslant  j\leqslant  2$ and  for any $t \geqslant 0$,
\begin{align}
\| \partial_{t}^j(\eta^\mm{a}, u^\mm{a})\|_3=c_7^j \delta e^{c_7 t} \|\chi_{n,n}(\tilde{\eta}^0,\tilde{u}^0)\|_3/n
\leqslant c_5 c_7^j \delta e^{c_7 t}, \label{appdesimtsofu1857}
\end{align}
and the following relations
\begin{equation}\label{linearizsded} \begin{cases}
\eta_t^\mm{a}=u^\mm{a} &\mbox{in } \Omega,\\[1mm]
\bar{\rho} u_t^\mm{a} = g\bar{\rho}(\mm{div}\eta^\mm{a} \mathbf{e}^3- \nabla \eta_3^\mm{a})+
\mm{div}\Upsilon( \eta^\mm{a},u^\mm{a} ) +\delta e^{{c_7 t}} \mathbf{R}^1( \tilde{\eta}^0 , \tilde{u}^0 ) /n&\mbox{in }  \Omega,\\[1mm]
 \llbracket   u^\mm{a}  \rrbracket  =\llbracket   \eta^\mm{a}  \rrbracket  =0,\  \llbracket     \Upsilon( \eta^\mm{a},u^\mm{a} )\mathbf{e}^3   \rrbracket+ \vartheta \Delta_{\mm{h}}\eta^\mm{a}_3 \mathbf{e}^3=\delta e^{{c_7 t}} \mathbf{R}^2( \tilde{\eta}^0 , \tilde{u}^0 )/n&\mbox{on }\Sigma,\\
(\eta^\mm{a}, u^\mm{a})=0 &\mbox{on }\partial\Omega\!\!\!\!-,\\
 (\eta^\mm{a},  u^\mm{a})|_{t=0}=\delta   \chi_{n,n}  (\tilde{\eta}^0,\tilde{u}^0)/n &\mbox{in }  \Omega,
\end{cases}
\end{equation}
where we have defined that
\begin{align*}
\mathbf{R}^1( \tilde{\eta}^0 , \tilde{u}^0 ):=&  g\bar{\rho}(\tilde{\eta}^0  \cdot\nabla \chi_{n,n} \mathbf{e}^3- \tilde{\eta}^0_3\nabla \chi_{n,n})+\mm{div}  ( P'(\bar{\rho})\bar{\rho}\tilde{\eta}^0  \cdot \nabla \chi_{n,n} \mathbb{I} \\
& +  \mu (\nabla \chi_{n,n} (\tilde{u}_0)^{\top }+  \tilde{u}^0 (\nabla \chi_{n,n})^{\top } )+\left(\varsigma-{2\mu}/{3}\right)\tilde{u}^0 \cdot \nabla \chi_{n,n} \mathbb{I}), \\
\mathbf{R}^2( \tilde{\eta}^0 , \tilde{u}^0 ):=&  -\llbracket \mu (\nabla \chi_{n,n} (\tilde{u}^0)^{\top }+ \tilde{u}^0(\nabla \chi_{n,n})^{\top } )+\left(\varsigma-{2\mu}/{3}\right) \tilde{u}^0 \cdot \nabla \chi_{n,n} \mathbb{I}    \rrbracket\nonumber \\
 &- \vartheta ( 2 \nabla_{\mm{h}} \chi_{n,n}\cdot \nabla_{\mm{h}}\tilde{\eta}^0_3  + \Delta_{\mm{h}}\chi_{n,n} \tilde{\eta}^0_3)\mathbf{e}^3.
\end{align*}
Thanks to \eqref{appesimtsofu1857} and trace estimate, we easily estimate that
\begin{align}
&  n^{-1/2}(\| \mathbf{R}^1( \tilde{\eta}^0 , \tilde{u}^0 )\|_{1}+| \mathbf{R}^2( \tilde{\eta}^0 , \tilde{u}^0 )|_{1/2})\nonumber \\
& \lesssim n^{-1/2}( \| (\tilde{\eta}^0  \cdot\nabla \chi_{n,n} ,\tilde{\eta}^0_3\nabla \chi_{n,n}, \nabla_{\mm{h}} \chi_{n,n}\cdot \nabla_{\mm{h}}\tilde{\eta}^0_3, \Delta_{\mm{h}}\chi_{n,n} \tilde{\eta}^0_3)\|_1\nonumber \\
&\qquad+\|(\tilde{\eta}^0  \cdot \nabla \chi_{n,n}  ,\nabla \chi_{n,n} (\tilde{u}_0)^{\top }, \tilde{u}^0 \nabla (\chi_{n,n})^{\top } ,\tilde{u}^0 \cdot \nabla \chi_{n,n} )\|_2 \lesssim  1.
\label{20222302522222}
\end{align}

Unfortunately, the initial data $\left(\eta^\mm{a}, u^\mm{a}\right)|_{t=0}$ does not satisfy the initial compatibility jump condition  of the (nonlinear) RT  problem \eqref{n0101nnnM} in general.
Therefore, next we modify the initial data $\left(\eta^\mm{a}, u^\mm{a}\right)|_{t=0}$.
\begin{pro}\label{lem:modfied} Let $(\tilde{\eta}^0,\tilde{u}^0)$ come from \eqref{2022306101947} in Proposition \ref{growingmodesolneriodic}.
 Then we can look for a constant  ${\delta}_3\in (0,\iota]$   depending on $(\tilde{\eta}^0,\tilde{u}^0)$
 such that for any given $\delta \in (0, {\delta}_3]$ and for any given $n$ satisfying \eqref{2022308081221}, there exist  an error function $u^\mm{r}\in  H_0^1\times H^2$   such that
\begin{enumerate}
  \item[(1)] the modified initial data
  \begin{equation}\label{mmmode04091215}(  {\eta}_0^\delta,{u}_0^\delta): =\delta
  n^{-1}   \chi_{n,n} (\tilde{\eta}^0,\tilde{u}^0)+ ( 0, \delta^2 u^\mm{r})\in (H_0^1\cap H^3)\times (H_0^1\cap H^2)
\end{equation}
satisfies the following compatibility jump condition of the RT  problem:
\begin{align}
\label{202307281417}
 \llbracket P(\bar{\rho}J^{-1}_{0,\delta})\mathbb{I}-\mathbb{S}_{\mathcal{A}_0^\delta}( u_0^\delta)    \rrbracket  \mathcal{A}_0^\delta\mathbf{e}^3 = \vartheta \mathcal{H}_0^\delta \mathcal{A}_0^\delta\mathbf{e}^3\mbox{ on }\Sigma,
 \end{align}
 where $J_{0,\delta}$, $ \mathcal{A}_{0}^\delta$ and $\mathcal{H}_0^\delta $ are defined by $J$, $\mathcal{A}$ and $\mathcal{H}^0$ with ${\eta}_0^\delta$ in place of $\eta$.
\item[(2)] $u^\mm{r} $ satisfies the following estimate:
\begin{equation}
\label{201702091755}
 \|u^\mm{r}\|_2 \leqslant  c_6,
 \end{equation}
where the constant $c_6\geqslant  1$ depends on the known physical parameters/functions, but is independent of $\delta$ and $n$.
\end{enumerate}\end{pro}
\begin{pf}Let
\begin{align*}
 \mathbf{F}^4({\eta}_0^\delta,{u}_0^\delta) = ( {\mathbf{N}}_1^4 ({\eta}_0^\delta,{u}_0^\delta), {\mathbf{N}}_2^4 ({\eta}_0^\delta,{u}_0^\delta), {\mathcal{N}} ({\eta}_0^\delta,{u}_0^\delta))^{\top } ,
\end{align*}
where $({\eta}_0^\delta,{u}_0^\delta)$ enjoys the mode \eqref{mmmode04091215}, and   the definitions of   ${\mathbf{N}}_1^4$, ${\mathbf{N}}_2^4$ and $ {\mathcal{N}}$ can be found in \eqref{2saf022306161627}  and \eqref{202306161628}.
Noting that \eqref{n0101nnnM}$_3$ is equivalent to \eqref{n0101nn1928M}$_2$ for $\|\eta\|_3\leqslant \delta$ with sufficiently small $\delta$, and  $\delta  \chi_{n,n}(\tilde{\eta}^0,\tilde{u}^0)/n$ belongs to $(H_0^1\cap H^4)^2$ and satisfies
 $$ \llbracket    \Upsilon(\chi_{n,n} \tilde{\eta}^0,\chi_{n,n} \tilde{u}^0  )  \rrbracket   \mathbf{e}^3+ \vartheta \Delta_{\mm{h}}(\chi_{n,n} \tilde{\eta}^0_3 ) \mathbf{e}^3  =   \mathbf{R}^2( \tilde{\eta}^0 , \tilde{u}^0 )
 \mbox{ on }\Sigma,$$
 therefore $u^\mm{r}$ in the mode \eqref{mmmode04091215} will enjoy the following conditions
 \begin{equation}
 \label{201702061320}
\begin{cases}
 \llbracket  u^\mm{r}  \rrbracket =0  &\mbox{on }\Sigma, \\
  \llbracket  \mathbb{S}(u^\mm{r} ) \mathbf{e}^3  \rrbracket  =\delta^{-2} \mathbf{F}^4 ({\eta}_0^\delta,{u}_0^\delta) -\delta^{-1}  \mathbf{R}^2( \tilde{\eta}^0 , \tilde{u}^0  )/n &\mbox{on }\Sigma, \\
u^\mm{r} =0&\mbox{on }\partial\Omega\!\!\!\!\!-,
\end{cases}\end{equation}
if $( {\eta}_0^\delta,{u}_0^\delta)$  given by \eqref{mmmode04091215} satisfies the compatibility jump condition
\eqref{202307281417} with sufficiently small $\delta$.

To look for such $u^\mm{r}$ satisfying \eqref{201702061320}, we consider the following stratified  Lam\'e problem for given $w\in H^2$:
\begin{equation}
\label{201702051905}
\begin{cases}
 \mu\Delta u+(\varsigma+\mu/3)\nabla \mm{div}u=0&\mbox{ in }  \Omega, \\[1mm]
  \llbracket  u  \rrbracket =0,&\\
    \llbracket  \mathbb{S}(u^\mm{r} )\mathbf{e}^3  \rrbracket   =\delta^{-2} \mathbf{F}^4 (\delta
  n^{-1}   \chi_{n,n}  \tilde{\eta}^0 ,\delta n^{-1}   \chi_{n,n} \tilde{u}^0+\delta^2 w)   -\delta^{-1}  \mathbf{R}^2( \tilde{\eta}^0 , \tilde{u}^0  )/n  &\mbox{ on }\Sigma, \\
u=0 &\mbox{ on }\partial\Omega\!\!\!\!\!-.
\end{cases}\end{equation} In view of the   theory of stratified Lam\'e problem in Lemma \ref{2022306231855}, there exists a solution $u$ to \eqref{201702051905}; moreover
\begin{equation}
\label{Ellipticestimate0839}
\|u\|_{2}\lesssim | ( \delta^{-2}\mathbf{F}^4 (\delta
  n^{-1}   \chi_{n,n}  \tilde{\eta}^0 ,\delta n^{-1}   \chi_{n,n}\tilde{u}^0+\delta^2 w)    ,\delta^{-1} \mathbf{R}^2( \tilde{\eta}^0 , \tilde{u}^0)/n)|_{1/2}.
\end{equation}
Following the arguments of \eqref{201807241645} and \eqref{060117sdfa34}, we have, for sufficiently small $\delta$,
  \begin{align}
& | \mathbf{F}^4 (\delta
  n^{-1}   \chi_{n,n}  \tilde{\eta}^0 ,\delta n^{-1}   \chi_{n,n}\tilde{u}^0+\delta^2 w)   |_{1/2}
 \nonumber \\
      &\lesssim \delta^2 \|n^{-1}   \chi_{n,n}  \tilde{\eta}^0 \|_{3} (\|n^{-1}   \chi_{n,n}  \tilde{\eta}^0 \|_3+\|n^{-1}   \chi_{n,n}  \tilde{u}^0+ \delta w\|_2)\nonumber \\
      &\lesssim \delta^2     (1+ \delta\| w\|_2),
\label{201702102217}
\end{align}
where we have used \eqref{appesimtsofu1857}  in the last inequality.

By \eqref{2022308081221}, we see that
\begin{align} n^{-1/2} \leqslant  \delta.
\label{202308081217}
\end{align}
Thanks to \eqref{20222302522222}, \eqref{201702102217} and \eqref{202308081217}, we can get from \eqref{Ellipticestimate0839} that
\begin{equation}
\label{201702141356}
\|u\|_{2}\leqslant  c_6 ( 1+\delta\|w\|_2)/2
\end{equation} for some constant $c_6 $.
 Therefore, one can construct an approximate function sequence $\{u_\mm{r}^m\}_{m=1}^\infty$, such that, for any $m\geqslant  2$,
\begin{equation}
\label{2017020519050845}
\begin{cases}
 \mu\Delta u_\mm{r}^m+(\varsigma+\mu/3)\nabla \mm{div}u_\mm{r}^m=0&\mbox{ in }  \Omega, \\[1mm]
  \llbracket  u_\mm{r}^m  \rrbracket =0,&\\   \llbracket \mathbb{S}(u_\mm{r}^m)\mathbf{e}^3\rrbracket= \delta^{-2}( \mathbf{F}^4 (\delta
  n^{-1}   \chi_{n,n} \tilde{\eta}^0,\delta n^{-1}   \chi_{n,n} \tilde{u}^0+\delta^2 u_\mm{r}^{m-1} )  + \delta^{-1} \mathbf{R}^2( \tilde{\eta}^0 , \tilde{u}^0 )/n   &\mbox{ on }\Sigma, \\
u_\mm{r}^m=0 &\mbox{ on }\partial\Omega\!\!\!\!\!-
\end{cases}\end{equation}
and $\|u_\mm{r}^1\|_2\leqslant  c_6 $.  Moreover, by \eqref{201702141356}, one has
$$ \|u_\mm{r}^m\|_{2}\leqslant   c_6 (1+\delta\|u_\mm{r}^{m-1}\|_{2})/2$$
for any $m\geqslant  2$, which implies that  \begin{equation}
 \label{201702090854}
 \|u_\mm{r}^m\|_{2}\leqslant  c_6
 \end{equation} for any $n$ satisfying \eqref{2022308081221}, and for any sufficiently small $\delta\leqslant  1/c_6$.

 Next we further show that  $u_\mm{r}^m$ is a Cauchy sequence in $H^2$. We define that  $u_{\mm{r}}^{m,\mm{d}}:=u_{\mm{r}}^m-u_{\mm{r}}^{m-1}$,
\begin{align*}\mathbb{D} _{\mm{h},0}^{m,\delta} :=&\delta^2(\llbracket(\mathbb{S}  (u_{\mm{r}}^{m,\mm{d}})   \mathbf{e}^3)\cdot\tilde{\mathbf{n}}_0^\delta \mathbf{n}_0^\delta   + (\mathbb{S}(u_{\mm{r}}^{m,\mm{d}}) \mathbf{e}^3 ) \cdot \mathbf{e}^3\tilde{\mathbf{n}}_0^\delta    \rrbracket
\\
                    &-  \Pi_{ \mathbf{n}_0^\delta } \llbracket   \mathbb{S}   (u_{\mm{r}}^{m,\mm{d}})   (J_0^\delta{\mathcal{A}}^\delta_0\mathbf{e}^3- \mathbf{e}^3)+
  \mathbb{S}_{\tilde{\mathcal{A}}_0^\delta}(u_{\mm{r}}^{m,\mm{d}})J^\delta \mathcal{A}_0^\delta\mathbf{e}^3  \rrbracket)_{\mm{h}}
\end{align*}
and
$$  \mathbb{D}_{3,0}^{m,\delta}           :=-\delta^2(
      \mathbb{S}_{\tilde{\mathcal{A}}^\delta_0} (u_{\mm{r}}^{m,\mm{d}}) \tilde{\mathbf{n}}_0^\delta\cdot \mathbf{n}_0^\delta+ (\mathbb{S} (u_{\mm{r}}^{m,\mm{d}} )    \tilde{\mathbf{n}}_0^\delta )\cdot  \mathbf{n}_0^\delta +( \mathbb{S} (u_{\mm{r}}^{m,\mm{d}})  \mathbf{e}^3)\cdot \tilde{\mathbf{n}}_0^\delta),
$$
where $J^\delta_0$,   ${\mathcal{A}}^\delta_0$, $\tilde{\mathcal{A}}^\delta_0$, $\mathbf{n}^\delta_0$ resp. $\tilde{\mathbf{n}}^\delta_0$ are defined as $J$, ${\mathcal{A}}$,  $\tilde{\mathcal{A}}$, $\mathbf{n} $ resp. $\tilde{\mathbf{n}}$ by with $\delta n^{-1} \chi_{n,n}\tilde{\eta}^0$ in place of $\eta$, and $\mathbb{S}_{\tilde{\mathcal{A}}^\delta} (u_{\mm{r}}^{m,\mm{d}} )$ is defined as $\mathbb{S}_{ \mathcal{A}} (u )$ in \eqref{2022306090918} with $\tilde{\mathcal{A}}^\delta_0$  resp.  $u_{\mm{r}}^{m,\mm{d}} $ in place of ${\mathcal{A}}$  resp. $u$.

Noting that
$$
 \begin{cases}
 \mu\Delta u_{\mm{r}}^{m+1,\mm{d}} +(\varsigma+\mu/3)\nabla \mm{div}u_{\mm{r}}^{m+1,\mm{d}} =0&\mbox{in }  \Omega, \\[1mm]
  \llbracket  u_{\mm{r}}^{m+1,\mm{d}}  \rrbracket =0,\  \llbracket \mathbb{S}  (u_{\mm{r}}^{m+1,\mm{d}}  ) \mathbf{e}^3\rrbracket=  \delta^{-2}\mathbb{D}^{m,\delta}_0 &\mbox{on }\Sigma, \\
u_{\mm{r}}^{m+1,\mm{d}}  =0 &\mbox{on }\partial\Omega\!\!\!\!\!-,
\end{cases}$$
thus we have
\begin{equation}
\label{Ellipticestimate0837}
\|u_{\mm{r}}^{m+1,\mm{d}} \|_{2}\lesssim \delta^{-2} |\mathbb{D}^{m,\delta}_0 |_{1/2},
\end{equation}
where $\mathbb{D}^{m,\delta}_0=((\mathbb{D}^{m,\delta}_{\mm{h},0})^{\top}, \mathbb{D}^{m,\delta}_{3,0})$.

In addition, similarly to \eqref{201702102217}, it is easy to estimate that, for sufficiently small $\delta$,
$$|\mathbb{D}^{m,\delta}_0 |_{1/2}\lesssim \delta^3\| u_{\mm{r}}^{m,\mm{d}}\|_2.$$
  Putting the above estimate into \eqref{Ellipticestimate0837} yields
$$
\| u_{\mm{r}}^{m+1,\mm{d}} \|_{2}\lesssim \delta\|u_{\mm{r}}^{m,\mm{d}}\|_2,
$$which presents that $\{u_\mm{r}^m\}_{m=1}^\infty$ is a Cauchy sequence in $H^2$  by choose a sufficiently small $\delta$.
Consequently, we can  use a compactness argument to get
a limit function $u^{\mm{r}}$ which solves \eqref{201702061320} by \eqref{2017020519050845}. Moreover $u^\mm{r}$ satisfies \eqref{201702091755}  by \eqref{201702090854} and the strong convergence of $\{u_\mm{r}^m\}_{m=1}^\infty$.   \hfill$\Box$
\end{pf}

\section{Error estimates}\label{20230721}

This section is devoted to the derivation of error estimates between the solutions of the (nonlinear)  RT problem {\eqref{n0101nnnM}} and the solutions
of the linearized RT problem  \eqref{linearizedxx}. To start with, let us introduce the estimate:
\begin{equation}
\label{201811201409}
|\phi|_{3}\leqslant |\phi|_{7/2}\leqslant \tilde{c}_3( \|\chi\|_4 ) \|\chi\|_4  \mbox{ for any }\chi\in H^{4,1/2}_{0,*}\mbox{ satisfying }\|\chi\|_4^2\leqslant 1,
\end{equation}
where $\phi(x_\mm{h}):=\tilde{\chi}_3((\tilde{\chi}_{\mm{h}})^{-1}(x_{\mm{h}}),0)$, $\tilde{\chi}:=\chi+y$, and \emph{the positive constant $\tilde{c}_3( \|\chi\|_4 ) $
is increasing with respect to $\|\chi\|_4$}. Please refer to \eqref{201809292130} in Section  \ref{2022308031213} for a derivation of \eqref{201811201409}.

Then we define
 \begin{align}
 \tilde{c}_4:=  c_5+c_6+\sqrt{\vartheta}\tilde{c}_3(c_5+c_6)+1\mbox{ and }c_3:= \min\{\delta_1/\tilde{c}_4 ,\delta_2/2\tilde{c}_4, \delta_3  \} \leqslant 1. \label{2022307282134}
 \end{align}
 Frow now on, we always assume  $\delta$ satisfies
\begin{align}
\label{2022309091827}
 \delta\in (0, c_3].
 \end{align}
Making use of \eqref{appesimtsofu1857}, \eqref{mmmode04091215}, \eqref{201702091755}, \eqref{201811201409}, Lemma \ref{201809012320}  and Minkowski's inequality in discrete form, we have
$$\sqrt{\| {\eta}_0^\delta\|_3^2+\|{u}_0^\delta\|_2^2+ \vartheta |d^\delta_0|_3^2}\leqslant \tilde{c}_4 \delta<  \delta_2  \in (0,\iota) $$
and
$$(\eta_0^\delta,u_0^\delta)\in  H_{0,*}^{3,1/2}\times (H^1_0\cap H^2), $$
where   $(\eta_0^\delta,u_0^\delta) $ is constructed by Proposition \ref{lem:modfied} with the condition $\delta\leqslant \delta_3$,  and $d^\delta_0(x_\mm{h}):=\mathbf{e}^3\cdot \zeta^\delta_0((\zeta_{\mm{h}}^0)^{-1}(x_{\mm{h}}),0)$ with $\zeta^\delta_0:=\eta^\delta_0+y $.
Hence, in view of Proposition  \ref{pro:0401nd},  there exists a unique local (nonlinear) solution
$(\eta, u )\in C^0([0,T_{\mm{loc}}),H^{3,1/2}_{0,*}\times H^2)$ to the RT  problem emanating
from the initial data $(  {\eta}_0^\delta,{u}_0^\delta)$; moreover
\begin{align}
\label{2022307262218}
\lim_{t\to T^-_{\mm{loc}}}(\|\eta (t)\|_3^2+ \| u (t)\|_2^2)>\delta_2^2.
\end{align}

Now we estimate the error between the solution $(\eta,u )$ and the solution $(\eta^\mm{a},u^{\mm{a}})$ provided by \eqref{0501}.
To this purpose, we define an error function
$(\eta^{\mathrm{d}}, u^{\mathrm{d}}):=(\eta, u)-(\eta^\mm{a},u^{\mm{a}})$.
Then, we can establish the following estimate of the error function.

\begin{pro}\label{lem:0401}
Let $\delta $ and $n$ satisfy \eqref{2022308081221} and \eqref{2022309091827}. For any given positive constant  $\beta$,  we  assume that,  for any $t\in [0,  T]$, \begin{align}
&\label{201702100940}
 \sqrt{ \|\eta\|_3^2+
\| u\|_2^2+
\| u_t\|_0^2  +\|u\|_{L^2((0,t),H^3)}^2 +\|u_\tau\|_{L^2((0,t),H^1)}^2}\leqslant \beta \delta e^{c_7 t},\\
&\label{201702100940n}\delta e^{c_7 t}\leqslant  1\mbox{ and }\|\eta(t)\|_3
  \leqslant  \delta_1,
\end{align}
then
there exists a constant $c $ such that   for  any $t\in [0,T]$,
\begin{equation}\label{ereroe}
\begin{aligned}
  \| (\eta^{\mm{d}},u^{\mm{d}})\|_1+\| u_t^{\mm{d}} \|_0+|( \eta^{\mm{d}}, u^{\mm{d}})|_{1/2}   \leqslant   c\sqrt{\delta^3e^{3c_7 t}},
 \end{aligned}  \end{equation}
\end{pro}
\begin{rem}Here and in the proof of Proposition \ref{lem:0401}, the generic constants $ c$ not only may
depend on the domain $\Omega$, and other known physical parameters/functions in the RT problem \eqref{n0101nnnM},  \emph{but also increasingly depend  on the given constant $\beta$}.
\end{rem}
\begin{pf}
We will break up the proof into four steps.

(1) \emph{First we assert that, for any $w\in H^1_0$},
\begin{align}
& g \llbracket\bar{\rho} \rrbracket |w|_0^2-\mathcal{I}(w)   \leqslant  {\Lambda^2} \|\sqrt{ \bar{\rho}} w\|_0^2+\Lambda\mathcal{U}(w) , \label{202230as60508057}
\end{align}\emph{
where $\Lambda$ is defined by \eqref{2308031344}.}
Next we verify the above assertion.

Let $\hat{f}$ be the horizontal Fourier transform of $f$, see \eqref{hftx} the definition. By Fubini's and Parseval's theorems, we have that $\hat{f} \in L^2 $ and
\begin{equation}\label{parseval}
\int  |{f(y)}|^2 \mm{d}y =\int_{\mathbb{R}^2}\int_{h_-}^{h^+} |{\hat{f}(\xi,x_3)}|^2 \mm{d}y_3\mm{d}\xi .
\end{equation}
Applying the horizontal Fourier transform  to $w$, we have
\begin{equation}
\label{202230731saf2108}
 \mm{div}{w} =\partial_3  \hat{w}_3-\mm{i} \xi_1 \hat{w}_1-\mm{i} \xi_2  \hat{w}_2
\end{equation}
and
\begin{equation}\nonumber
 \mathbb{D} w = \begin{pmatrix}
       -2\mm{i}\xi_1  \hat{w}_1 & -\mm{i}(\xi_1\hat{w}_2+ \xi_2  \hat{w}_1) & \partial_3 \hat{w}_1-\mm{i} \xi_1  \hat{w}_3    \\
      -\mm{i}(\xi_1\hat{w}_2+ \xi_2  \hat{w}_1)  & - 2 \mm{i} \xi_2\hat{w}_2  & \partial_3 \hat{w}_2  -  \mm{i} \xi_2  \hat{w}_3 \\
  \partial_3  \hat{w}_1    - \mm{i} \xi_1  \hat{w}_3  & \partial_3 \hat{w}_2 -\mm{i} \xi_2  \hat{w}_3    & 2\partial_3  \hat{w}_3
      \end{pmatrix}  .
\end{equation}
In particular, we can further compute out that
\begin{align}
| \mathbb{D} w|^2/2=&
 2(|\xi_1  \hat{w}_1 |^2 +|\xi_2\hat{w}_2|^2+ |\partial_3  \hat{w}_3|)+|\mm{i}(\xi_1\hat{w}_2+ \xi_2  \hat{w}_1)|^2\nonumber \\
&+| \partial_3 \hat{w}_1 -\mm{i} \xi_1  \hat{w}_3|^2+| \partial_3 \hat{w}_2-  \mm{i} \xi_2  \hat{w}_3 |^2 \nonumber \\
=&
 |\mm{i}(\xi_1\hat{w}_2- \xi_2  \hat{w}_1)|^2+| \partial_3 \hat{w}_1-\mm{i} \xi_1  \hat{w}_3|^2+| \partial_3 \hat{w}_2-  \mm{i} \xi_2  \hat{w}_3 |^2   \nonumber \\
&+|\partial_3  \hat{w}_3+\mm{i} \xi_1 \hat{w}_1+\mm{i} \xi_2  \hat{w}_2|^2+|\partial_3  \hat{w}_3-\mm{i} \xi_1 \hat{w}_1-\mm{i} \xi_2  \hat{w}_2|^2 . \label{23080314191}
\end{align}

Exploiting \eqref{parseval} and \eqref{202230731saf2108}, we have
\begin{align}
  \|\sqrt{\bar{\rho}}w\|_0^2= & \int_{\mathbb{R}^2}\int_{h_-}^{h_+} \bar{\rho}
  (|\hat{w}_1|^2 +|\hat{w}_2|^2+| \hat{w}_3|^2 )\mm{d}y_3\mm{d} \xi  \nonumber
\end{align}
and
\begin{align}
 g \llbracket\bar{\rho} \rrbracket |w|_0^2-\mathcal{I}(w) = &(g \llbracket\bar{\rho} \rrbracket- \vartheta \abs{\xi}^2) \abs{\hat{w}_3}^2_0 - \left\|
 \sqrt{P'(\bar{\rho}) \bar{\rho}}\left|\frac{g}{P'(\bar{\rho})}\hat{w}_3+\mm{i}\xi_1 \hat{w}_1+ \mm{i}\xi_2 \hat{w}_2- \partial_3 \hat{w}_3\right|\right\|^2_0 \nonumber \\
=& \int_{\mathbb{R}^2}\Bigg( ( g \llbracket\bar{\rho} \rrbracket -  \vartheta\abs{\xi}^2 ) \abs{\hat{w}_3}^2 \nonumber \\
& - \int_{h_-}^{h_+} {P'(\bar{\rho})\bar{\rho}}\abs{  \frac{g}{P'(\bar{\rho})}\hat{w}_3 +\mm{i}\xi_1 \hat{w}_1+ \mm{i}\xi_2 \hat{w}_2 - \partial_3 \hat{w}_3 }^2 \mm{d}y_3 \Bigg)\mm{d} \xi. \nonumber
\end{align}
Similarly, by \eqref{parseval}--\eqref{23080314191}, we can compute out that
\begin{align}
\mathcal{U}(w)  = & \int_{\mathbb{R}^2} \int_{h_-}^{h_+}
 \bigg(\varsigma | \partial_3 \hat{w}_3-\mm{i}\xi_1 \hat{w}_1 - \mm{i}\xi_2 \hat{w}_2|^2 \nonumber \\
&  + \mu
  ( |\mm{i}(\xi_1\hat{w}_2- \xi_2  \hat{w}_1)|^2 +|\mm{i}\partial_3 \hat{w}_1+ \xi_1  \hat{w}_3|^2+|  \mm{i}\partial_3 \hat{w}_2+  \xi_2  \hat{w}_3 |^2  \nonumber \\
&+|\partial_3  \hat{w}_3+\mm{i} \xi_1 \hat{w}_1+\mm{i} \xi_2  \hat{w}_2|^2+|\partial_3  \hat{w}_3-\mm{i} \xi_1 \hat{w}_1-\mm{i} \xi_2  \hat{w}_2|^2/3)\bigg) \mm{d}y_3\mm{d}\xi . \nonumber
\end{align}

Let $\varphi(y_3)= -\mm{i} \hat{w}_1(\xi,y_3)$, $\theta(y_3)=-\mm{i} \hat{w}_2(\xi,y_3)$  and $\psi(y_3) = \hat{w}_3(\xi,y_3)$, then we have
\begin{align}
 \begin{cases}
   \|\sqrt{\bar{\rho}}w\|_0^2=  \int_{\mathbb{R}^2} \mathcal{J}(\varphi,\theta,\psi;\xi) \mm{d} \xi ,&\\
 g \llbracket\bar{\rho} \rrbracket |w|_0^2-\mathcal{I}(w)  =  \int_{\mathbb{R}^2}  E(\varphi,\theta,\psi;\xi)\mm{d}\xi,&\\\mathcal{U}(w) =\int_{\mathbb{R}^2}D(\varphi,\theta,\psi;\xi) \mm{d}\xi ,&
 \end{cases}
\label{2022306051007}
\end{align}
where $\mathcal{J}$, $E$ and $D$ are defined by  \eqref{202230614sfd2048},   \eqref{202202307312119} and \eqref{2022306031524}  with  the notations $|\cdot| $ representing  the moduli of  complex functions.
By the expressions of $\mathcal{J}$, $E$ and $D$, we see that
\begin{equation} \begin{cases}
\mathcal{J}(\varphi,\theta,\psi;\xi) = \mathcal{J}(\Re \varphi,\Re \theta,\Re \psi;\xi) +\mathcal{J}(\Im \varphi,\Im \theta,\Im \psi;\xi), &\\ E(\varphi,\theta,\psi;\xi) = E(\Re \varphi,\Re \theta,\Re \psi;\xi) +E(\Im \varphi,\Im \theta,\Im \psi;\xi), &\\
 D(\varphi,\theta,\psi;\xi) = D(\Re \varphi,\Re \theta,\Re \psi;\xi) +D(\Im \varphi,\Im \theta,\Im \psi;\xi).& \end{cases}\label{2023080316061}
\end{equation}

Noting that $\Re\varphi$, $\Re\theta$,  $\Re\psi\in H_0^1(h_-,h_+)$, thus
exploiting  \eqref{2023052safda91431x} and \eqref{2308031344}, we have, for $\abs{\xi} \in (0,\abs{\xi}_{\mm{c}})$,
\begin{align}
 E(\Re\varphi,\Re\theta,\Re\psi;\xi)
 \leqslant &{\lambda^2(\xi)}  \mathcal{J}(\Re\varphi,\Re\theta,\Re\psi;\xi)  + {\lambda(\xi)}  D(\Re\varphi,\Re\theta,\Re\psi;\xi)\nonumber \\
\leqslant & {\Lambda^2} \mathcal{J}(\Re\varphi,\Re\theta,\Re\psi;\xi)  + \Lambda D(\Re\varphi,\Re\theta,\Re\psi;\xi) .
\label{202306051004}
\end{align}
If $\abs{\xi} \geqslant \xi_{\mm{c}}$ for $\vartheta>0$, the expression for $E$ is non-positive, so the above inequality holds automatically. In addition,
for $\xi=0$, we can compute out that
$$ E(\Re\varphi,\Re\theta,\Re\psi;\xi):= -\int_{h_-}^{h_+}
  P'(\bar{\rho})\bar{\rho}|\Re\psi' |^2  \mm{d}y_3. $$
This means \eqref{202306051004} also holds for $\xi=0$.
Therefore we get \eqref{202306051004}  for all $\xi\in \mathbb{R}^2$.

Obviously we also have
\begin{align}
E(\Im\varphi,\Im\theta,\Im\psi;\xi)\leqslant {\Lambda^2}  \mathcal{J}(\Im\varphi,\Im\theta,\Im\psi;\xi)  + {\Lambda}  D(\Im\varphi,\Im\theta,\Im\psi;\xi)\mbox{ for any }\xi\in \mathbb{R}^2.
\label{2023080031608}
\end{align}
By \eqref{2023080316061}--\eqref{202306051004} for any $\xi\in \mathbb{R}^2$ and \eqref{2023080031608}, we arrive at $$E ( \varphi,\theta, \psi;\xi)\leqslant {\Lambda^2}  \mathcal{J}(\varphi,\theta,\psi;\xi)  + {\Lambda}  D(\varphi,\theta,\psi;\xi) . $$
Integrating each side of the above inequality over all $\xi \in\mathbb{R}^{2}$, and using   \eqref{2022306051007}  then proves  \eqref{202230as60508057}.

(2) \emph{Secondly we will derive the energy inequality satisfied by $(\eta^{\mathrm{d}},\eta^{\mathrm{d}})$}.

It is easy to see from the RT  problem and the linear problem \eqref{linearizsded} that $(\eta^{\mathrm{d}}, u^{\mathrm{d}})$ satisfies the following error problem:
\begin{equation}\label{201702052209}\begin{cases}
\eta_t^{\mathrm{d}}=u^{\mathrm{d}} &\mbox{in } \Omega,\\[1mm]
{\bar{\rho}}J^{-1} u_t^\mathrm{d}-\mm{div}_{\mathcal{A} }( {P}'(\bar{\rho})\bar{\rho}\mm{div}\eta^{\mathrm{d}}  \mathbb{I} +\mathbb{S}_{\mathcal{A}}(u^{\mathrm{d}}))=g\bar{\rho}(\mm{div}\eta^{\mathrm{d}} \mathbf{e}^3-\nabla \eta_3^{\mathrm{d}}) \mathbf{e}^3+\mathbf{R}^3 & \mbox{in }\Omega, \\
 \llbracket  u^\mathrm{d}   \rrbracket =0,\
  \llbracket   {P}'(\bar{\rho})\bar{\rho}\mm{div}\eta^{\mathrm{d}}\mathbb{I} +\mathbb{S}_{\mathcal{A}}(u^{\mathrm{d}})\rrbracket  J  \mathcal{A}\mathbf{e}^3 +\vartheta \Delta_{\mm{h}}\eta_3^{\mathrm{d}}\mathbf{e}^3  = {\mathbf{R}^4}  &\mbox{on }\Sigma, \\
(\eta^\mathrm{d}, u^\mathrm{d}) =0&\mbox{on }\partial\Omega\!\!\!\!\!\!-,\\
(\eta^{\mm{d}},  u^{\mm{d}})|_{t=0}=(0,\delta^2 u_\mm{r})  &\mbox{in }  \Omega,
 \end{cases}\end{equation}
 where we have defined that
\begin{align}
\mathbf{R}^3:=&\mathbf{N}^1-{\bar{\rho}}(J^{-1} -1) u_t^\mathrm{a}  + \mm{div}_{\tilde{\mathcal{A}}}  (  P'(\bar{\rho})\bar{\rho} \mm{div}\eta^\mm{a}\mathbb{I} + \mathbb{S}_{\mathcal{A}}(u^\mm{a}) )\nonumber \\
&+ \mm{div}\mathbb{S}_{\tilde{\mathcal{A}}}(u^\mm{a}) -\delta e^{{c_7 t}} \mathbf{R}^1( \tilde{\eta}^0 , \tilde{u}^0 ) /n \label{20203072608561}
\end{align}
and
\begin{align}
\mathbf{R}^4 :=  &\mathbf{N}^2-   \llbracket   {P}'(\bar{\rho})\bar{\rho}\mm{div}\eta^\mathrm{a}\mathbb{I} +\mathbb{S} (u^\mathrm{a})\rrbracket  ( J\mathcal{A}\mathbf{e}^3-\mathbf{e}^3) - \llbracket \mathbb{S}_{\tilde{\mathcal{A}}}(u^\mathrm{a})\rrbracket  J\mathcal{A }\mathbf{e}^3   -\delta e^{{c_7 t}} \mathbf{R}^2( \tilde{\eta}^0 , \tilde{u}^0 )/n.
\label{2020307260856}
\end{align}

 Similarly to \eqref{n0101nnnn2026m}, we
apply $\partial_t $ to \eqref{201702052209}$_2$--\eqref{201702052209}$_4$ to get
\begin{equation}\label{201702052221}\begin{cases}
\bar{\rho}J^{-1}     u_{tt}^{\mathrm{d}}-\mm{div}_{\ml{A} }
( {P}'(\bar{\rho})\bar{\rho}\mm{div}u^{\mathrm{d}} \mathbb{I} +\partial_t\mathbb{S}_{\mathcal{A}}(u^{\mathrm{d}}))   )
 =
g\bar{\rho}(\mm{div}u^{\mathrm{d}} \mathbf{e}^3-
\nabla u_3^{\mathrm{d}} ) +   \mathbf{R}^5&\mbox{ in }  \Omega,\\[1mm]
 \llbracket   {P}'(\bar{\rho})\bar{\rho}\mm{div}u^{\mathrm{d}}\mathbb{I} +\partial_t \mathbb{S}_{\mathcal{A}}(u^{\mathrm{d}}) \rrbracket  J  \mathcal{A}  \mathbf{e}^3+\vartheta \Delta_{\mm{h}}u_3\mathbf{e}^3 =\mathbf{R}^6,\   \llbracket   u^{\mathrm{d}}_t  \rrbracket =0,&\mbox{ on }\Sigma,\\
   u^{\mathrm{d}}_t=0 &\mbox{ on }\partial\Omega\!\!\!\!-,
\end{cases}\end{equation}
where we have define that
  \begin{align}
\mathbf{R}^5:=& \mathbf{R}_t^3 +  \mm{div}_{   \ml{A}_t }
( {P}'(\bar{\rho})\bar{\rho}\mm{div}u^{\mathrm{d}} \mathbb{I}+
\mathbb{S}_{\mathcal{A} } ( u^{\mathrm{d}})
 )-\bar{\rho} J^{-1}_t
 u_t^{\mm{d}}   , \label{2020307260856d} \\
  \mathbf{R}^6:=&\mathbf{R}^4_t-  \llbracket   {P}'(\bar{\rho})\bar{\rho}\mm{div}\eta^{\mathrm{d}}\mathbb{I} +  \mathbb{S}_{\mathcal{A}}(u^{\mathrm{d}}) \rrbracket  \partial_t(J  \mathcal{A}  \mathbf{e}^3). \label{2020307260856d1}
\end{align}

Similarly to \eqref{060817561757m}, we can deduce from \eqref{201702052221} that
\begin{equation} \label{060817561757mi}
\begin{aligned}
\frac{1}{2}\frac{\mm{d}}{\mm{d}t}\left(\|\sqrt{\bar{\rho} } u_t^{\mathrm{d}} \|^2_0+\mathcal{I}  ( u^{\mathrm{d}})-g \llbracket\bar{\rho} \rrbracket |u^{\mathrm{d}}|_0^2
\right)+\mathcal{U}_{\sqrt{J}\mathcal{A} }( u^{\mathrm{d}}_t )= R_1   ,
\end{aligned}
\end{equation}
where
 \begin{align*}
R_1:= &\int (J  \mathbf{R}^5+
g\bar{\rho}(J-1)(\mm{div}u^{\mathrm{d}} \mathbf{e}^3-
\nabla u_3^{\mathrm{d}} ) )\cdot  u_t^{\mathrm{d}}+(
 (1-J){P}'(\bar{\rho}) \bar{\rho} \mm{div}
u^{\mathrm{d}} \mathbb{I}\\
&-J
\mathbb{S}_{\mathcal{A}_t } ( u^{\mathrm{d}}) ): \nabla_{\mathcal{A} }u_t^{\mathrm{d}} -
{P}'(\bar{\rho}) \bar{\rho} \mm{div}
u^{\mathrm{d}} \mm{div}_{\tilde{\ml{A}} } u^{\mathrm{d}}_t)\mm{d}y-\int_\Sigma  \mathbf{R}^6 \cdot  u_t^{\mathrm{d}} \mm{d}y_{\mm{h}}.
\end{align*}
Recalling $u^\mm{d}(0)=\delta^2 u^{\mm{r}}$, we integrate \eqref{060817561757mi} in time from $0$ to $t$ to get
\begin{equation}\label{0314}
\begin{aligned}
\|\sqrt{\bar{\rho} }u_t^\mm{d}\|^2_{0}+ 2\int_0^t \mathcal{U}(u^{\mathrm{d}}_t)\mm{d}\tau =g \llbracket\bar{\rho} \rrbracket |u^{\mathrm{d}}|_0^2-\mathcal{I}(u^{\mathrm{d}}(t))
 +\sum_{i=1}^3\mathcal{R}_i,
\end{aligned}\end{equation}
where
\begin{align*}
&\mathcal{R}_{1}:=2\int_0^t R_1 (\tau)\mm{d}\tau,\ \mathcal{R}_2:= \|\sqrt{\bar{\rho} } u_t^{\mathrm{d}}\|^2_0 \bigg|_{t=0}
+\mathcal{I}(\delta^2 u^{\mm{r}})-g \llbracket\bar{\rho} \rrbracket |\delta^2 u^{\mm{r}}|_0^2,\nonumber \\
&\mbox{ and }
\mathcal{R}_3=  - 2\int_0^t \int  (\mathbb{S}( u^{\mathrm{d}}):
 \nabla_{J {\mathcal{A}}-\mathbb{I}}    u^{\mathrm{d}}+ \mathbb{S}_{J{\mathcal{A}}-\mathbb{I}}( u^{\mathrm{d}}):
 \nabla_{ {\mathcal{A}}}   u^{\mathrm{d}})\mm{d}y \mm{d}\tau.
\end{align*}

(3) \emph{Thirdly we will estimate for the three rest terms $\mathcal{R}_1$--$\mathcal{R}_3$.}

Noting the assumption $  \|\eta\|_3  \leqslant  \delta_1$  in \eqref{201702100940n}, which makes sure that all estimates in Lemmas \ref{lem:201612041032}--\ref{lem:0933} can be satisfied by $(\eta,u)$ here, thus we have
\begin{align}
&1\lesssim \inf_{y\in {\Omega}}\{J,J^{-1}\} \mbox{ (see \eqref{Jdetemrinat})},  \label{Jdsetemrixzcvnat11}
\\
    & \| (   J^{-1 }-1,  \tilde{\mathcal{A}}  ) \|_2\lesssim  \| \eta\|_3 \mbox{ (see \eqref{Jdetemrinatneswn} and \eqref{prtislsafdsfsfds})}, \label{Jdsetemrixzcvnat} \\
& \| \mathbf{N}^1\|_0\lesssim  \| \eta \|_{3}  \|\eta\|_{2 } \mbox{ (referring to \eqref{06011711jumpv})},\label{Jdsetemrixzcvnat1} \\
 &
\begin{cases}     \|( J-1,J\mathcal{A}\mathbf{e}^3-\mathbf{e}^3 ) \|_2\lesssim  \| \eta\|_3 \mbox{ (see   \eqref{Jdetemrinatneswn}, \eqref{201702022056} and \eqref{201612071026})},\\
   \|\partial_t(J^{-1},\mathcal{A},J\mathcal{A}\mathbf{e}^3)\|_1\lesssim  \| u\|_{2} \mbox{ (see   \eqref{Jtestismtsnn}, \eqref{prtislsafdsfs} and \eqref{201612071401xx})},\\
\|\mathbf{N}_t^1\|_0 \lesssim \|\eta\|_3\|u\|_2,\   |  \mathbf{N}^2_t |_{1/2}
    \lesssim   \|\eta\|_3\|u\|_3.  \mbox{ (see   \eqref{20223080800630} and \eqref{20210880633})}.
 \end{cases}\label{Jdvnat}
\end{align}
Thanks to  \eqref{appdesimtsofu1857}, \eqref{201702100940} and  \eqref{Jdsetemrixzcvnat},  it is easy to see that
\begin{align}
\mathcal{R}_3\lesssim\int_0^t \delta^3e^{3c_7 \tau}\mm{d}\tau\lesssim \delta^3e^{3c_7 t} .
\label{202308081700}
\end{align}

Applying $\|\cdot \|$ to \eqref{201702052209}$_2$ yields
$$\begin{aligned}
 \|\bar{\rho}J^{-1} u_t^{\mathrm{d}}\|_0 =\|\mm{div}_{\mathcal{A}} ({P}'(\bar{\rho})\bar{\rho}\mm{div}\eta^{\mathrm{d}} \mathbb{I} +\mathbb{S}_{\mathcal{A}}(u^{\mathrm{d}}) )+g\bar{\rho}(\mm{div}\eta^{\mathrm{d}} \mathbf{e}^3-\nabla \eta_3^{\mathrm{d}}
 )+\mathbf{R}^3)\|_0.
 \end{aligned}  $$
Thus, using   \eqref{Jdsetemrixzcvnat11} and  \eqref{Jdsetemrixzcvnat}, we obtain
\begin{equation}\label{appesimtsofu12}       \| u_t^{\mathrm{d}}\|^2_0
  \lesssim\|( \eta^{\mm{d}}, u^{\mm{d}})\|_2^2+  \|\mathbf{R}^3
  \|_0^2. \end{equation}
In addition, making use of \eqref{appdesimtsofu1857},  \eqref{20222302522222}, \eqref{201702100940}, \eqref{Jdsetemrixzcvnat} and \eqref{Jdsetemrixzcvnat1}, we can estimate that
\begin{align*}
\| \mathbf{R}^3\|_0 =&\|\mathbf{N}^1-{\bar{\rho}}(J^{-1}-1 ) u_t^\mathrm{a} +\mm{div}_{\tilde{\mathcal{A}}}  ( P'(\bar{\rho})\bar{\rho} \mm{div}\eta^\mm{a}\mathbb{I} + \mathbb{S}_{\mathcal{A}}(u^\mm{a}) )\nonumber \\
& +\mm{div}\mathbb{S}_{\tilde{\mathcal{A}}}(u^{\mm{a}})-\delta e^{{c_7 t}} \mathbf{R}^1( \tilde{\eta}^0 , \tilde{u}^0 ) /n\|_0^2 \lesssim\delta^2 e^{2\Lambda t}+\delta n^{-1/2}e^{c_7 t} .
\end{align*}
Thus, putting the above   estimate    into \eqref{appesimtsofu12}  and taking then $t=0$,
we can have
\begin{equation}\label{estimeR1x}
 \|  u_t^{\mathrm{d}}\|^2_0 \bigg|_{t=0}\lesssim  \delta^4 +\delta^2 n^{-1}+ \|  \delta^2 u^{\mm{r}}\|_2^2
\lesssim   \delta^4 , \nonumber
 \end{equation}
where we have used \eqref{2022308081221} and \eqref{201702091755}. Similarly,
\begin{equation}\nonumber
\mathcal{I}(\delta^2 u^{\mm{r}})-g \llbracket\bar{\rho} \rrbracket | \delta^2 u^{\mm{r}}|_0^2
\lesssim  \delta^4.
 \end{equation}
Putting the above two estimates together yields
\begin{align}\label{estimeR1}
\mathcal{R}_2\lesssim  \delta^4\lesssim  \delta^3e^{3c_7 t}.
\end{align}

Now we turn to the estimate of $\mathcal{R}_1$. Recalling the definitions of $ \mathbf{N}^1$, $\mathbf{N}^2$ in \eqref{2022307260853}, $\mathbf{R}^3$ in  \eqref{20203072608561}, $\mathbf{R}^4$ in \eqref{2020307260856}, $\mathbf{R}^5$ in \eqref{2020307260856d} and $\mathbf{R}^6$ in \eqref{2020307260856d1}, it is easy to see that
 \begin{align*}
& R_1(t) \lesssim  (1+\|J-1\|_{2}) \| \mathbf{R}^5 \|_0\|  u_t^{\mathrm{d}}\|_0+
((1+\|\tilde{\mathcal{A}}\|_2)(\| J -1\|_{2}\| u^{\mm{d}}\|_1\\
&\qquad \quad +(1+\|J-1\|_2)\|{\mathcal{A} }_t\|_1\|u^{\mathrm{d}}\|_2)
+ \|\tilde{\ml{A}}\|_{2}\|u^{\mathrm{d}}\|_1+| \mathbf{R}^6|_{1/2})\|u_t^{\mathrm{d}} \|_1,\\
&\|\mathbf{R}^5\|_0\lesssim \| \mathbf{R}_t^3 \|_0+
 \|\ml{A}_t\|_1 (1+\|\tilde{\mathcal{A} }\|_{2})\| u^{\mathrm{d}}\|_3
 +\|J^{-1}_t\|_{1}\|
 u_t^{\mm{d}}  \|_1,\\
 &| \mathbf{R}^6|_{1/2}\lesssim |\mathbf{R}^4_t|_{1/2}+
        (\|\eta^\mm{d}\|_3 +(1+ \|\tilde{\mathcal{A}}\|_{2} )\|u^{\mathrm{d}}\|_3 )\|\partial_t(J\mathcal{A}\mathbf{e}^3)\|_1,\\
&\|\mathbf{R}^3_t\|_0\lesssim\|\mathbf{N}^1_t\|_0+\|J^{-1}-1\|_2\|u_{tt}^\mathrm{a} \|_0 +\|J^{-1}_t\|_1\|u_t^\mathrm{a} \|_1
+ \|{\tilde{\mathcal{A}}}_t\|_1  (1+\|\tilde{\mathcal{A}}\|_2)\|(\eta^\mm{a},u^\mm{a})\|_3  \nonumber \\
 &\qquad \qquad + \|{\tilde{\mathcal{A}}}\|_2
(  \| u^\mm{a}\|_2 +(1+\|\tilde{\mathcal{A}}\|_2)\|u^\mm{a}_t\|_2)+\delta e^{{c_7 t}} \|\mathbf{R}^1( \tilde{\eta}^0 , \tilde{u}^0 )\|_0/n  , \\
&\|\mathbf{R}^4_t\|_{1/2} \lesssim |\mathbf{N}^2_t|_{1/2}+   \| J\mathcal{A}\mathbf{e}^3-\mathbf{e}^3\|_2  \|(u^\mathrm{a},u^\mathrm{a}_t)\|_2
\nonumber \\
&\qquad \qquad+ (1+\| J\mathcal{A}\mathbf{e}^3-\mathbf{e}^3\|_2)(  \|\tilde{\mathcal{A}}\|_2\|u^\mathrm{a}_t\|_2 +\|\tilde{\mathcal{A}}_t\|_1\|u^\mathrm{a}\|_3)
+(\|\eta^\mathrm{a} \|_3
\nonumber \\
& \qquad \qquad+(1+\|\tilde{\mathcal{A}}\|_2)\|u^\mathrm{a}\|_3) \|\partial_t(J\mathcal{A}\mathbf{e}^3)\|_1    +\delta e^{{c_7 t}} \|\mathbf{R}^2( \tilde{\eta}^0 , \tilde{u}^0 )\|_1/n .
\end{align*}
Thus, making use of \eqref{appdesimtsofu1857}, \eqref{201702100940},  \eqref{Jdsetemrixzcvnat}, \eqref{Jdvnat} and the condition $\|\eta\|_3\leqslant 1$, the above seven estimates reduce to
 \begin{align*}
& R_1(t) \lesssim   \delta^3 e^{3c_7 t} +\delta e^{ c_7 t}(\| \mathbf{R}^5\|_0+  | \mathbf{R}^6|_{1/2})+
(\delta^2 e^{2 c_7 t}+ | \mathbf{R}^6|_{1/2})\|u_t\|_{1}  ,\\
&\|\mathbf{R}^5\|_0\lesssim  \delta^2 e^{2 c_7 t}+\delta e^{c_7 t}(\|u \|_3   +\| u_t\|_1) + \| \mathbf{R}_t^3 \|_0  ,\\
  & | \mathbf{R}^6|_{1/2}\lesssim  \delta^2  e^{2c_7 t}+ \delta  e^{ c_7 t}\|u \|_3+ |\mathbf{R}^4_t|_{1/2}
, \\
&  \|\mathbf{R}^3_t\|_0\lesssim\delta^2 e^{2c_7 t} + \delta e^{{c_7 t}}n^{-1}\|\mathbf{R}^2( \tilde{\eta}^0 , \tilde{u}^0 )\|_1 +\|\mathbf{N}^1_t\|_0, \\
&|\mathbf{R}^4_t|_{1/2} \lesssim \delta^2  e^{2c_7 t} +\delta e^{{c_7 t}}n^{-1}|\mathbf{R}^2( \tilde{\eta}^0 , \tilde{u}^0 )|_{1/2}+  |\mathbf{N}^2_t|_{1/2}\nonumber \\
& \|  \mathbf{N}^1_t\|_0 \lesssim \delta^2 e^{2c_7 t}\mbox{ and } |{\mathbf{N}^2_t} |_{1/2}\lesssim \delta e^{c_7 t}\|u\|_3 ,
\end{align*}
 which, together with \eqref{2022308081221} and \eqref{20222302522222}, imply
 \begin{align*}
& R_1(t) \lesssim \delta^3 e^{3c_7 t} + \delta^2 e^{2c_7 t}(\|u\|_3+\|u_t\|_{1})+
\delta e^{c_7 t}\|u \|_{3}\|u_t\|_{1} .
\end{align*}
 Integrating the above inequality over $(0,t)$, and then using \eqref{201702100940} and H\"older's inequality, we
 further have \begin{equation}\label{201602071414MH}
\mathcal{R}_1\lesssim \delta^3 e^{3c_7 t}. \end{equation}

Now summing up   the three estimates \eqref{202308081700}, \eqref{estimeR1}  and \eqref{201602071414MH} yields
\begin{align*}
 \sum_{i=1}^3\mathcal{R}_1 \lesssim \delta^3 e^{3c_7 t}.
\end{align*}
Putting the above estimate into \eqref{0314}, one obtains
\begin{align}
 \| \sqrt{\bar{\rho} }u_t^\mm{d}\|^2_{0}+2 \int_0^t \mathcal{U}_{\mathcal{A}} (u_t^\mm{d})\mm{d}\tau  \leqslant g \llbracket\bar{\rho} \rrbracket |u^{\mathrm{d}}|_0^2-\mathcal{I}(u^{\mathrm{d}})+ c\delta^3 e^{3c_7 t}.
 \label{202308111657}
 \end{align}

(4) \emph{Finally we are in the position to the derivation of \eqref{ereroe} by applying the largest growth rate $\Lambda$ in \eqref{202230as60508057}.}

Thanks to \eqref{202230as60508057}, we have
\begin{align}
&g \llbracket\bar{\rho} \rrbracket |u^{\mathrm{d}}|_0^2-\mathcal{I}(u^{\mathrm{d}})  \leqslant {\Lambda^2} \|\sqrt{ \bar{\rho}}u^{\mathrm{d}}\|_0^2
+\Lambda  \mathcal{U}(u^{\mathrm{d}}) . \nonumber
\end{align}
Combining the above  inequality and \eqref{202308111657} together implies
\begin{equation}
\label{new0311}
\begin{aligned}  \|\sqrt{\bar{\rho}} u_t^\mm{d}\|^2_{0}+2\int_0^t \mathcal{U}_{\mathcal{A}}(u^{\mathrm{d}}_t)\mm{d}\tau \leqslant  {\Lambda^2} \|\sqrt{\bar{\rho}} u^\mm{d}\|_{0}^2 + {\Lambda}\mathcal{U}(u^{\mathrm{d}})+  c\delta^3 e^{3c_7 t}.
\end{aligned}\end{equation}

 We apply Newton--Leibniz's formula, Cauchy--Schwarz's inequality and the fact $u^{\mm{d}}(0)=\delta^2 u^{\mm{r}}$ to find that
\begin{align}
 \Lambda\mathcal{U}(u^{\mathrm{d}})
 = &2\Lambda \int_0^t ((  \varsigma-2\mu/3)\mm{div}u^{\mathrm{d}}\mm{div}u^{\mathrm{d}}_{\tau}+\mu \mathbb{D} u : \partial_\tau \mathbb{D} u  \mathrm{d}\tau  \nonumber \\
 &  +\Lambda\mathcal{U}(\delta^2 u^\mm{r}) \leqslant \Lambda^2\int_0^t\mathcal{U}(u^{\mathrm{d}})\mm{d}\tau+\int_0^t\mathcal{U}(u^{\mathrm{d}}_\tau )\mathrm{d}\tau+c\delta^3 e^{3c_7 t}. \nonumber
 \end{align}
 Combining \eqref{new0311} with the above estimate, one gets
 \begin{equation}\label{inequalemee}
 \frac{1}{\Lambda}\|\sqrt{\bar{\rho}}  u_t^{\mathrm{d}}\|^2_{0 }+
\mathcal{U}(u^{\mathrm{d}})  \leqslant    {\Lambda}\|\sqrt{\bar{\rho}} u^{\mathrm{d}}\|^2_{0 }+2 {\Lambda}\int_0^t
\mathcal{U}(u^{\mathrm{d}})\mm{d}\tau + c\delta^3 e^{3c_7 t}.
 \end{equation}
In addition,
\begin{equation*}
\frac{\mm{d}}{\mm{d}t}\|\sqrt{\bar{\rho}} u^\mm{d} \|^2_{0}= 2\int
\bar{\rho} u^\mm{d} \cdot  u^\mm{d}_t \mm{d}y
\leqslant  \frac{1}{\Lambda}\|\sqrt{ \bar{\rho} }  u_t^\mm{d} \|^2_{0}
+\Lambda\|\sqrt{ \bar{\rho}} u^\mm{d} \|^2_{0} . \end{equation*}
If we put the previous two estimates together, we get the differential inequality
\begin{equation}\label{growallsinequa}
 \frac{\mm{d}}{\mm{d}t} \|\sqrt{\bar{\rho}} u^\mm{d}\|^2_{0}+\mathcal{U}(u^{\mathrm{d}}) \leqslant  2\Lambda\left( \|\sqrt{\bar{\rho}} u^\mm{d}(t)\|^2_{0}
 +\int_0^t\mathcal{U}(u^{\mathrm{d}})\mathrm{d}\tau\right) +c\delta^3 e^{3c_7 t}.
\end{equation}
Recalling $u^{\mm{d}}(0)=\delta^2u^{\mm{r}}$ and $c_7  =\lambda({\xi^1})\in (2\Lambda/3,\Lambda)$ in Proposition \ref{growingmodesolneriodic}, one can apply Gronwall's inequality to \eqref{growallsinequa} to conclude that
\begin{align}
\|\sqrt{  \bar{\rho}} u^{\mathrm{d}}\|^2_{0}+  \int_0^t\mathcal{U}(u^{\mathrm{d}})\mm{d}\tau \lesssim  & e^{2\Lambda t}\left(\int_0^t   \delta^3 e^{(3c_7 -2\Lambda)\tau}\mm{d}\tau
+\|\sqrt{ \bar{\rho}}\delta^2 u^{\mathrm{r}}\|^2_{0}\right)
\lesssim   \delta^3e^{3c_7 t} ,\label{estimerrvelcoity}
 \end{align}Moreover,  we can further deduce from \eqref{new0311}, \eqref{inequalemee},  \eqref{estimerrvelcoity} and Korn's inequality that
\begin{equation}\label{uestimate1n}
\|u^{\mathrm{d}}\|_{1 }^2+\| u_t^{\mathrm{d}}\|^2_0 +\|u^{\mathrm{d}}_\tau\|^2_{L^2((0,t),H^1)}\lesssim\delta^3 e^{3c_7 t}.
\end{equation}

Finally it follows from  \eqref{201702052209}$_1$, \eqref{201702052209}$_5$ and \eqref{uestimate1n}  that
\begin{equation}\begin{aligned}\label{erroresimts}
 \|\eta^{\mathrm{d}}\|_{1 }\lesssim  &  \int_0^t \|u^{\mathrm{d}}\|_{1} \mm{d}\tau
\lesssim  \sqrt{ \delta^3e^{3\Lambda t}} .
\end{aligned}\end{equation}
Summing up the two estimates \eqref{uestimate1n} and \eqref{erroresimts}, and then using trace estimate, we
obtain the desired estimate \eqref{ereroe}.  \hfill$\Box$
\end{pf}

\section{Existence of escape times}\label{sec:030845}

Now we are in a position to show Theorem \ref{thm:0202}.
Let $\delta$ satisfy \eqref{2022309091827}, $(\eta,u)$ be the strong  solution constructed in Section \ref{20230721} with an existence time $[0,T_{\mm{loc}})$, and $(\eta^{\mm{d}},u^{\mm{d}})$ be defined  in Section \ref{20230721}.
Let $\epsilon _0\in (0,1]$ be a constant, which will be defined in \eqref{defined}.  We further define
 \begin{align}\label{times}
& T^\delta:=c_7^{-1}\mm{ln}({\epsilon _0}/{\delta})>0, \mbox{ i.e., }
 \delta e^{c_7 T^\delta }=\epsilon _0,\\
&T^*:=\sup\left\{t\in(0,{T}_{\mm{loc}})\left|~{\sqrt{\|\eta (\tau)\|_3^2+
\| u (\tau)\|_2^2 + \vartheta |d(\tau)|_3^2}}\leqslant  2\tilde{c}_4c_3\mbox{ for any }\tau\in [0,t)\right.\right\},\nonumber\\
&  T^{**}:=\sup\left\{t\in (0,T_{\mm{loc}})\left|~\left\|(\eta,u)(\tau)\right\|_0\leqslant  2 \tilde{c}_4\delta e^{c_7 t}\mbox{ for any }\tau\in [0,t)\right\}\right..\nonumber
 \end{align}
 Noting that, by \eqref{2022307282134} and \eqref{2022309091827},
\begin{align} \left.\sqrt{\|\eta(t)\|_3^2+\|u(t)\|_2^2+\vartheta  |d(t)|_3^2 }\right|_{t=0}
&=\sqrt{\|\eta_0^\delta\|_3^2 +\|u_0^\delta\|_2^2+ \vartheta |d^\delta_0|_3^2 } \leqslant \tilde{c}_4\delta<2\tilde{c}_4c_3\leqslant  \delta_2,
\label{201809121553}\end{align}
 thus $T^*>0$ by  \eqref{2022307262218}.  Similarly, we also have  $T^{**}>0$.
Moreover, we can easily see that
 \begin{align}\label{0502n1}
&\sqrt{\|\eta (T^*)\|_3^2+ \| u (T^*)\|_2^2+ \vartheta |d(T^*)|_3^2 }= 2\tilde{c}_4c_3\leqslant \delta_2\mbox{ if }T^*<\infty ,\\
\label{0502n111}  & \left\|(\eta,u) (T^{**}) \right\|_0
=2 \tilde{c}_4\delta e^{c_7 T^{**}}\quad\mbox{ if }T^{**}<T_{\mm{loc}}.
\end{align}

We denote ${T}^{\min}:= \min\{T^\delta ,T^*,T^{**}\}$, thus ${T}^{\min}<T_{\mm{loc}}$ by \eqref{0502n1} and Proposition \ref{pro:0401nd}.
 Noting that, by \eqref{2022307282134},
$$
  \sup_{0\leqslant  t\leqslant  T^{\min}}\sqrt{\|\eta(t)\|_3^2 +\|u(t)\|_2^2+\vartheta |d(t)|_3^2 }  \leqslant \tilde{c}_4\delta\leqslant \tilde{c}_4c_3\leqslant \delta_1,$$
thus, in view of both the second and third assertions in Proposition \ref{pro:0401nd},  we deduce from the estimate \eqref{2016121521430850} satisfied by $(\eta,u,d,\vartheta\eta^2)$ that for all $t\in
[0, {T}^{\min}]$,
\begin{align}
& \tilde{\mathcal{E}}(t)+ \vartheta( \|  \eta^2(t)  \|_{2,1}^2 +c^{-1}  |d(t)|_3^2)
+ c^{-1}\int_0^t\mathcal{D} (\tau)\mm{d}\tau\nonumber \\
&
\leqslant  c\delta^2e^{2c_7 t} +c_7\int_0^t(\tilde{\mathcal{E}}(\tau)+\vartheta\|  \eta^2(\tau)  \|_{2,1}^2 )\mm{d}\tau.
 \label{0503} \end{align}
   Applying Gronwall's inequality to the above estimate, we deduce that
$$
 \tilde{\mathcal{E}}(t)+\vartheta (|d(t)|_3^2+ \|  \eta^2(t)  \|_{2,1}^2 )\lesssim \delta^2\left(e^{2 c_7 t}+
 \int_0^t e^{c_7 (t+\tau)} \mm{d}\tau\right)\lesssim\delta^2 e^{2c_7 t}.
 $$
 Putting the above estimate to  \eqref{0503}, we get
$$ \tilde{\mathcal{E}}(t)+\vartheta( |d(t)|_3^2+ \|  \eta^2(t)  \|_{2,1}^2 )+\int_0^t\mathcal{D} (\tau)\mm{d}\tau \lesssim \delta^2 e^{2 c_7 t},$$
which, together with \eqref{2018008121027} satisfied by $(\eta,u)$, yields that
\begin{align}
&\sqrt{\|\eta\|_3^2+\|u\|_2^2+ \vartheta |d|_3^2 +\|u_t\|_0^2+\|u\|_{L^2((0,t),H^3)}^2+\|u_\tau\|_{L^2((0,t),H^1)}^2}\nonumber \\
& \leqslant  \tilde{c}_5 \delta e^{c_7 t}\leqslant  \tilde{c}_5\epsilon _0\mbox{ on }[0, {T}^{\min}].
\label{201702092114}
\end{align}

Let $\beta=\tilde{c}_5$ in \eqref{201702100940} and then we denote the constant $c$ in \eqref{ereroe} in Proposition \ref{lem:0401} by $\tilde{c}_6$.
Now we define that
 \begin{equation}\label{defined}
\epsilon_0:=\min\left\{\frac{ {c}_3}{\tilde{c}_5},
\frac{\tilde{c}_4^2}{4\tilde{c}_6^2},
\frac{c_8^2}{\tilde{c}_6^2},1 \right\}>0,
 \end{equation}
Noting that $(\eta,u)$ satisfies \eqref{201702092114}, where $\delta ^{c_7t}\leqslant 1$  and $\epsilon _0\leqslant c_3/ {\tilde{c}_5} \leqslant \delta_1/{\tilde{c}_5}$ (i.e., $\tilde{c}_5\epsilon_0\leqslant \delta_1\leqslant 1$) by the definitions of $c_3$ and $\epsilon_0$, thus, by Proposition \ref{lem:0401} with $\beta=\tilde{c}_5$ and $c=\tilde{c}_6$, we immediately see that
\begin{equation}\label{ereroelast}
  \|(\eta^{\mm{d}} ,   u^{\mm{d}}) \|_1+ |(\eta^{\mm{d}}, u^{\mm{d}})|_0  \leqslant  \tilde{c}_6\sqrt{\delta^3e^{3c_7 t}}\mbox{ on }[0, {T}^{\min}],  \end{equation}where
$(\eta^{\mathrm{d}}, u^{\mathrm{d}})=(\eta, u)-(\eta^\mm{a},u^{\mm{a}})$.
Consequently, we further have the  relation
\begin{equation}
\label{201702092227}
T^\delta =T^{\min},
\end{equation}
 which can be showed by contradiction as follows:

If $T^{\min}=T^*$, then $T^*<\infty$. Recalling $\epsilon_0\leqslant  {c}_3/\tilde{c}_5$ and $\tilde{c}_4 \geqslant 1$, we can deduce from \eqref{201702092114} that
 \begin{equation*}\begin{aligned}
 \sqrt{{\|\eta (T^*)\|_3^2+
\| u (T^*)\|_2^2}+|d (T^*)|_3^2}\leqslant  \tilde{c}_5\epsilon_0\leqslant \tilde{c}_4{c}_3 < 2 \tilde{c}_4 c_3,
 \end{aligned} \end{equation*}
which contradicts \eqref{0502n1}. Hence, $T^{\min}\neq T^*$.

If $T^{\min}=T^{**}$, then $T^{**}<T^*\leqslant T_{\mm{loc}}$.
 Noting that $\sqrt{\epsilon_0}\leqslant \tilde{c}_4/2\tilde{c}_6$ and $c_5 \leqslant \tilde{c}_4$, we obtain from  \eqref{appdesimtsofu1857}, \eqref{ereroe}, \eqref{times} and the fact
 $\varepsilon_0\leqslant \tilde{c}_4^2/4 \tilde{c}_5^2$ that \begin{equation*}\begin{aligned}
 \| (\eta,u)  (T^{**})\|_0
&\leqslant  \|(\eta^\mm{a}, u^\mm{a})(T^{**})\|_0+\|(\eta^\mm{d}, u^\mm{d}) (T^{**})\|_0\leqslant   \delta e^{{c_7 T^{**}}}(
\tilde{c}_4+ \tilde{c}_6\sqrt{\delta e^{c_7 T^{**}}})\\
&\leqslant   \delta e^{{c_7 T^{**}}}(\tilde{c}_4+ \tilde{c}_6\sqrt{\epsilon _0})
\leqslant  3c_4 \delta e^{\tilde{c}_4 T^{**}}/2<2\tilde{c}_4\delta e^{c_7 T^{**}},
 \end{aligned} \end{equation*}
which contradicts \eqref{0502n111}. Hence, $T^{\min}\neq T^{**}$. We immediately see that \eqref{201702092227} holds.

Noting that $\sqrt{\epsilon_0}\leqslant c_8/\tilde{c}_6$,
making use of \eqref{2023061632013}, \eqref{0501}, \eqref{times} and  \eqref{ereroelast}, we can  deduce that
 \begin{align*}
 \|\omega_3(T^\delta)\|_{0}\geqslant  & \|\omega^{\mathrm{a}}_3(T^\delta )\|_{0}-\|\omega^{\mm{d}}_3(T^\delta )\|_{0 }
  \geqslant  \delta e^{c_7 T^\delta }( \|\chi_{n,n}{\omega}_3/n\|_{0}- \tilde{c}_6\sqrt{\delta e^{c_7 T^\delta }}) \\
 = & (2c_8 -\tilde{c}_6\sqrt{\epsilon _0})\epsilon _0 \geqslant  c_8\epsilon _0 ,
 \end{align*}
 where   $\omega=\tilde{\eta}^{0}$ or $\tilde{u}^{0}$.
 Similarly, we also can verify that $
\|\omega_{\mm{h}}(T^\delta)\|_{0}$,  $|{\omega}_3(T^\delta)|_{0}\geqslant  c_8\epsilon _0  $.
This completes the proof of  Theorem \ref{thm:0202} by taking $\epsilon = c_8\epsilon _0$ and $c_k= L_k$, where $L_k$ is provided by \eqref{20223072522255} for $k=1$ and $2$.
\appendix
\section{Analysis tools}\label{sec:09}
\renewcommand\thesection{A}
This Appendix is devoted to listing some mathematical analysis tools, which have been used in the previous sections.
It should be remarked that in this appendix we still adapt the simplified mathematical notations in Section \ref{202306071229}.
For simplicity, we still use the notation $a\lesssim b$ to mean that $a\leqslant cb$ for some constant $c>0$,
where the positive constant $c$ may depend on the domain and other given parameters in the lemmas below.
\begin{lem}
\label{201806171834}
Embedding inequality (see \cite[Theorem 7.58]{RAADSBS} and \cite[Theorems 4.12]{ARAJJFF}): Let $D\subset \mathbb{R}^3$ be a domain satisfying the cone condition, then
\begin{align}
&\label{esmmdforinftfdsdy}\|f\|_{L^p(D)}\lesssim \| f\|_{H^1(D)}\;\mbox{ for }2\leqslant p\leqslant 6,\\
&\label{esmmdforinfty}\|f\|_{C^0(\overline{D})}= \|f\|_{L^\infty(D)}\lesssim\| f\|_{H^2(D)},\\
& \|\phi\|_{L^{4}(\mathbb{R}^2)}\lesssim |\phi|_{1/2}  \label{201811231df301},\\& \|\phi\|_{L^{\infty}(\mathbb{R}^2)}\lesssim\|\phi\|_{W^{1,4}(\mathbb{R}^2)} \lesssim |\phi|_{2}  \label{201811sdf231301} ,
\end{align}
where $\overline{D}$ denotes the closure of $D$.
\end{lem}
\begin{lem}
Interpolation inequality in $H^j$ (see \cite[5.2 Theorem]{ARAJJFF}): Let $D$ be a domain in $\mathbb{R}^n$ satisfying the cone condition, then
 for any $0\leqslant j< i$, $\varepsilon>0$,
\begin{equation}
\|f\|_{H^j(D)}\lesssim\|f\|_{L^2(D)}^{1-\frac{j}{i}}\|f\|_{H^i(D)}^{\frac{j}{i}}\leqslant c(\varepsilon,j)\|f\|_{L^2(D)} +\varepsilon\|f\|_{H^i(D)},  \label{201807291850}
\end{equation}
where the constant $c(\varepsilon,j)$ depends on the domain, $j$ and $\varepsilon$.
\end{lem}

\begin{lem}\label{xfsddfsf20180508}
Product estimates: Let $D\subset \mathbb{R}^3$ be a   domain  satisfying the cone condition. The functions $f$, $g$ are
defined on $D$, and $\phi$, $\varphi$ are defined on $\mathbb{R}^2$.
 \begin{enumerate}[(1)]
   \item Product estimates in $H^i$: \begin{align}
\label{fgestims}
  \|fg\|_{H^i(D)}\lesssim \begin{cases}
                      \|f\|_{H^1(\Omega)}\|g\|_{H^1(D)} & \hbox{for }i=0; \\
  \|f\|_{H^i(D)}\|g\|_{H^2(D)} & \hbox{for }0\leqslant i\leqslant 2;  \\
                    \|f\|_{H^2(D)}\|g\|_{H^i(D)}+\|f\|_{H^i(D)}\|g\|_{H^2(D)}& \hbox{for }i=3.
                    \end{cases}
 \end{align}
   \item Product estimates in $H^{s}(\mathbb{R}^2)$: \begin{align}
&\label{06041605} |\phi\varphi|_{1/2}\lesssim
 \|\phi\|_{W^{1,4}(\mathbb{R}^2)}
|\varphi|_{1/2}  \lesssim
 |\phi|_{3/2}
|\varphi|_{1/2} \end{align}
 \end{enumerate}
\end{lem}
 \begin{rem}In particular, by   \eqref{06041605}, we further have\begin{align} |\phi\varphi|_{j+1/2}\lesssim   |\phi|_{j+1/2}|\varphi|_{3/2}+|\varphi|_{j+1/2}|\phi|_{3/2}\mbox{ for }j=1\mbox{ and }2.\label{202307231325} \end{align}\end{rem}
 \begin{pf}
The product estimate \eqref{fgestims} can be shown by using H\"older's inequality and the embedding inequalities \eqref{esmmdforinftfdsdy}--\eqref{esmmdforinfty}. The estimate
\eqref{06041605} can be obtained by following the proof of \cite[Lemma A.2]{guo2013decay}  and using  the estimates \eqref{201811231df301}--\eqref{201811sdf231301}.
\hfill $\Box$
\end{pf}
\begin{lem}\label{xfsddfsf201805072234Ponc}Poinc\'are's inequality (see \cite[Lemma A.10]{guo2013decay}):  It holds that
\begin{equation}
\|w\|_{0} \lesssim  \|\partial_3 w\|_{0}\mbox{ for any }w\in H^1_0.
\label{2022306232041}
\end{equation}
 \end{lem}
\begin{lem}\label{xfsddfsf201805072234}
 Korn's inequality:
\begin{equation}
\|w\|_1\lesssim \|\mathbb{D}w-2\mm{div}w \mathbb{I}/3 \|_{0}\mbox{ for any }w\in H_0^1.
\label{2022306241238}
\end{equation}
\end{lem}
\begin{pf}  It's well-known that, see \cite[Theorem 10.16]{FENA},
\begin{align}
\label{202230623230241}
\|\nabla w\|_0\lesssim \|\mathbb{D}w-2\mm{div}w \mathbb{I}/3\|_0\mbox{ for any }w\in W^{1,2}(\mathbb{R}^3).
\end{align}
Now we defined that
$$\tilde{w}:=\begin{cases}
0&\mbox{for }y\in \mathbb{R}^3\backslash\Omega; \\
w&\mbox{for }y\in \Omega .
\end{cases}$$
By \cite[Theorems 5.29]{ARAJJFF},  $\tilde{w}\in W^{1,2}(\mathbb{R}^3)$, and thus we can use \eqref{202230623230241} to derive that
$$\|\nabla w\|_0\lesssim \|\mathbb{D}w-2\mm{div}w \mathbb{I}/3\|_0 . $$
Since $w|_{\partial \Omega\!\!\!\!-}=0$, we further derive from Poinc\'are's inequality \eqref{2022306232041} and the above estimate.
\end{pf}
\begin{lem}
\label{xfsddfsf201805072254}
 Trace estimate: It holds that, for given $i\geqslant 0$,
\begin{align}
\|f|_{y_3=a}\|_{H^{i+1/2}(\mathbb{R}^2)}\lesssim \|f\|_{\underline{1+i},0}^{1/2}\| f\|_{\underline{i},1}^{1/2 } \lesssim \|f\|_{\underline{1+i},0}^{1/2}\|\partial_3f\|_{\underline{i},0}^{1/2}
+\|f\|_{\underline{1+i},0}    \label{2022306232141}
\end{align}
 for any $f\in H^{1+i}$ and for any $a\in (h_-,h_+)$.
\end{lem}
\begin{rem}
By the trace estimate, we further have
\begin{equation}
\label{201808051835}
|f|_{i+1/2}\lesssim \|f\|_{\underline{i},1}\mbox{ for any }f\in H^{i+1}.
\end{equation}
\end{rem}
\begin{pf}
It suffices to consider the case $i=0$. Let $\mathbb{R}^3_{h_-,+}=\mathbb{R}^2\times (h_-,+\infty)$ and $\mathbb{R}^3_{h_+,-}=\mathbb{R}^2\times (-\infty,h_+)$. It is well-known that,   for any $f\in H^1(\mathbb{R}^3_{h_-,+})$, $\varphi\in H^1(\mathbb{R}^3_{h_+,-})$,  $a\in [h_-,+\infty)$ and $b\in (-\infty,h_+]$,
 \begin{align}
\|f|_{y_3=a}\|_{H^{1/2}(\mathbb{R}^2)}\lesssim\sum_{|\alpha|\leqslant 1} \|\partial_{\mathrm{h}}^\alpha f\|_{L^2(\mathbb{R}^3_{h_-,+})}^{1/2}\| f\|_{H^1(\mathbb{R}^3_{h_-,+})}^{1/2}
  \label{202230623aa2141}
\end{align}
and \begin{align}
\|\varphi|_{y_3=b}\|_{H^{1/2}(\mathbb{R}^2)}\lesssim\sum_{|\alpha|\leqslant 1} \|\partial_{\mathrm{h}}^\alpha \varphi\|_{L^2(\mathbb{R}^3_{h_+,-})}^{1/2}\| \varphi\|_{H^1(\mathbb{R}^3_{h_+,-})}^{1/2},
  \label{202230623aa2141x}
\end{align} see Theorem 5.7 in \cite[Section 1.5.3]{CSXI2020}.
Thus using cut-off functions depending on $y_3$ to cut off the extension function of $f$, we  easily obtain \eqref{2022306232141} with $i=0$ from \eqref{202230623aa2141} and \eqref{202230623aa2141x}.
\hfill $\Box$
\end{pf}
\begin{lem}
\label{201807312051}
 Estimates involving the fractional differential operators:
\begin{align}
\label{201806291928}
& \int_{\mathbb{R}^2}\|\mathfrak{D}_{\mf{h}}^{3/2} \varphi\|_0^2\mm{d}\mf{h}\lesssim \|\varphi\|_{1}^2 \mbox{ for any }\varphi\in H^1,\\
& \int_{\mathbb{R}^2} |\mathfrak{D}_{\mf{h}}^{s} \phi|^2_{1/2}\mm{d}\mf{h}\lesssim |\phi|_{s+1/2}^2 \mbox{ for any }\phi:=\phi(y_1,y_2,0)\in H^{s+1/2}(\mathbb{R}^2),  \label{201811132119}
\end{align}
\end{lem}
\begin{pf}    In view of  Fubini's theorem and the trace estimate \eqref{2022306232141}, we can
obtain
\begin{align*} \int_{\mathbb{R}^2}
 \|\mathfrak{D}_{\mf{h}}^{3/2} \varphi\|_0^2\mm{d}\mf{h}= \int_{h_-}^{h_+} \int_{\mathbb{R}^2}
\int_{\mathbb{R}^2}|\mathfrak{D}_{\mf{h}}^{3/2} \varphi|_0^2\mm{d}y_{\mm{h}}\mm{d}\mf{h}\mm{d}y_3\lesssim  \int_{h_-}^{h_+}|\varphi|^2_{1/2}\mm{d}y_3\lesssim \|\varphi\|_{1}^2,
\end{align*}
which yields  \eqref{201806291928}.

 It is easy to see that
\begin{align}\int_{\mathbb{R}^2}\sin^2(\xi\cdot \mathbf{h}/2) |\mf{h}|^{-3}\mm{d}\mf{h} = & \int_0^\infty\int_0^{2\pi}\sin^2(|\xi|r\cos(\xi,\mathbf{h})/2) |r|^{-2}\mm{d}\theta\mm{d}r \nonumber \\
\leqslant & 2\pi\left(\int_{r<|\xi|^{-1}}r^{-2}\sin^2(|\xi|r/2)\mm{d}r+\int_{|\xi|^{-1}<r<1}
r^{-2}\mm{d}r+\int_{1<r}
r^{-2}\mm{d}r\right)\nonumber \\
\leqslant & 2\pi   |\xi|\left(1+\int_{0}^1r^{-2}\sin^2(r/2)\mm{d}r\right) \leqslant 4\pi|\xi|,
  \end{align}
  where $\cos(\xi,\mathbf{h})$ denotes  the cosine of angle between $\xi$  and $\mathbf{h}$.
 By  the above estimate, Parseval's theorem of Fourier transform defined on $\mathbb{R}^2$ and  the equivalent definition of the norm $|\cdot|_{1/2}$ by Fourier transform in  \cite[Theorems 7.63]{RAADSBS}, we have
\begin{align} \int_{\mathbb{R}^2}
|\mathfrak{D}_{\mf{h}}^{s} {\phi}|^2_{1/2}\mm{d}\mf{h} \lesssim&\int_{\mathbb{R}^2}\int_{\mathbb{R}^2}
(1+|\xi|^{2})^{s}| \widehat{\mathfrak{D}_{\mf{h}}^{3/2} {\phi}}|^2\mm{d}\xi\mm{d}\mf{h}\nonumber \\
\lesssim& \int_{\mathbb{R}^2}\int_{\mathbb{R}^2}(1+|\xi|^{2})^{s}|e^{-\xi\cdot \mathbf{h}~\mm{i}}
-1|^2|\hat{\phi}|^2 |\mf{h}|^{-3}\mm{d}\xi\mm{d}\mf{h}\nonumber \\ \lesssim & \int_{\mathbb{R}^2}(1+|\xi|^{2})^{s} |\hat{\phi}|^2 \int_{\mathbb{R}^2}\sin^2(\xi\cdot \mathbf{h}/2) |\mf{h}|^{-3}\mm{d}\mf{h}\mm{d}\xi\nonumber \\
 \lesssim& \int_{\mathbb{R}^2}
(1+|\xi|^{2})^{s+1/2} |\hat{\phi}|^2  \mm{d}\xi \lesssim |\phi|_{s+1/2 },
\end{align}
which yields \eqref{201811132119}.
\hfill $\Box$
\end{pf}
\begin{lem}\label{20180812}
Dual estimates: If $\varphi$ and $\psi\in H^{1/2}$, then
\begin{equation}
\label{201808121426}
\left|\int_{\mathbb{R}^2}\partial_{\mm{h}}\varphi \psi\mm{d}y_\mm{h}\right|\lesssim |\varphi|_{1/2} |\psi|_{1/2}.
\end{equation}In particular,
\begin{equation}
\label{201808121247}
|\partial_{\mm{h}}\varphi|_{-1/2}\lesssim |\varphi|_{1/2}.
\end{equation}
\end{lem}
\begin{pf} Using Parseval's theorem of  Fourier transform,   the equivalent definition of the norm $|\cdot|_{1/2}$ by Fourier transform in  \cite[Theorems 7.63]{RAADSBS}  and H\"older's inequality, we immediately get
\eqref{201808121426}, which immediately yields \eqref{201808121247}.
 \hfill $\Box$
\end{pf}

\begin{lem}\label{2022307241635} Let $f\in H^{1}$, then there exists a function $F$ such that
$$F\in C^0(\mathbb{R}^+,H^{1}(\mathbb{R}^2))\cap L^2(\mathbb{R}^+,H^{3/2}),\ F_t\in L^2(\mathbb{R}^+,H^{1/2}) \mbox{ and }F|_{t=0}=f.$$
Moreover,
$$\|F\|_{L^\infty(\mathbb{R}^+,H^{1}(\mathbb{R}^2))}+\|F\|_{L^2(\mathbb{R}^+,H^{3/2})}+\|F_t
\|_{L^2(\mathbb{R}^+,H^{1/2})}\lesssim |f|_{1}.$$
\end{lem}
\begin{pf}
Let $\varphi \in C_0^\infty(\mathbb{R})$ be such that
$\varphi (0)=1$.
We then define  $\hat{F}(\xi,t)=\varphi(t\langle\xi\rangle)\hat{f} (\xi) $,
 where $\hat{\cdot}$ denotes the Fourier transform defined on $\mathbb{R}^2$ and
$\langle\xi\rangle:=\sqrt{1+|\xi|^2}$. By construction, ${F}(\cdot,0)=
f$ and $\partial_t\hat{F} (\xi,t)=\varphi'(t\langle\xi\rangle )\hat{f} (\xi)\langle\xi\rangle$.

 We can estimate that	
\begin{equation}
 | F(\cdot,t) |_{1}^2  =
\int_{\mathbb{R}^{2}}
\big|\varphi(t\langle\xi\rangle)\big|^{2}\big|\hat{f} (\xi)\big|^{2} \langle\xi\rangle^2 \mm{d}\xi \leqslant\big\|\varphi
\big\|_{L^{\infty}}^{2}\big|f \big|_{1}^{2} \label{202220202012}
\end{equation}	
and
\begin{equation}
 | \partial_t F(\cdot,t) |_{0}^2  =
\int_{\mathbb{R}^{2}}
\big|\varphi'(t\langle\xi\rangle)\big|^{2}\big|\hat{f} (\xi)\big|^{2} \langle\xi\rangle^2 \mm{d}\xi \leqslant\big\|\varphi'
\big\|_{L^{\infty}}^{2}\big|f \big|_{1}^{2}. \nonumber
\end{equation}	
We easily further observe from the above two estimates  that
\begin{align}
\label{2022230811120}
F\in C^0(\mathbb{R}^+_0,H^{1})\mbox{ and }F_t\in C^0(\mathbb{R}^+_0,L^2)
\end{align}  by using Lebesgue's dominated convergence theorem. Similarly, for $i=0$, and $1$,
\begin{align}
\int_0^\infty|\partial_t^i F |_{3/2-i }^{2}\mm{d}t& =\int_{0}^{\infty}\int_{\mathbb{R}^{2}}
\langle\xi\rangle^{2i}\big|\varphi ^{(i)}(t\langle \xi \rangle)\big|^{2}\big|\hat{f} \langle\xi\rangle\big|^{2}\langle\xi\rangle^{3-2i}\mm{d}\xi \mm{d}t \nonumber   \\
&=\int_{0}^{\infty}\int_{\mathbb{R}^{2}}\big|\varphi ^{(i)}(t\langle\xi\rangle)\big|^{2}\big|\hat{f} (\xi)\big|^{2}\langle\xi\rangle^{3}\mm{d}\xi \mm{d}t \nonumber \\
&=\int_{\mathbb{R}^{2}}\left|\hat{f} (\xi)\right|^{2} \langle\xi\rangle ^{3}\bigg(\int_{0}^{\infty}\big|\varphi ^{(i)}(t\langle\xi\rangle)\big|^{2}\mm{d}t\bigg)\mm{d}\xi  \nonumber \\
&=\int_{\mathbb{R}^{2}}\big|\hat{f} (\xi)\big|^{2}\langle\xi\rangle^{3}\bigg(\frac{1}{\langle\xi\rangle}\int_{0}^{\infty}\big|\varphi ^{(i)}(r)\big|^{2}\mm{d}r\bigg)\mm{d}\xi  \nonumber  \\
&=\left\|\varphi ^{(i)}\right\|_{L^{2}(\mathbb{R}^+)}^{2}\int_{\mathbb{R}^{2}}\left|\hat{f} (\xi)\right|^{2}\langle\xi\rangle^2 \mm{d}\xi=\left\|\varphi ^{(i)}\right\|_{L^{2}(\mathbb{R})}^{2}\left|f \right|_{1}^{2},  \label{2020308112034}
\end{align}
where $\varphi^{(0)}=\varphi$ and $\varphi^{(1)}=\varphi'$.
Consequently, we immediately see from \eqref{202220202012}--\eqref{2020308112034} that $F$ satisfies the desired conclusion stated in Lemma \ref{2022307241635} (referring to \cite[Lemma A.10]{ZHAOYUI}). \hfill $\Box$
\end{pf}

\begin{lem}
\label{201809012320} Homeomorphism theorem:
Let $k\geqslant 3$. There is a constant $\iota$ depending on $\Omega$, such that for any $\eta\in H_0^1\cap H^{k}$ satisfying
$\|\eta\|_3\leqslant \iota\in (0,1]$, we have
 \begin{align}
&\det\nabla\zeta(y),\ \det\nabla_{\mm{h}}\zeta_{\mm{h}} (y_{\mm{h}},0),\ H^{\mm{d}} \geqslant  {1}/{2}\mbox{ for any }y_{\mm{h}}\in \mathbb{R}^2, \label{201803121601xx21082109}  \\
&\label{201803121601xx2108} \zeta_{\mm{h}} (y_{\mm{h}},0):\mathbb{R}^2 \to \mathbb{R}^2\mbox{ is a }C^{k-2}\mbox{-diffeomorphic mapping}, \\
&\label{201803121601xx} \zeta  : \overline{\Omega}\to \overline{\Omega} \mbox{ is a homeomorphism mapping},  \\
& \zeta _\pm : \Omega_\pm \to \zeta _\pm(\Omega_\pm)\mbox{ are }C^{k-2}\mbox{-diffeomorphic mappings}, \label{2018031adsadfa21601xx}
\end{align}
where $\zeta:=\eta+y$, $(\zeta_{\mm{h}} )^{-1}$ denotes the inverse mapping of $\zeta_{\mm{h}}(y_{\mm{h}},0)$, $\zeta_{\mm{h}}$
represents the first two components of $\zeta :={\eta} +y$ and $H^{\mm{d}} = |\partial_1\zeta|^2|\partial_2\zeta|^2-|\partial_1\zeta\cdot \partial_2\zeta|^2$.
\end{lem}
\begin{pf} Please refer to \cite[Proposition 5.2]{JFJSZWC}. \hfill $\Box$
\end{pf}
\begin{lem}  \label{xfsddfsf201805072212}
 Existence theory of the stratified Lam\'e problem with Dirichlet boundary conditions:
 let $i\geqslant  0$, $\mathbf{F}^1\in H^{i }$ and $\mathbf{F}^2\in H^{i+3/2} $, then there exists a unique solution $u\in H^{i+2} $
of the following  Lam\'e problem:
 \begin{equation}
 \begin{cases}
\mu\Delta u+(\varsigma+\mu/3)\nabla \mm{div}u= \mathbf{F}^1&\mbox{in }  \Omega, \\[1mm]
u_+=u_-=\mathbf{F}^2&\mbox{on }  \Sigma,\\
 u=0 &\mbox{on }\partial\Omega\!\!\!\!\!-    ;
\end{cases}\label{201808311539}
\end{equation}
moreover,
\begin{equation}
\label{xfsddfsf201705141252}
\|u\|_{i+2}\lesssim
\|\mathbf{F}^1\|_{i}+|\mathbf{F}^2|_{i+3/2}.
\end{equation}
\end{lem}
\begin{pf} Both the  results of  existence and regularity for unique solutions of    the following  horizontally periodic stratified Lam\'e problem
\begin{equation}
\label{2023306232004}
\begin{cases}
\mu\Delta u+(\varsigma+\mu/3)\nabla \mm{div}u= {\mathbf{F}}^1&\mbox{in }\Omega_{L_1,L_2} , \\[1mm]
u_+=u_-=\mathbf{F}^2 &\mbox{on }  \Sigma_{L_1,L_2},\\
u=0 &\mbox{on }\partial\Omega\!\!\!\!\! -
\end{cases}\end{equation}
had been proved by Jang--Tice--Wang in \cite[Proposition E.4]{jang2016compressible}, see \eqref{2022307222447} for the definitions of
$\Omega_{L_1,L_2}$ and $\Sigma_{L_1,L_2}$.
 Following Jang--Tice--Wang's augment, we easily extend the results of the periodic case in \cite[Proposition E.4]{jang2016compressible} to our non-periodic case stated in Lemma \ref{xfsddfsf201805072212}. \hfill $\Box$
\end{pf}

\begin{lem}\label{2022306231855}
 Existence theory of the stratified Lam\'e problem with jump conditions:  let $i\geqslant  0$, $\mathbf{F}^1\in H^{i }$ and $\mathbf{F}^2\in H^{i+1/2}$, then there exists a unique solution $u\in H^{i+2}$
of the following stratified Lam\'e problem:
\begin{equation*}\begin{cases}
\mu\Delta u+(\varsigma+\mu/3)\nabla \mm{div}u= \mathbf{F}^1&\mbox{in }  \Omega, \\[1mm]
  \llbracket  u  \rrbracket =0,\  \llbracket \mathbb{S}(u)\mathbf{e}^3\rrbracket= \mathbf{F}^2&\mbox{on }\Sigma, \\
u=0 &\mbox{on }\partial\Omega\!\!\!\!\!\!-.
\end{cases}\end{equation*}
Moreover,
\begin{equation}
\label{Ellipticestimate}
\|u\|_{i+2}\lesssim
\|\mathbf{F}^1\|_{i}+|\mathbf{F}^2|_{i+1/2}.
\end{equation}
\end{lem}
\begin{pf} Both the  results of  existence and regularity for unique solutions of   the  following  horizontally periodic  stratified Lam\'e problem
\begin{equation}
\label{202330623sdfa2004xx}
\begin{cases}
\mu\Delta u+(\varsigma+ \mu/3)\nabla \mm{div}u=  {\mathbf{F}}^1&\mbox{in }\Omega_{L_1,L_2}  , \\[1mm]
  \llbracket  u  \rrbracket =0,\  \llbracket \mathbb{S}(u) \mathbf{e}^3\rrbracket= {\mathbf{F}}^2&\mbox{on }\Sigma_{L_1,L_2} , \\
   \mathbb{S}(u)   \mathbf{e}^3= \mathbf{F}^3&\mbox{on }\Sigma_+^{L_1,L_2}:=2\pi L_1\mathbb{T}\times 2\pi L_1\mathbb{T}\times \{h_+\}, \\
u=0 &\mbox{on }\partial\Omega\!\!\!\!\!  -
\end{cases}\end{equation}
had been provided  by Jang--Tice--Wang  in \cite[Lemma A.10]{JJTIWYJ}. Following Jang--Tice--Wang's augment, we easily extend the results of the periodic case in  \cite[Lemma A.10]{JJTIWYJ} to our non-periodic case stated in Lemma \ref{2022306231855}. \hfill $\Box$
\end{pf}

\renewcommand\thesection{Appendix B}
\section{Regularity for the solutions in Eulerian coordinates}\label{secsd:09}
\renewcommand\thesection{B}

This section is devoted to the derivation of regularity for solutions, which are obtained from Theorem \ref{thm:0202} by the inverse transform of Lagrangian coordinates, in Eulerian coordinates. In particular, we provide the derivation of \eqref{20181114101015}  and \eqref{20181114101015xx}  for $\vartheta>0$ in Section \ref{2022307291526}.
In what follows, we always assume that the solution $(\eta,u)\in C^0([0,T],H_{0,*}^{3,1/2} \times H^2)$ is provided by Theorem \ref{thm:0202} for given $\delta$. In addition, for the sake of the simplicity, we define that $\Omega^T_\pm:=\Omega_\pm\times (0,T)$, $\Omega^T:=\Omega\times (0,T)$, $\mathbb{R}^2_T:=\mathbb{R}^2\times (0,T)$  and  the closure of a set $S$ by  $\overline{S}$.

\subsection{Homeomorphism}\label{201811191634}
Let $\zeta=\eta+y$, $\tilde{y}:=(y,t)$, $\tilde{y}_{\mm{h}} :=(y_{\mm{h}},t)$, $\tilde{x}=(x,t)$ and $\tilde{x}_{\mm{h}}=(x_{\mm{h}},t)$. Since $\eta(t)\in H^{3,1/2}_{0,*}$ and $\eta\in C([0,T],H^3)$, we have
\begin{align}
\label{201806031606}
&\nabla_x \zeta^{-1} =(\nabla_y \zeta)^{-1}|_{y=\zeta^{-1} }=\mathcal{A}^{\mm{T}}|_{y=\zeta^{-1} }\mbox{ in }\Omega,\\
& \tilde{\zeta}: \overline{\Omega^{T}} \mapsto  \overline{\Omega^{T}} \mbox{ is a bijective mapping}, \label{201809300836} \\
 & \tilde{\zeta}_\pm: X \mapsto \tilde{\zeta}_\pm(X)\mbox{ are bijective mappings}, \label{201809300834}\\
& \det \nabla_{\tilde{y}}\tilde{\zeta}(\tilde{y})\geqslant1/2\mbox{ in }\overline{\Omega_\pm^T}, \label{201810011535} \\
&  \bar{\zeta}: Y \to  Y\mbox{ is a bijective mapping},\label{2018081518dsfgs5sdfa6}\\
&\det \nabla_{\tilde{y}_\mm{h}}\bar{\zeta}(\tilde{y}_{\mm{h}})=\det \nabla_{y_{\mm{h} }} \zeta_{\mm{h}}(y_{\mm{h}},0,t)\geqslant 1/2, \label{2018081518dsfgs5sdfa6xx}
 \end{align}
 where $\tilde{\zeta} (\tilde{y}):=(\zeta (y,t),t)$, $\bar{\zeta}(\tilde{y}_{\mm{h}}):=(\zeta_{\mm{h}}(y_{\mm{h}},0,t),t)$, $X =\overline{\Omega_\pm^{T}}$ or $\Omega_\pm^{T}$, and $Y =\overline{\mathbb{R}^2_{T}} $ or ${\mathbb{R}^2_{T}}$.

We denote the inverse functions of $\tilde{\zeta}_\pm(\tilde{y})$ resp. $\bar{\zeta}(\tilde{y}_{\mm{h}})$ by $\tilde{\zeta}_{\pm}^{-1}(\tilde{x})$
resp. $\bar{\zeta}^{-1}(\tilde{x}_{\mm{h}})$.
Recalling the regularity
\begin{align}
&  \zeta-y\in C^0([0,T],H^3) \mbox{ and }\zeta_t=u\in C^0([0,T],H^2),\label{201810011540x}
\end{align}
we use the embedding theorem $H^{k+2}(\Omega_\pm)\hookrightarrow C^{k}(\overline{\Omega}_\pm)$ for $k\geqslant 0$ to get
\begin{align}  & \label{201809271654}
\tilde{\zeta}_\pm \in C^1(\overline{\Omega_\pm^T}),\\
& \label{201809271654xx}   \tilde{\zeta} \in C^1(\overline{\mathbb{R}^T_2}).
 \end{align}
By \eqref{n0101nnnM}$_4$, we get \begin{equation}
\label{20180927165xx4}
\tilde{\zeta},\ u, \nabla_{y_{\mm{h}}} {\zeta} \in C^0(\overline{\Omega ^T}).
 \end{equation}

We further derive
 from \eqref{201809300836} and the continuity of $\tilde{\zeta}$ in \eqref{20180927165xx4} that (referring to (5.57) in \cite{JFJSZWC})
\begin{equation}  \label{201810012335}
\tilde{\zeta} (\tilde{y}) :\overline{\Omega ^T}\to  \overline{\Omega ^T}\mbox{ is a homeomorphism mapping}.
\end{equation}
  Similarly to \eqref{2018031adsadfa21601xx} and \eqref{201810012335}, we have by \eqref{201809300834}, \eqref{201810011535}
  and \eqref{201809271654} that
 \begin{align}
& \tilde{\zeta}_\pm(\tilde{y}):\overline{ \Omega_\pm^T} \to \tilde{\zeta}_\pm( \overline{\Omega_\pm^T})\mbox{ are   } C^1\mbox{-diffeomorphic mappings}  .\label{201808151856sdfaxxxxx}
\end{align}
Moreover, $\nabla_{\tilde{x}}\tilde{\zeta}^{-1}_\pm
=(\nabla_{\tilde{y}}\tilde{\zeta}_\pm)^{-1}|_{\tilde{y}=\tilde{\zeta}_\pm^{-1}}$. In particular,
\begin{equation}
\label{2018060321637}
\partial_t \zeta^{-1}_\pm=-((\nabla_y \zeta_\pm)^{-1}u_\pm)|_{y=\zeta^{-1}_\pm}=-(\mathcal{A}^{\top}_\pm u_\pm)|_{y=\zeta^{-1}_\pm}.
\end{equation}

Similarly  to \eqref{201806031606} and \eqref{201808151856sdfaxxxxx}, we deduce from \eqref{2018081518dsfgs5sdfa6}--\eqref{2018081518dsfgs5sdfa6xx}  and \eqref{201809271654xx}  that
\begin{align}
&\bar{\zeta}(y_{\mm{h}},t): \overline{\mathbb{R}^2_{T}} \to   \overline{\mathbb{R}^2_{T}}\mbox{ is a }C^1\mbox{-diffeomorphic mapping},\label{2018081518dsfgs56}\\
&
\nabla_{x_\mm{h}} (\zeta_{\mm{h}})^{-1}(x_{\mm{h}},t)=(\nabla_{y_\mm{h}} \zeta_{\mm{h}}(y_{\mm{h}},0,t))^{-1}|_{y_{\mm{h}}=(\zeta_{\mm{h}})^{-1}(x_{\mm{h}},t)}\nonumber \\
& =\left.\left(\frac{1}{ \det \nabla_{y_{\mm{h} }} \zeta_{\mm{h}}(y_{\mm{h}},0,t)  } \left( \begin{array}{cc }
  \partial_2 \zeta_2 (y_{\mm{h}},t)    &    -\partial_2 \zeta_1 (y_{\mm{h}},t)  \\
  -\partial_1 \zeta_2 (y_{\mm{h}},t)  &    \partial_1 \zeta_1 (y_{\mm{h}},t)
                                                    \end{array}\right) \right)\right|_{y_{\mm{h}}=(\zeta_{\mm{h}})^{-1}(x_{\mm{h}},t)}
.
\label{201809291022}\end{align}

\subsection{Regularity of solutions in Eulerian coordinates}\label{201811142125}

Let $a:=\|\eta\|_{C^0([0,T],H^3)}$, where $a\in  (0,\iota]$. In what follows, the notation
\begin{equation*}
A \lesssim_a B\mbox{  means }A \leqslant c(a) B,
 \end{equation*} where  $
 c(a)$ denotes   a generic positive  constant, which may vary from line to line,  depends on $a$, $\Omega$ and increases with respect to $a$.

Thanks to the homeomorphism properties of $\zeta$ and $\zeta_{\mm{h}}$, the following definitions make sense.
\begin{align}&\rho:= (\bar{\rho}J^{-1})|_{y={\zeta}^{-1}(x,t)},\  {v}:={u}({\zeta}^{-1}(x,t),t),\ \tilde{\rho}:= \bar{\rho}|_{y_3={\zeta}^{-1}_3(x,t))}, \nonumber\\
&d:=\zeta_3((\zeta_{\mm{h}})^{-1}(x_{\mm{h}},t),0,t) \in (h_-,h_+),\
\nu:=(-\partial_{x_1} d,-\partial_{x_2} d,1 )^{\top}/\sqrt{1+|\nabla_{x_{\mm{h}}} d|^2}, \label{202223020155} \\
&
  \Sigma(t):=\{(x_{\mm{h}},x_3)~|~x_{\mm{h}}\in \mathbb{R}^2,\ x_3:=d(x_{\mm{h}},t)\},\nonumber \\
&\Omega_+(t):=\{(x_{\mm{h}},x_3)~|~x_{\mm{h}}\in \mathbb{R}^2  ,\ d(x_{\mm{h}},t)< x_3<h_+\},\nonumber \\
&\Omega_-(t):=\{(x_{\mm{h}},x_3)~|~x_{\mm{h}}\in \mathbb{R}^2 ,\ h_-< x_3<d(x_{\mm{h}},t)\}. \nonumber
\end{align}

\subsubsection{Motion domains}
By \eqref{201803121601xx2108} and the  fact $$ \zeta_3((\zeta_{\mm{h}})^{-1}(x_{\mm{h}},t),0,t)= \eta_3((\zeta_{\mm{h}})^{-1}(x_{\mm{h}},t),0,t),$$ for any given $x_{\mm{h}}$, there is $y_{\mm{h}}$, such that
\begin{equation}
\label{xx201809282115}
x_{\mm{h}}:=\zeta_{\mm{h}}(y_{\mm{h}},0,t)=\eta_{\mm{h}}(y_{\mm{h}},0,t)\mbox{ and }
d(x_{\mm{h}},t)=\eta_3(y_{\mm{h}},0,t).
\end{equation}
This means $ \Sigma(t)\subset\zeta(\Sigma,t )$. Similarly, we also have $\zeta(\Sigma,t )\subset \Sigma(t)$. Consequently, we arrive at that, for any given $t\geqslant 0$,
\begin{equation}
\label{201810011624}
 \zeta(\cdot,t): \{y_3=0\}\to \{x_3=d(x_{\mm{h}},t)\}\mbox{ is a bijective mapping}.
\end{equation}
Moreover, we further deduce that
\begin{align}
&\Omega_\pm(t) =\zeta_\pm(\Omega_\pm,t),\ \Omega\!\!\!\!-=\Omega (t)\cup \Sigma(t),\ \Omega_+(t)\cap \Omega_-(t)=\emptyset \mbox{ and }\Omega_\pm(t)\cap \Sigma(t)=\emptyset,
\quad t\geqslant 0, \nonumber
\end{align}
where $\Omega(t):=\Omega_+(t)\cup \Omega_-(t)$.
In what follows, if $f$ is defined in $\Omega(t)$, resp. $\Omega$,  then $f_\pm:=f|_{\Omega_\pm(t)}$, resp. $f|_{\Omega_\pm}$.

\subsubsection{Regularity and motion equations of $(\rho,v)$}
 Thanks to the continuity of $(\nabla_y\zeta_\pm, u)$ in \eqref{201809271654} and \eqref{20180927165xx4}--\eqref{201808151856sdfaxxxxx},
$v\in C^0(\overline{\Omega^T})$ and   $(\rho_\pm,P_\pm(\rho_\pm))$ $\in (C^0(\overline{\Omega^T_\pm}))^2$. Obviously, $v=0$ on $\partial\Omega\!\!\!\!\!-$ due to $u=\eta=0$ on $\partial\Omega\!\!\!\!\!-$.

By transform of Lagrangian coordinates (i.e., $x=\zeta(y)$) and the regularity of $(\zeta,u)$, we can bound that
for any given $t\geqslant 0$,
\begin{align}  \label{201811210953}
\|v(t)\|_{L^2(\Omega (t))}^2= & \int |u|^2J\mm{d}y   \lesssim_a\|u(t)\|_0^2.
\end{align}
Noting that $\partial_{x_i} v= (\mathcal{A}_{il}\partial_{y_l} u )|_{y=\zeta^{-1}(x)}$ for $1\leqslant i\leqslant 3$, we have
\begin{align}
\label{201811210954}
\|\partial_{x_i}v(t)\|_{L^2(\Omega (t))}^2= &
\int|
\mathcal{A}_{il}\partial_{y_l} u |^2J\mm{d}y \lesssim_a\|u(t)\|_1^2 .
\end{align}

Similarly, we can further derive that for any $1\leqslant i,j,k\leqslant 3$,
$$ \|\partial_{x_ j}\partial_{x_i}v(t)\|_{L^2(\Omega (t))}  \lesssim_a  \|u(t)\|_2, \quad
\|\partial_{x_ k}\partial_{x_ j}\partial_{x_i}v(t)\|_{L^2(\Omega(t))}\lesssim_a \|u(t)\|_3 $$
by virtue of the relations
 \begin{align*}
 \partial_{x_j}\partial_{x_i}v=(\mathcal{A}_{jm} \partial_{y_m}(\mathcal{A}_{il}\partial_{y_l} u ))|_{y=\zeta^{-1}}\mbox{ and }\partial_{x_k}\partial_{x_j}\partial_{x_i}v=(\mathcal{A}_{kn} \partial_n (\mathcal{A}_{jm}
\partial_{y_m}(\mathcal{A}_{il}\partial_{y_l} u )))|_{y=\zeta^{-1}}.
\end{align*}
Therefore,
\begin{align}
& \|v(t)\|_{H^2(\Omega(t))}\lesssim _a \|u\|_2\mbox{ for any  }t\in [0,T],\label{2018092137}\\
&\|v(t)\|_{H^3(\Omega(t))}\lesssim _a \|u\|_3\mbox{ for a.e. }t\in (0,T).\label{2018092137zxx} \end{align}

Recalling \eqref{201611050816}, we have
\begin{align} P_\pm(\rho_\pm(t))-P_\pm(\tilde{\rho}_\pm(t))= ( R_P-{P}'(\bar{\rho})\bar{\rho}\mm{div}\eta)_\pm|_{y={\zeta}^{-1}_\pm(x,t)}.
\label{1141558xx}  \end{align}
Thus, similarly  to \eqref{201612071026} and \eqref{2018092137}, we can also get
\begin{align*}
  &\|P_+(\rho_+(t))-P_+(\tilde{\rho}_+(t))\|_{H^{2}(\Omega_+(t))}+
  \|P_-(\rho_-(t))-P_-(\tilde{\rho}_-(t))\|_{H^{2}(\Omega_-(t))}\\
  &\lesssim _a
  \|R_P-{P}'(\bar{\rho})\bar{\rho}\mm{div}\eta \|_2\lesssim _a\|\eta(t)\|_{3}.
 \end{align*} Similarly, we can also verify that
 $$\| \rho (t) -\tilde{\rho} (t)\|_{H^{2}(\Omega_\pm(t))} \lesssim _a
  \|\eta(t)\|_{3}$$
 and
 $\rho_\pm(t)$, $P_\pm(\rho_\pm(t))  \in H^2 (\Omega_\pm(t)\cap B_r)$ for any given $r>0$, where $B_r=\{x~|~x\in \mathbb{R}^3~|~|x|\leqslant r\}$.

By the definitions of $\rho$ and $v$, we can employ  \eqref{201806031606}, \eqref{2018060321637}, the relation $J_t=J\mm{div}_{\mathcal{A}}u$
and the chain rule of differentiation to derive that
\begin{align}
\rho_t+ \mm{div}(\rho v) = ((\bar{\rho}J^{-1})_t+ \bar{\rho}J^{-1}\mm{div}_{\mathcal{A}}u)|_{y={\zeta}^{-1}(x,t))}=0.\label{201811142017cc}
\end{align}
Similarly, we can easily get from \eqref{n0101nnnM}$_2$ that $(\rho,v) $ satisfies
  $$ \rho v_t+\rho v\cdot\nabla v+\mm{div}\mathcal{S}=-\rho g \mathbf{e}^3.$$

\subsubsection{Regularity and motion equations at interface}\label{2022308031213}

 For any given $t\geqslant 0$, $\Xi(x ,t):=((\zeta_{\mm{h}})^{-1}(x_{\mm{h}},t),x_3):\Omega_+\to \Omega_+$ is a bijective mapping.
We denote the inverse function of $\Xi(x ,t)$ by $\Xi^{-1}(y ,t)$.  Then $\Xi$ and $\Xi^{-1}$ enjoy the following properties:
\begin{align}
& \Xi^{-1}(y ,t) = (\zeta_{\mm{h}}(y_{\mm{h}},0,t),y_3 ), \nonumber \\
& \nabla_x \Xi(x ,t) \nonumber \\
&=
\left(\frac{1}{\det \nabla_{y_{\mm{h} }} \zeta_{\mm{h}}(y_{\mm{h}},0,t)}
\left.\left( \begin{array}{ccc}
  \partial_2 \zeta_2 (y_{\mm{h}},t)    &    -\partial_2 \zeta_1 (y_{\mm{h}},t)  & 0 \\
  -\partial_1 \zeta_2 (y_{\mm{h}},t)  &    \partial_1 \zeta_1 (y_{\mm{h}},t)  & 0 \\
 0 & 0 &\Theta(y_{\mm{h}},t)
                                                    \end{array}\right)\right)
\right|_{y_{\mm{h}}=(\zeta_{\mm{h}})^{-1}(x_{\mm{h}},t)}
,\nonumber \\
& \label{2020307saf302204} 1/2\leqslant  \det \nabla_{y_{\mm{h} }} \zeta_{\mm{h}}(y_{\mm{h}},0,t)  =\det  \nabla_y  \Xi^{-1}(y ,t)\lesssim_a1,
\end{align}

Similarly to \eqref{2018092137zxx}, we easily get that
\begin{align}
\|\eta_3^+(y,t)|_{y=\Xi(x ,t)}\|_{H^3(\Omega_+)}\lesssim _a \|\eta(t)\|_3,
\label{20223080322257}
\end{align}
where $\eta_3^+$ denote the third component of $\eta_+$.
 Exploiting the trace estimate,
we can deduce from \eqref{20223080322257} that
\begin{equation}  \label{201809292130}
|d(x_\mm{h},t)=(\eta_3^+(y,t)|_{y=\Xi(x ,t)})|_{x_3=0}|_{5/2}\lesssim _a \|\eta(t)\|_3.
\end{equation}
Moreover, from \eqref{201809291022}  we can obtain
\begin{align}
 \nabla_{x_{\mm{h}}}d= & (\nabla_{x_\mm{h}}(\zeta_{\mm{h}})^{-1}(x_{\mm{h}},t))^{\top} \nabla_{y_\mm{h}}\zeta_3
(y_{\mm{h}},0,t)|_{y_{\mm{h}}=(\zeta_{\mm{h}})^{-1}(x_{\mm{h}},t)} \nonumber \\
=& ((\nabla_{y_\mm{h}} \zeta_{\mm{h}}(y_{\mm{h}},0,t))^{-\top} \nabla_{y_\mm{h}}\zeta_3(y_{\mm{h}},0,t))|_{y_{\mm{h}}=
(\zeta_{\mm{h}})^{-1}(x_{\mm{h}},t)}   \label{201810012320}\end{align}
and
\begin{align}
&\nabla_{x_{\mm{h}}}\partial_{x_i}d =((\nabla_{y_\mm{h}} \zeta_{\mm{h}}(y_{\mm{h}},0,t))^{-\top}\nabla_{y_\mm{h}} ( (\nabla_{y_\mm{h}} \zeta_{\mm{h}}(y_{\mm{h}},0,t))^{-\top}
\nabla_{y_\mm{h}}\zeta_3(y_{\mm{h}},0,t))_i) |_{y_{\mm{h}}
=(\zeta_{\mm{h}})^{-1}(x_{\mm{h}},t)} ,  \label{201810012320xx}
\end{align}
where $(f)_i$ denotes the $i$-th component of $f$ for $i=1$  and $2$.
In addition, making use of the continuity of $(\zeta,\nabla_{y_{\mm{h}}} {\zeta})$ in \eqref{2018081518dsfgs56}, \eqref{20180927165xx4}
and \eqref{201810012320}, we see that
$\nabla_{\mm{h}}d$ is continuous on $\overline{\mathbb{R}^2_T}$.

Analogously to \eqref{201809292130}, we obtain $|1- ( 1+|\nabla_{x_{\mm{h}}} d(x_{\mm{h}},t)|^2 )^{-1/2}|_{3/2}\lesssim _a 1$,
which, together with \eqref{202307231325} with $j=1$ and \eqref{201809292130}, yields
\begin{equation}
\nu(x_{\mm{h}},t)-\mathbf{e}^3 \in H^{3/2} \;\;\mbox{  for any }t\geqslant 0.
\label{201811201627}
\end{equation}

Obviously, one has
\begin{align}
&u_3|_{y_3=0}=\partial_t \eta_3(y_{\mm{h}},0,t)=\partial_t d (\zeta_{\mm{h}}(y_{\mm{h}},0,t),t)\nonumber \\
&\qquad \ \ \ = (u_1|_{y_3=0} \partial_1d+u_2 |_{y_3=0} \partial_2d +d_t)|_{x_{\mm{h}}=\zeta_{\mm{h}}(y_{\mm{h}},0,t)},
\label{201810011601}\\
&
  \label{201810011617}
\mbox{ and } v(x_{\mm{h}},d(x_{\mm{h}},t),t)|_{x_{\mm{h}}= \zeta_{\mm{h}} (y_{\mm{h}},t)}= v(\zeta_{\mm{h}}(y_{\mm{h}},0,t),\zeta_3(y_{\mm{h}},0,t),t)
=u(y_{\mm{h}},0,t).
\end{align}
Plugging $y_{\mm{h}}= (\zeta_{\mm{h}})^{-1}(x_{\mm{h}},t)$ into \eqref{201810011601} and
then  using \eqref{201810011617}, we get  for any given $t\geqslant 0$.
\begin{equation}\label{0103nxx}
 d_t+v_1 \partial_1d+v_2 \partial_2d=v_3 \mbox{ on }\Sigma(t).
 \end{equation}

In terms of \eqref{201810011617} and the definition of $\Xi(x ,t)$,
$$v_+(x_{\mm{h}}, d(x_{\mm{h}},t), t) = u_+((\zeta_{\mm{h}})^{-1}(x_{\mm{h}},t),0,t)=(u_+(y,t)|_{y=\Xi(x ,t)})|_{x_3=0}.$$
Thus similarly to \eqref{201809292130}, one gets
\begin{equation}
|v_+(x_{\mm{h}}, d(x_{\mm{h}},t), t)|_{3/2} = |(u_+(y,t)|_{y=\Xi(x ,t)})|_{x_3=0} |_{3/2} \lesssim \| u_+(y,t)|_{y=\Xi(x ,t)}  \|_{2} \lesssim _a\| u\|_2.
\label{2018090232130}
\end{equation}
Thus from \eqref{202307231325} with $j=1$, \eqref{201809292130}, \eqref{0103nxx} and \eqref{2018090232130} it follows that $
d_t\in H^{3/2}$.

Using  relations \eqref{20230626} and \eqref{201810012320}, we obtain after a straightforward calculation that
$$  \begin{aligned}
J\mathcal{A}\mathbf{e}^3|_{y_3=0}= & (-\partial_1\eta_3+\partial_1 \eta_2\partial_2\eta_3-\partial_1\eta_3\partial_2\eta_2,\\
&
-\partial_2\eta_3+\partial_1 \eta_3\partial_2\eta_1-\partial_1\eta_1\partial_2\eta_3,
\det\nabla_{y_\mm{h}}\zeta_{\mm{h}})^{\mm{T}}|_{y_3=0} \\
=&\det\nabla_{y_\mm{h}}\zeta_{\mm{h}}(y_{\mm{h}},0,t)  (- \partial_{x_1}d, - \partial_{x_2}d, 1  )^{\mm{T}}|_{x_{\mm{h}}=\zeta_{\mm{h}}(y_{\mm{h}},0,t)},
\end{aligned} $$
which yields
\begin{align}
\label{201809282114}
&\nu|_{x_{\mm{h}}=\zeta_{\mm{h}}(y_{\mm{h}},0,t)}=\left.\frac{(-\partial_{x_1} d,-\partial_{x_2} d,1)^{\mm{T}}}{\sqrt{1+
|\nabla_{x_\mm{h}}d|^2}} \right|_{x_{\mm{h}}=\zeta_{\mm{h}}(y_{\mm{h}},0,t)} = \left. \frac{J\mathcal{A}\mathbf{e}^3}{|J\mathcal{A}\mathbf{e}^3|}\right|_{y_3=0}.
\end{align}
With the help of   \eqref{201810012320} and \eqref{201810012320xx}, we further obtain that
\begin{align}
\label{20180928212154}
\mathcal{C}|_{x_{\mm{h}}=\zeta_{\mm{h}}(y_{\mm{h}},0,t)}
=\mathcal{H}|_{y_3=0}.
\end{align}

In addition, we have that for  a.e. $t\in (0,T)$,
  \begin{align}
  \left(\mathcal{S}_\pm|_{\{x_3=d(x_{\mm{h}},t)\}}\nu\right)|_{x_{\mm{h}}=
\zeta_{\mm{h}}(y_{\mm{h}},0,t)}
  &=(( P(\bar{\rho}J^{-1})  \mathbb{I}-\mathbb{S}_{\mathcal{A}}(u)  )|_{\Omega_\pm} J \mathcal{A}\mathbf{e}^3   /|J\mathcal{A}\mathbf{e}^3|)|_{y_3=0}.
\label{201809282115}
\end{align}
Consequently, making use of \eqref{n0101nnnM}$_3$,  \eqref{xx201809282115}, \eqref{201810011624},
\eqref{201809282114}--\eqref{201809282115} and the fact $J\geqslant 1$, we get \begin{align}\label{202307251554}
 (\mathcal{S}_+|_{\Sigma(t)} - \mathcal{S}_-|_{\Sigma(t)}) \nu= \vartheta \mathcal{C} \nu\mbox{ on }\Sigma(t).
 \end{align}

\subsection{Higher regularity of $d$ for $\vartheta> 0$}\label{2022308081027}

Similarly to \eqref{201810011617},
\begin{align}((P_\pm(\tilde{\rho}_+(t)))|_{\{x_3=d(x_{\mm{h}},t)\}})|_{x_{\mm{h}}=
\zeta_{\mm{h}}(y_{\mm{h}},0,t)}=
(P_\pm(\bar{\rho}_\pm))|_{y_3=0} ,
\label{202230802120145}
\end{align}
which, together with \eqref{201611051547}$_2$, yields
\begin{align}
  P_+(\tilde{\rho}_+ ) |_{\Sigma(t)} - P_-(\tilde{\rho}_- )  |_{\Sigma(t)}=0.
\end{align}
Thanks to  \eqref{202307251554} and the above identity, we have
\begin{align}\label{20230725sasfaf1554}
 ((\mathcal{S}_+-P_+(\tilde{\rho}_+ ) )|_{\Sigma(t)} -( \mathcal{S}_--P_-(\tilde{\rho}_- ) ))|_{\Sigma(t)}) \nu= \vartheta \mathcal{C} \nu\mbox{ on }\Sigma(t).
 \end{align}
Multiplying \eqref{20230725sasfaf1554} by $\nu$, we find that
\begin{equation}
\label{201811131906}
((\mathcal{S}_+-P_+(\tilde{\rho}_+ ) )|_{\Sigma(t)} -( \mathcal{S}_--P_-(\tilde{\rho}_- ) ))|_{\Sigma(t)}) \nu\cdot \nu  =\vartheta \mathcal{C}\mbox{ on }\mathbb{R}^2,
\end{equation}
where $\mathcal{C}$ can be rewritten as follows \cite{JJTIWYJ}
\begin{align}\nonumber
\mathcal{C} =\mm{div}_{\mm{h}}\left(\frac{\nabla_{\mm{h}}d}{\sqrt{1+|\nabla_{\mm{h}}d|^2}}\right).
\end{align}

Since $|d|_{5/2}\lesssim \|\eta\|_3$, thanks to \eqref{1141558xx}, \eqref{201809282114}, \eqref{201809282115}, \eqref{202230802120145}  and the regularity theory of elliptic equation in      \cite[Lemma 3.2]{TATNLTE1995},  we obtain from \eqref{201811131906} that for sufficiently small $\delta$,
\begin{align}
 |d|_{7/2}\lesssim & |  ((\mathcal{S}_+-P_+(\tilde{\rho}_+ ))|_{\Sigma(t)} - (\mathcal{S}_--P_-(\tilde{\rho}_- ))|_{\Sigma(t)}) \nu\cdot \nu |_{3/2}  \nonumber \\
= &\left|( \llbracket   (R_P-{P}'(\bar{\rho})\bar{\rho}\mm{div}\eta)\mathbb{I} -\mathbb{S}_{\mathcal{A}}(u) \rrbracket J^2\mathcal{A}\mathbf{e}^3\cdot \mathcal{A}\mathbf{e}^3/|J\mathcal{A}\mathbf{e}^3|^2) |_{y=((\zeta_{\mm{h}})^{-1}(x_{\mm{h}},t),0)}\right|_{3/2} .\label{201811211302xx}
\end{align}

Employing the same arguments as for \eqref{2018090232130}, we obtain
\begin{align}
& \left|   \llbracket\mathbb{S}_{\mathcal{A}}(u)  J^2\mathcal{A}\mathbf{e}^3\cdot \mathcal{A}\mathbf{e}^3/|J\mathcal{A}\mathbf{e}^3|  \rrbracket  |_{y_{\mm{h}}=(\zeta_{\mm{h}})^{-1}(x_{\mm{h}},t)}\right|_{3/2} \nonumber \\
&   \lesssim _a\|\mathbb{S}_{\mathcal{A}}(u)  J^2\mathcal{A}\mathbf{e}^3\cdot \mathcal{A}\mathbf{e}^3/|J\mathcal{A}\mathbf{e}^3|\|_2 \lesssim _a \|\nabla u\|_2.
\nonumber
\end{align}
Similarly,
\begin{align}
& \left|   \llbracket   R_P -{P}'(\bar{\rho})\bar{\rho}\mm{div}\eta  \rrbracket |_{y_{\mm{h}}=(\zeta_{\mm{h}})^{-1}(x_{\mm{h}},t)}\right|_{3/2} \lesssim _a   \|  \eta\|_3 .\nonumber
\end{align}
Thus, inserting   the above two estimates into \eqref{201811211302xx}, one gets
\begin{align}\label{2018111401554}
 |d|_{7/2}\lesssim _a  \|(\eta,u)\|_3 .
\end{align}

Analogously to \eqref{2018090232130}, we can obtain
$$ |v_\pm(x_{\mm{h}},d(x_{\mm{h}},t))|_{5/2}\lesssim _a\|u\|_3 \;\;\mbox{ for a.e. }t>0. $$
Utilizing  \eqref{202307231325} with $j=2$,  \eqref{2018090232130} and the above inequality, one can derive
from \eqref{0103nxx}  that
\begin{align}
|d_t|_{5/2}\lesssim_a  & |v^3_+|_{5/2}+|\partial_{x_1} d v^1_++\partial_{x_2} d v^2_+|_{5/2}\nonumber  \\
\lesssim_a  &  |v^3_+|_{5/2}+|d|_{7/2}| v_+|_{3/2}+|d|_{5/2}|v_+|_{5/2}\nonumber \\
\lesssim_a &  \|u\|_{3} + |d|_{7/2}\|u\|_{2}. \label{201811141556}
\end{align}
Moreover, from \eqref{2018111401554} and \eqref{201811141556} one gets
\begin{align}
|d|_{3}^2\lesssim_a &  |d^0|_{3}^2+\int_0^t|d|_{7/2}|d_{\tau}|_{5/2}\mm{d}\tau \nonumber \\
\lesssim_a   &   |d^0|_{3}^2+\int_0^t \left( 1+\|u\|_2   \right) \|(\eta,u)\|_3^2 \mm{d}\tau .
\label{201811141623}
\end{align}
By the regularity of $(\eta,u)$ and \eqref{2018111401554}--\eqref{201811141623}, we have  $
d \in C([0,T],H^{3}(\mathbb{R}^2))\cap L^2((0,T),H^{7/2} )$ and $d_t\in L^2((0,T),H^{5/2} )$ for $\vartheta> 0$.

Finally,
\begin{align}
|\mathcal{H}|_{1} \lesssim _a   |\mathcal{C}|_1 \lesssim _a |d|_3, \label{201811141623xx}
\end{align}
where we have used \eqref{esmmdforinfty}, \eqref{201803121601xx2108}, \eqref{201809291022}, \eqref{2020307saf302204} and  \eqref{20180928212154}  to derive the first inequality by following the argument of \eqref{201811210953} and \eqref{201811210954}, and  we employ \eqref{201811231df301} and \eqref{201811sdf231301}  to infer the second inequality.

\vspace{4mm} \noindent\textbf{Acknowledgements.}
The research of Fei Jiang was supported by NSFC (Grant Nos.  12022102 and 12231016) and the Natural Science Foundation of Fujian
Province of China (Grant Nos. 2020J02013 and 2022J01105), and the research of Song Jiang by National Key R\&D Program (2020YFA0712200),
National Key Project (GJXM92579), the Sino-German Science Center (Grant No. GZ 1465) and the ISF--NSFC joint research program
(Grant No. 11761141008).

\renewcommand\refname{References}
\renewenvironment{thebibliography}[1]{%
\section*{\refname}
\list{{\arabic{enumi}}}{\def\makelabel##1{\hss{##1}}\topsep=0mm
\parsep=0mm
\partopsep=0mm\itemsep=0mm
\labelsep=1ex\itemindent=0mm
\settowidth\labelwidth{\small[#1]}%
\leftmargin\labelwidth \advance\leftmargin\labelsep
\advance\leftmargin -\itemindent
\usecounter{enumi}}\small
\def\newblock{\ }
\sloppy\clubpenalty4000\widowpenalty4000
\sfcode`\.=1000\relax}{\endlist}

\end{document}